\documentclass[reqno]{amsart}

\usepackage[usenames,dvipsnames]{pstricks}
\usepackage{epsfig}

\usepackage{color}
\date{\today}

\usepackage[latin1]{inputenc}
\usepackage[T1]{fontenc}
\usepackage{amsmath, amsthm}
\usepackage[english]{babel}
\usepackage{amssymb}
\usepackage{latexsym}
\usepackage{graphicx}
\usepackage[all]{xy}

\newtheorem{apart}{}[section]
\newtheorem{teo}[apart]{{\bf Theorem}}
\newtheorem{prop}[apart]{{\bf Proposition}}
\newtheorem{lem}[apart]{{\bf Lemma}}
\newtheorem{ej}[apart]{{\bf Example}}

\newtheorem{rem}[apart]{{\bf Remark}}

\newtheorem{defin}[apart]{{\bf Definition}}

\newtheorem{cor}[apart]{{\bf Corollary}}

\newenvironment{dem}
{{\noindent{\bf Proof: }}\newline }
{$\hfill\Box$\vspace{0.25cm}}

\newcommand{\ot}{\otimes}
\newcommand{\co}{\circ}

\begin{document}

\begin{center}

{\huge{\bf Integrals and  crossed products over weak Hopf algebras}}

\end{center}

\ \\
\begin{center}
{\bf J.N. Alonso \'Alvarez$^{1}$, J.M. Fern\'andez Vilaboa$^{2}$, R.
Gonz\'{a}lez Rodr\'{\i}guez$^{3}$}
\end{center}

\ \\
\hspace{-0,5cm}$^{1}$ Departamento de Matem\'{a}ticas, Universidad
de Vigo, Campus Universitario Lagoas-Marcosende, E-36280 Vigo, Spain
(e-mail: jnalonso@ uvigo.es)
\ \\
\hspace{-0,5cm}$^{2}$ Departamento de \'Alxebra, Universidad de
Santiago de Compostela.  E-15771 Santiago de Compostela, Spain
(e-mail: josemanuel.fernandez@usc.es)
\ \\
\hspace{-0,5cm}$^{3}$ Departamento de Matem\'{a}tica Aplicada II,
Universidad de Vigo, Campus Universitario Lagoas-Marcosende, E-36310
Vigo, Spain (e-mail: rgon@dma.uvigo.es)
\ \\

{\bf Abstract} In this paper we present the general theory of cleft extensions for a cocommutative
weak Hopf algebra $H$. For a weak left $H$-module algebra we obtain a bijective correspondence
between the isomorphisms classes of $H$-cleft extensions $A_{H}\hookrightarrow A$, where $A_{H}$
is the subalgebra of coinvariants, and the equivalence classes of crossed systems for $H$ over $A_H$.
Finally, we establish a bijection between the set of equivalence classes  of crossed systems with
a fixed weak $H$-module algebra structure and the second cohomology group
$H_{\varphi_{Z(A_H)}}^2(H, Z(A_H))$, where   $Z(A_H)$ is the center of  $A_H$.

\vspace{0.5cm}

{\bf Keywords.} Monoidal category, weak Hopf algebra, cleft extension, weak crossed product, Sweedler cohomology for weak Hopf algebras.

{\bf MSC 2010:} 18D10, 16T05.

\section{Introduction}  Weak Hopf algebras (or quantum groupoids in the
terminology of Nikshych and Vainerman \cite{NV}) have been
introduced  by B\"{o}hm, Nill and Szlach\'anyi \cite{bohm1} as a new
generalization of Hopf algebras and groupoid algebras. Roughly
speaking, a weak Hopf algebra $H$ in a symmetric monoidal category
is an object which has both algebra and coalgebra structures with
some relations between them and that possesses an antipode
$\lambda_{H}$ which does not necessarily satisfy the usual convolution equalities with the identity morphism.  The main differences with other Hopf algebraic
constructions, such as quasi-Hopf algebras and rational Hopf
algebras, are the following: weak Hopf algebras are coassociative
but the coproduct is not required to preserve the unity morphism or, equivalently, the counity is not an algebra morphism. Some
motivations to study weak Hopf algebras come from their connection
with the theory of algebra extensions, the important applications in
the study of dynamical twists of Hopf algebras and their link with
quantum field theories and operator algebras (see \cite{NV}). It is well-known that groupoid algebras of finite groupoids provides examples of weak Hopf algebras. If $G$ is a finite groupoid (a category with a finite number of objects such that each morphism is invertible) then the groupoid algebra over a commutative ring $R$ is an example of cocommutative weak Hopf algebra. This weak Hopf algebra, denoted by $RG$, is generated by the morphisms $g$ in $G$ and the product of two morphisms is defined by the composition if it exists and $0$ in other case.  The coalgebra structure is defined by the coproduct $\delta_{RG}(g)=g\ot g$, and the counit $\varepsilon_{RG}(g)=1$ and the antipode by $\lambda_{RG}(g)=g^{-1}$.  There are more interesting examples of cocommutative weak Hopf algebras, for example recently Bulacu in \cite{Bul1, Bul2} proved that
Cayley-Dickson and Clifford algebras are commutative and cocommutative weak Hopf algebras  in some suitable symmetric monoidal categories of graded vector spaces.

Like in the Hopf algebra setting,  it is possible to define a theory of crossed products for weak Hopf algebras. The key  to extend the crossed product constructions of the Hopf world  to the weak setting is the  use of
idempotent morphisms combined with the ideas giving in \cite{tb-crpr}. In
\cite{nmra4, mra-preunit},  the authors defined a product on $A\otimes V$, for an
algebra $A$ and an object $V$ both living in a strict monoidal
category $\mathcal C$ where every idempotent splits. In order to obtain that
product we must consider two morphisms $\psi_{V}^{A}:V\otimes
A\rightarrow A\otimes V$ and $\sigma_{V}^{A}:V\otimes V\rightarrow
A\otimes V$ that satisfy some  twisted-like and cocycle-like
conditions. Associated to these morphisms it is possible to define
an idempotent morphism $\nabla_{A\otimes V}:A\otimes V\rightarrow
A\otimes V$, that becomes the identity in the classical case. The
image of this idempotent inherits the associative product from
$A\otimes V$. In order to define a unit for $Im (\nabla_{A\otimes
V})$, and hence to obtain an algebra structure, we require the existence of a
preunit $\nu:K\rightarrow A\otimes V$ and, under these conditions, it is possible to characterize
weak crossed products with an unit as products on $A\otimes V$
that are morphisms of left $A$-modules with preunit.
Finally, it is convenient to observe that, if the preunit is an unit,
the idempotent becomes the identity and we recover the classical
examples of the Hopf algebra setting. The theory presented in
\cite{nmra4, mra-preunit} contains as a particular instance the one
developed by Brzezi\'nski in \cite{tb-crpr} as well as all the crossed product constructed in the weak setting, for example the ones defined in \cite{Caen}, \cite{ana1} and \cite{LS}. Recently, G. B\"ohm showed
in \cite{bohm} that a monad in the weak version of the Lack and
Street's 2-category of monads in a 2-category is identical to a
crossed product system in the sense of \cite{nmra4} and also in
\cite{mra-partial-unif} we can find that  unified crossed products
\cite{AM1} and partial crossed products \cite{partial} are
particular instances of weak crossed products. An interesting example of weak crossed product comes from the theory of weak cleft extensions associated to weak Hopf algebras. This notion was introduced in \cite{nmra1} and in \cite{nmra4} we show that it provides an example of weak crossed product satisfying the  weak twisted and the cocycle conditions. These crossed products are deeply connected with Galois theory as we can see in the intrinsic characterization of weak cleftness in terms of weak Galois extensions with normal basis obtained in \cite{nmra2}.  We want to point  that, when we particularize this weak cleft theory to the Hopf algebra setting we obtain a more general notion  than the usual  of cleft extension (see Definition 7.2.1 of \cite{nmra2}) because in this case  the uniqueness  of the cleaving morphism is not guaranteed.

In the Hopf setting  the theory of crossed products arise as a generalization of the classical smash products and by the results obtained by Doi and Takeuchi in \cite{doi3} we know that every cleft extension $D\hookrightarrow A$ with cleaving morphism $f$ such that $f(1_H)=1_A$ induces a crossed product
$D\sharp_{\sigma}H$ where $\sigma:H\ot H\rightarrow D$ is a suitable convolution invertible morphism ( a normal 2-cocycle). Conversely, in \cite{blat-susan} we can find the reverse result, that is, if $D\sharp_{\sigma}H$ is a crossed product, the extension $D\hookrightarrow D\sharp_{\sigma}H$ is cleft. On the other hand, in \cite{Moss} Sweedler introduced the cohomology of a cocommutative Hopf algebra $H$ with coefficients in a commutative $H$-module algebra $A$. We will denote these cohomology groups as  $H_{\varphi_{A}}(H^{\bullet},A)$ where $\varphi_{A}$ is a fixed action of $H$ over $A$. In \cite{Moss} we can find an interesting interpretation of the second cohomology group $H_{\varphi_{A}}^2(H,A)$ in terms of extensions: This  group classifies the set of equivalence classes of cleft
extensions, i.e., classes of equivalent crossed products determined by a 2-cocycle. This result was extended by Doi \cite{doi1} proving that, in the non commutative case, there exists a bijection between the
isomorphism classes of $H$-cleft extensions $D$ of $A$ and equivalence classes of crossed systems for $H$ over $A$ with a fixed action. If $H$ is cocommutative the equivalence is described by $H_{\varphi_{\mathcal{Z}(A)}}^{2}(H,\mathcal{Z}(A))$ where $\mathcal{Z}(A)$ is the center of $A$.

The aim of this paper is to extend the preceding results to the  cocommutative weak Hopf algebra setting completing the program initiated in \cite{NikaRamon6}. To do it, in the second section, we introduce the notion of $H$-cleft extension for a  weak Hopf algebra $H$ and we prove that this kind of  extensions are examples of weak cleft extensions as the ones introduced in \cite{nmra1} and satisfying that, when we particularize  to the Hopf setting the classical notion used in the papers of Doi and Takeuchi is obtained. Also, we prove that, under cocommutative conditions, we can assume that the associated cleaving morphism is a total integral.  In the third section, assuming that $H$ is cocommutative, we prove that it is possible to identify the set of crossed systems in a weak setting as the set of weak crossed products induced by a weak left action and a convolution invertible
twisted normal 2-cocycle and then, as a consequence, we obtain the main result of this section that assures the following: If $(A,\rho_{A})$ be a right $H$-comodule algebra, there exists a bijective correspondence between the equivalence  classes of $H$-cleft extensions $A_{H}\hookrightarrow A$ and the equivalence classes of crossed systems for $H$ over $A_{H}$ where $A_{H}$ denotes the subalgebra of coinvariants in the weak setting. Finally, in the fourth section we generalize  the result obtained by Doi and Takeuchi about  the characterization of equivalence classes of crossed systems using the second Sweedler cohomology group. To obtain this generalization we must use the cohomology theory of algebras over weak Hopf algebras that we have developed in  \cite{NikaRamon6}. The main result contained in \cite{NikaRamon6} (see Theorem 3.11) asserts that if $(A,\varphi_{A})$ is a commutative left $H$-module algebra, there exists a bijection between the second cohomology group, denoted by $ H^{2}_{\varphi_{A}}(H,A)$, and the equivalence classes of  weak crossed products $A\ot_{\alpha} H$ where $\alpha:H\ot H\rightarrow A$ satisfy the 2-cocycle  and the normal conditions. Then, by this bijection and using the results of the previous sections, we obtain the description of the bijection between the isomorphism classes of $H$-cleft extensions $A_{H}\hookrightarrow B$ and the equivalence classes of crossed systems for $H$ over $A_H$ in terms of  $H_{\varphi_{\mathcal{Z}(A_{H})}}^{2}(H,\mathcal{Z}(A_{H}))$.

\section{Integrals over weak Hopf algebras}

From now on ${\mathcal C}$ denotes a strict symmetric category with
tensor product denoted by $\ot$ and unit object $K$. With $c$ we
will denote the natural isomorphism of symmetry and we also assume
that ${\mathcal C}$ has equalizers. Then, under these conditions,
every idempotent morphism $q:Y\rightarrow Y$ splits, i.e., there
exist an object $Z$ and morphisms $i:Z\rightarrow Y$ and
$p:Y\rightarrow Z$ such that $q=i\circ p$ and $p\circ i =id_{Z}$. We
denote the class of objects of ${\mathcal C}$ by $\vert {\mathcal C}
\vert $ and for each object $M\in \vert {\mathcal C}\vert$, the
identity morphism by $id_{M}:M\rightarrow M$. For simplicity of
notation, given objects $M$, $N$, $P$ in ${\mathcal C}$ and a
morphism $f:M\rightarrow N$, we write $P\ot f$ for $id_{P}\ot f$ and
$f \ot P$ for $f\ot id_{P}$.

An algebra in ${\mathcal C}$ is a triple $A=(A, \eta_{A}, \mu_{A})$
where $A$ is an object in ${\mathcal C}$ and
 $\eta_{A}:K\rightarrow A$ (unit), $\mu_{A}:A\otimes A
\rightarrow A$ (product) are morphisms in ${\mathcal C}$ such that
$\mu_{A}\circ (A\otimes \eta_{A})=id_{A}=\mu_{A}\circ
(\eta_{A}\otimes A)$, $\mu_{A}\circ (A\otimes \mu_{A})=\mu_{A}\circ
(\mu_{A}\otimes A)$. We will say that an algebra $A$ is commutative
if $\mu_{A}\co c_{A,A}=\mu_{A}$.

Given two algebras $A= (A, \eta_{A}, \mu_{A})$ and $B=(B, \eta_{B},
\mu_{B})$, $f:A\rightarrow B$ is an algebra morphism if
$\mu_{B}\circ (f\otimes f)=f\circ \mu_{A}$ and  $ f\circ \eta_{A}=
\eta_{B}.$

If $A$, $B$ are algebras in ${\mathcal C}$, the object $A\otimes B$
is an algebra in
 ${\mathcal C}$ where
$\eta_{A\otimes B}=\eta_{A}\otimes \eta_{B}$ and
\begin{equation}
\label{it-prod} \mu_{A\otimes B}=(\mu_{A}\otimes \mu_{B})\circ
(A\otimes c_{B,A}\otimes B).
\end{equation}

For an algebra $A$ we define the center of $A$ as a subalgebra
$Z(A)$  of $A$ with inclusion algebra morphism
$i_{Z(A)}:Z(A)\rightarrow A$ satisfying

$$\mu_{A}\co (A\ot i_{Z(A)})=\mu_{A}\co c_{A,A}\co (A\ot i_{Z(A)})$$
\vspace{0.05cm}

and if $f:B\rightarrow A$ is a morphism such that $\mu_{A}\co (A\ot
f)=\mu_{A}\co c_{A,A}\co (A\ot f)$, there exists an unique morphism
$f^{\prime}:B\rightarrow Z(A)$ satisfying $i_{Z(A)}\co f^{\prime}=f.$ As a consequence, we obtain that $
Z(A)$ is a commutative algebra. For example, if ${\mathcal C}$ is a
closed category and $\alpha_{A}$ and $\beta_{A}$ are the unit and the counit,
respectively, of the ${\mathcal C}$-adjunction $A\otimes -\dashv
[A,-]:{\mathcal C}\rightarrow {\mathcal C}$, the center of $A$ can
be obtained by the  following equalizer diagram:

$$
\setlength{\unitlength}{3mm}
\begin{picture}(30,4)
\put(3,2){\vector(1,0){4}} \put(11,2.5){\vector(1,0){10}}
\put(11,1.5){\vector(1,0){10}} \put(0,2){\makebox(0,0){$Z(A)$}} \put(9,2){\makebox(0,0){$A$}} \put(24,2){\makebox(0,0){$[A,
A]$}} \put(5.5,3){\makebox(0,0){$i_{Z(A)}$}}
\put(16,3.5){\makebox(0,0){$\vartheta_{A}$}}
\put(16,0.15){\makebox(0,0){$\theta_{A}$}}
\end{picture}
$$
\vspace{0.05cm}

where $\vartheta_{A}=[A, \mu_{A}]\circ \alpha_{A}(A)$ and
$\theta_{A}=[A, \mu_{A}\circ c_{A,A}]\circ \alpha_{A}(A)$.

A coalgebra in ${\mathcal C}$ is a triple ${D} = (D,
\varepsilon_{D}, \delta_{D})$ where $D$ is an object in ${\mathcal
C}$ and $\varepsilon_{D}: D\rightarrow K$ (counit),
$\delta_{D}:D\rightarrow D\otimes D$ (coproduct) are morphisms in
${\mathcal C}$ such that $(\varepsilon_{D}\otimes D)\circ
\delta_{D}= id_{D}=(D\otimes \varepsilon_{D})\circ \delta_{D}$,
$(\delta_{D}\otimes D)\circ \delta_{D}=
 (D\otimes \delta_{D})\circ \delta_{D}.$ We will say that $D$ is
 cocommutative if $c_{D,D}\co \delta_{D}=\delta_{D}$ holds.

 If ${D} = (D, \varepsilon_{D},
 \delta_{D})$ and
${ E} = (E, \varepsilon_{E}, \delta_{E})$ are coalgebras,
$f:D\rightarrow E$ is a coalgebra morphism if $(f\otimes f)\circ
\delta_{D} =\delta_{E}\circ f$ and $\varepsilon_{E}\circ f
=\varepsilon_{D}.$

When $D$, $E$ are coalgebras in ${\mathcal C}$, $D\otimes E$ is a
coalgebra in ${\mathcal C}$ where $\varepsilon_{D\otimes
E}=\varepsilon_{D}\otimes \varepsilon_{E}$ and
\begin{equation}
\label{it-coprod} \delta_{D\otimes E}=(D\otimes c_{D,E}\otimes
E)\circ( \delta_{D}\otimes \delta_{E}).
\end{equation}

If $A$ is an algebra, $B$ is a coalgebra and $\alpha:B\rightarrow
A$, $\beta:B\rightarrow A$ are morphisms, we define the convolution
product by
$$\alpha\wedge \beta=\mu_{A}\circ (\alpha\otimes
\beta)\circ \delta_{B}.$$

 Let  $A$ be an algebra. The pair
$(M,\varphi_{M})$ is a left $A$-module if $M$ is an object in
${\mathcal C}$ and $\varphi_{M}:A\otimes M\rightarrow M$ is a morphism
in ${\mathcal C}$ satisfying $\varphi_{M}\circ(\eta_A\otimes M)
=id_{M}$, $\varphi_{M}\circ (A \ot \varphi_{M})=\varphi_{M}\co
(\mu_{A}\ot M)$. Given two right ${A}$-modules $(M,\varphi_{M})$
and $(N,\varphi_{N})$, $f:M\rightarrow N$ is a morphism of right
${A}$-modules if $\varphi_{N}\circ (A\ot f)=f\circ \varphi_{M}$.

 Let  $C$ be a coalgebra. The pair
$(M,\rho_{M})$ is a right $C$-comodule if $M$ is an object in
${\mathcal C}$ and $\rho_{M}:M\rightarrow M\ot C$ is a morphism in
${\mathcal C}$ satisfying $(M\otimes \varepsilon_{C})\co
\rho_{M}=id_{M}$, $(M\ot \rho_{M})\co \rho_{M}=(M\ot \delta_{C})\co
\rho_{M}$. Given two right ${C}$-comodules $(M,\rho_{M})$ and
$(N,\rho_{N})$, $f:M\rightarrow N$ is a morphism of right
${C}$-comodules if $(f\otimes C)\co \rho_{M}=\rho_{N}\co f$.

By weak Hopf algebras  we understand the objects introduced in
\cite{bohm1}, as a generalization of ordinary Hopf algebras. Here we
recall the definition of these objects in a monoidal symmetric setting.

\begin{defin}
\label{wha} {\rm A weak Hopf algebra $H$  is an object in ${\mathcal
C}$ with an algebra structure $(H, \eta_{H},\mu_{H})$ and a
coalgebra structure $(H, \varepsilon_{H},\delta_{H})$ such that the
following axioms hold:
\begin{itemize}
\item[(a1)] $\delta_{H}\circ \mu_{H}=(\mu_{H}\otimes \mu_{H})\circ
\delta_{H\otimes H},$ \item[(a2)]$\varepsilon_{H}\circ \mu_{H}\circ
(\mu_{H}\otimes H)=(\varepsilon_{H}\otimes \varepsilon_{H})\circ
(\mu_{H}\otimes \mu_{H})\circ (H\otimes \delta_{H}\otimes H)$ \item[
]$=(\varepsilon_{H}\otimes \varepsilon_{H})\circ (\mu_{H}\otimes
\mu_{H})\circ (H\otimes (c_{H,H}\circ\delta_{H})\otimes H),$
\item[(a3)]$(\delta_{H}\otimes H)\circ \delta_{H}\circ
\eta_{H}=(H\otimes \mu_{H}\otimes H)\circ (\delta_{H}\otimes
\delta_{H})\circ (\eta_{H}\otimes \eta_{H})$ \item[ ]$=(H\otimes
(\mu_{H}\circ c_{H,H})\otimes H)\circ (\delta_{H}\otimes
\delta_{H})\circ (\eta_{H}\otimes \eta_{H}).$

\item[(a4)] There exists a morphism $\lambda_{H}:H\rightarrow H$
in ${\mathcal C}$ (called the antipode of $H$) satisfying:
\begin{itemize}
\item[(a4-1)] $id_{H}\wedge \lambda_{H}=((\varepsilon_{H}\circ
\mu_{H})\otimes H)\circ (H\otimes c_{H,H})\circ ((\delta_{H}\circ
\eta_{H})\otimes H),$ \item[(a4-2)] $\lambda_{H}\wedge
id_{H}=(H\otimes(\varepsilon_{H}\circ \mu_{H}))\circ (c_{H,H}\otimes
H)\circ (H\otimes (\delta_{H}\circ \eta_{H})),$
\item[(a4-3)]$\lambda_{H}\wedge id_{H}\wedge
\lambda_{H}=\lambda_{H}.$
\end{itemize}
\end{itemize}

Note that, in this definition, the conditions (a2), (a3) weaken the
conditions of multiplicativity of the counit, and comultiplicativity
of the unit that we can find in the Hopf algebra definition. On the
other hand, axioms (a4-1), (a4-2) and (a4-3) weaken the properties
of the antipode in a Hopf algebra. Therefore, a weak Hopf algebra is
a Hopf algebra if an only if the morphism $\delta_{H}$
(comultiplication) is unit-preserving or if and only if the counit
is a homomorphism of algebras. }
\end{defin}

\begin{apart}
{\rm If $H$ is a weak Hopf algebra in ${\mathcal C}$, the antipode
$\lambda_{H}$ is unique, antimultiplicative, anticomultiplicative
and leaves the unit $\eta_{H}$ and the counit $\varepsilon_{H}$
invariant:
\begin{equation} \lambda_{H}\circ \mu_{H}=\mu_{H}\circ
(\lambda_{H}\otimes \lambda_{H})\circ
c_{H,H};\;\;\;\;\delta_{H}\circ \lambda_{H}=c_{H,H}\circ
(\lambda_{H}\otimes \lambda_{H})\circ \delta_{H};
\end{equation}
\begin{equation}
\lambda_{H}\circ \eta_{H}=\eta_{H};\;\;\;\;\varepsilon_{H}\circ
\lambda_{H}=\varepsilon_{H}.
\end{equation}

If we define the morphisms $\Pi_{H}^{L}$ (target), $\Pi_{H}^{R}$
(source), $\overline{\Pi}_{H}^{L}$ and $\overline{\Pi}_{H}^{R}$ by
\begin{itemize}
\item[ ]$\Pi_{H}^{L}=((\varepsilon_{H}\circ \mu_{H})\otimes
H)\circ (H\otimes c_{H,H})\circ ((\delta_{H}\circ \eta_{H})\otimes
H);$ \item[ ]$\Pi_{H}^{R}=(H\otimes(\varepsilon_{H}\circ
\mu_{H}))\circ (c_{H,H}\otimes H)\circ (H\otimes (\delta_{H}\circ
\eta_{H}));$ \item[ ]$\overline{\Pi}_{H}^{L}=(H\otimes
(\varepsilon_{H}\circ \mu_{H}))\circ ((\delta_{H}\circ
\eta_{H})\otimes H);$ \item[
]$\overline{\Pi}_{H}^{R}=((\varepsilon_{H}\circ \mu_{H})\otimes
H)\circ(H\otimes (\delta_{H}\circ \eta_{H})).$
\end{itemize}
it is straightforward to show (see \cite{bohm1}) that they are
idempotent and $\Pi_{H}^{L}$, $\Pi_{H}^{R}$ satisfy the equalities
\begin{equation}
\label{propiedadesgeneralespiLR}
\Pi_{H}^{L}=id_{H}\wedge
\lambda_{H};\;\;\;\;\;\Pi_{H}^{R}=\lambda_{H}\wedge id_{H}.
\end{equation}
and then
\begin{equation}
\label{id} \Pi_{H}^{L}\wedge \Pi_{H}^{L}=\Pi_{H}^{L},\;\;\;
\Pi_{H}^{R}\wedge \Pi_{H}^{R}=\Pi_{H}^{R}.
\end{equation}

 Moreover, we have that
\begin{equation}
\label{composiciones1}
\Pi_{H}^{L}\circ \overline{\Pi}_{H}^{L}=\Pi_{H}^{L};\;\;\;\;
\Pi_{H}^{L}\circ
\overline{\Pi}_{H}^{R}=\overline{\Pi}_{H}^{R};\;\;\;\;
\Pi_{H}^{R}\circ
\overline{\Pi}_{H}^{L}=\overline{\Pi}_{H}^{L};\;\;\;\;
\Pi_{H}^{R}\circ \overline{\Pi}_{H}^{R}=\Pi_{H}^{R};
\end{equation}
\begin{equation}
\label{composiciones2}
\overline{\Pi}_{H}^{L}\circ
\Pi_{H}^{L}=\overline{\Pi}_{H}^{L};\;\;\;\;
\overline{\Pi}_{H}^{L}\circ \Pi_{H}^{R}=\Pi_{H}^{R};\;\;\;\;
\overline{\Pi}_{H}^{R}\circ \Pi_{H}^{L}=\Pi_{H}^{L};\;\;\;\;
\overline{\Pi}_{H}^{R}\circ \Pi_{H}^{R}=\overline{\Pi}_{H}^{R}.
\end{equation}

Also it is easy  to show the formulas
\begin{equation}\label{composiciones}
\Pi_{H}^{L}=\overline{\Pi}_{H}^{R}\circ \lambda_{H}=\lambda_{H}
\circ\overline{\Pi}_{H}^{L};\;\;\;\;\Pi_{H}^{R}=
\overline{\Pi}_{H}^{L}\circ \lambda_{H}=\lambda_{H} \circ
\overline{\Pi}_{H}^{R};
\end{equation}
\begin{equation}
\Pi_{H}^{L}\circ \lambda_{H}=\Pi_{H}^{L}\circ \Pi_{H}^{R}=
\lambda_{H}\circ \Pi_{H}^{R};\;\;\;\;\Pi_{H}^{R}\circ
\lambda_{H}=\Pi_{H}^{R}\circ \Pi_{H}^{L}= \lambda_{H}\circ
\Pi_{H}^{L}.
\end{equation}

For the morphisms target an source we have the following identities:

\begin{equation}
\label{PiLmu} \mu_H\co (H\ot \Pi^{L}_{H})=((\varepsilon_H\co\mu_H)\ot H)\co (H\ot c_{H,H})\co (\delta_H\ot H),
\end{equation}

\begin{equation}
\label{PiRmu} \mu_H\co (\Pi^{R}_{H}\ot H)=(H\ot (\varepsilon_H\co\mu_H))\co (c_{H,H}\ot H)\co (H\ot \delta_H),
\end{equation}

\begin{equation}
\label{PiLbarramu} \mu_H\co (H\ot \overline{\Pi}^{L}_{H})=(H\ot (\varepsilon_H\co\mu_H))\co (\delta_H\ot H),
\end{equation}

\begin{equation}
\label{PiRbarramu} \mu_H\co (\overline{\Pi}^{R}_{H}\ot H)=((\varepsilon_H\co\mu_H)\ot H)\co (H\ot \delta_H),
\end{equation}

\begin{equation}
\label{deltaPIL}
(H\ot \Pi^{L}_{H})\co \delta_H=(\mu_H\ot H)\co (H\ot c_{H,H})\co ((\delta_H\co\eta_H)\ot H),
\end{equation}

\begin{equation}
\label{deltaPIR}
(\Pi^{R}_{H}\ot H)\co \delta_H=(H\ot \mu_H)\co (c_{H,H}\ot H)\co (H\ot (\delta_H\co\eta_H)),
\end{equation}

\begin{equation}
\label{deltaPILbarra}
(\overline{\Pi}^{L}_{H}\ot H)\co \delta_H=(H\ot \mu_H)\co ((\delta_H\co\eta_H)\ot H),
\end{equation}

\begin{equation}
\label{deltaPIRbarra}
(H\ot \overline{\Pi}^{R}_{H})\co \delta_H=(\mu_H\ot H)\co (H\ot (\delta_H\co\eta_H)),
\end{equation}

and

\begin{equation}
\label{PiLmuPiL}
\Pi^{L}_{H}\circ \mu_{H}\circ (H\ot
\Pi^{L}_{H})=\Pi^{L}_{H}\circ \mu_{H},
\end{equation}

\begin{equation}
\label{PiRmuPiR}
\Pi^{R}_{H}\circ
\mu_{H}\circ (\Pi^{R}_{H}\ot H)=\Pi^{R}_{H}\circ \mu_{H},
\end{equation}

\begin{equation}
\label{PiLdeltaPiL}
(H\ot \Pi^{L}_{H})\circ \delta_{H}\circ
\Pi^{L}_{H}=\delta_{H}\circ \Pi^{L}_{H},
\end{equation}

\begin{equation}
\label{PiRdeltaPiR}
( \Pi^{R}_{H}\ot H)\circ \delta_{H}\circ \Pi^{R}_{H}=\delta_{H}\circ \Pi^{R}_{H}.
\end{equation}
}
\end{apart}

\begin{defin}
\label{H-comodalg} {\rm Let $H$ be a weak Hopf algebra. We will say
that a right $H$-comodule $(A, \rho_{A})$ is a right $H$-comodule algebra if it satisfies
$$\rho_A\co\mu_A=\mu_{A\ot H}\co (\rho_{A}\ot \rho_{A})$$

and any of the following equivalent conditions hold:

\begin{itemize}

\item[(b1)] $(A\ot\Pi^{L}_{H})\co \rho_A=(\mu_A\ot H)\co (A\ot c_{H,A})\co((\rho_A\co\eta_A)\ot A).$

\item[(b2)] $(A\ot\overline{\Pi}^{R}_{H})\co \rho_A=(\mu_A\ot H)\co(A\ot (\rho_A\co\eta_A)).$

\item[(b3)] $(A\ot\Pi^{L}_{H})\co \rho_A\co\eta_A=\rho_A\co\eta_A.$

\item[(b4)] $(A\ot\overline{\Pi}^{R}_{H})\co \rho_A\co\eta_A=\rho_A\co\eta_A.$

\item[(b5)] $(\rho_A\ot H)\co \rho_A\co\eta_A=(A\ot\mu_H\ot H)\co (\rho_A\ot \delta_H)\co (\eta_A\ot\eta_H).$

\item[(b6)] $(\rho_A\ot H)\co \rho_A \co \eta_A=(A\ot (\mu_H\co c_{H,H})\ot H)\co (\rho_A\ot \delta_H)\co (\eta_A\ot\eta_H).$

\end{itemize}

If $(A,\rho_{A})$ is a right $H$-comodule algebra, the triple $(A,H,\Gamma_{A}^{H})$ is a right-right
weak entwining structure (see \cite{Caen}) where
\begin{equation}
\label{weak-Gamma} \Gamma_{A}^{H} = (A\ot \mu_{H})\co (c_{H,A}\ot H)\co
(H\ot \rho_{A}).
\end{equation}
Therefore the following identies hold:
\begin{equation}
\label{a)} \Gamma_{A}^{H} \circ (H\otimes \mu_A ) = (\mu_A \otimes
H)\circ (A\otimes \Gamma_{A}^{H} )\circ (\Gamma_{A}^{H} \otimes A),
\end{equation}
\begin{equation}
\label{b)} (A\otimes \delta_H)\circ \Gamma_{A}^{H} = (\Gamma_{A}^{H} \otimes
H)\circ (H\otimes \Gamma_{A}^{H} )\circ (\delta_H\otimes A),
\end{equation}
\begin{equation}
\label{c)} \Gamma_{A}^{H} \circ (H\otimes \eta_A) = (e_{A}\otimes
H)\circ \delta_H,
\end{equation}
\begin{equation}
\label{d)} (A\otimes \varepsilon_H)\circ \Gamma_{A}^{H} = \mu_A\circ
(e_{A}\otimes A),
\end{equation}
where
\begin{equation}
\label{e-gamma}
 e_{A} = (A\otimes \varepsilon_H)\circ \Gamma_{A}^{H} \circ
(H\otimes \eta_A).
\end{equation}
 We denote by ${\mathcal
M}_A^H(\Gamma_{A}^{H})$ the category of weak entwined modules, i.e., the
objects $M$ in ${\mathcal C}$ together with two morphisms
$\phi_M:M\otimes A\rightarrow A$ and $\rho_M:M\rightarrow M\otimes
H$ such that $(M, \phi_M)$ is a right $A$-module, $(M, \rho_M)$ is a
right $H$-comodule and such that the following equality
\begin{equation}
\label{m)} \rho_M\circ \phi_M = (\phi_M\otimes H)\circ (M\otimes
\Gamma_{A}^{H} ) \circ (\rho_M\otimes A)
\end{equation}
holds.

Then, if $(A,\rho_{A})$ is a right $H$-comodule algebra, $(A, \mu_{A},\rho_{A})$ is an object of ${\mathcal M}_A^H(\Gamma_{A}^{H})$.

 If $(A, \rho_A)$ is a right $H$-comodule algebra, we define the subalgebra of coinvariants of $A$ as the equalizer:

$$
\setlength{\unitlength}{3mm}
\begin{picture}(30,4)
\put(3,2){\vector(1,0){4}} \put(11,2.5){\vector(1,0){10}}
\put(11,1.5){\vector(1,0){10}} \put(1,2){\makebox(0,0){$A_{H}$}}
\put(9,2){\makebox(0,0){$A$}} \put(24,2){\makebox(0,0){$A\ot H$}}
\put(5.5,3){\makebox(0,0){${i_A}$}} \put(16,3.5){\makebox(0,0){$
\rho_A$}} \put(16,0.5){\makebox(0,0){$\zeta_{A}$}}
\end{picture}
$$
where $\zeta_{A}=(\mu_A\ot H)\co (A\ot
c_{H,A})\co((\rho_A\co\eta_A)\ot A)$. Note that, by (b1), we have
\begin{equation}
\label{zetaA} \zeta_{A}=(A\ot\Pi^{L}_{H})\co \rho_A.
\end{equation}
and also, by (\ref{composiciones1}) and (\ref{composiciones2}),
$$
\setlength{\unitlength}{3mm}
\begin{picture}(30,4)
\put(3,2){\vector(1,0){4}} \put(11,2.5){\vector(1,0){10}}
\put(11,1.5){\vector(1,0){10}} \put(1,2){\makebox(0,0){$A_{H}$}}
\put(9,2){\makebox(0,0){$A$}} \put(24,2){\makebox(0,0){$A\ot H$}}
\put(5.5,3){\makebox(0,0){${i_A}$}} \put(16,3.5){\makebox(0,0){$
\rho_A$}}
\put(16,0.5){\makebox(0,0){$(A\ot\overline{\Pi}^{R}_{H})\co
\rho_A$}}
\end{picture}
$$
is an equalizer diagram.

It is not difficult to see that $(A_H, \eta_{A_H}, \mu_{A_H})$ is an algebra, being $\eta_{A_H}$ and $\mu_{A_H}$ the factorizations through the equalizer $i_A$ of the morphisms $\eta_A$ and
 $\mu_A\co(i_A\ot i_A)$, respectively.

 For example, the weak Hopf algebra $H$ is a right $H$-comodule algebra with right comodule structure giving by $\rho_H=\delta_H$ and subalgebra of coinvariants $H_H$, the image of the idempotent morphism $\Pi^{L}_{H}$, which we will denote by $H_{L}$.

}
\end{defin}

\begin{defin}
\label{integral}
{\rm  Let $H$ be a weak Hopf algebra and let $(A, \rho_A)$ be a right $H$-comodule algebra. We define an integral as a morphism of right $H$-comodules $f:H\rightarrow A$. If moreover $f\co \eta_H=\eta_A$ we will say that the integral is total.

An integral $f:H\rightarrow A$ is convolution invertible if there exists a morphism $f^{-1}:H\rightarrow A$ (called the convolution inverse of $f$) such that

\begin{itemize}

\item[(c1)] $f^{-1}\wedge f=e_{A}.$

\item[(c2)] $f\wedge f^{-1}=(A\ot(\varepsilon_H\co\mu_H))\co ((\rho_A\co \eta_A)\ot H).$

\item[(c3)] $f^{-1}\wedge f\wedge f^{-1}=f^{-1}.$

\end{itemize}

where $e_{A}$ is the morphism defined in (\ref{e-gamma}).

Trivially, the inverse is unique because if $h:H\rightarrow A$ satisfies (c1)-(c3),
\begin{equation}\label{igualdadinversa}
h=h\wedge f\wedge h=h\wedge f\wedge f^{-1}=f^{-1}\wedge f\wedge f^{-1}=f^{-1}.
\end{equation}

Moreover, using condition (c1), if $f$ is an integral convolution invertible we get that $f\wedge f^{-1}\wedge f=f.$

Finally, when $f$ is a total integral we can rewrite equality (c1) as $f^{-1}\wedge f=f\co \Pi^{R}_{H}$ and (c2) as $f\wedge f^{-1}=f\co \overline{\Pi}^{L}_{H}.$
}
\end{defin}

\begin{ej}
\label{antipodoesintegral} {\rm
Let $H$ be a weak Hopf algebra such that $\Pi^{L}_{H}=\overline{\Pi}^{L}_{H}$ (equivalently, $\Pi^{R}_{H}=\overline{\Pi}^{R}_{H}$). Then the identity $id_H$ is a total integral convolution invertible with inverse $\lambda_H$. Note that this equality is always true in the Hopf algebra setting. In our case it holds, for example, if $H$ is a cocommutative weak Hopf algebra.
}
\end{ej}

\begin{defin}
\label{Cleft}
{\rm Let $H$ be a weak Hopf algebra and $(A,\rho_{A})$ a right $H$-comodule algebra.
We say that $A_{H}\hookrightarrow A$ is a $H$-cleft extension if
there exists an integral $f:H\rightarrow A$ convolution invertible and such
that the morphism $f\wedge f^{-1}$ factorizes through the equalizer $i_A$. In what follows,
the morphism $f$ will be called  a cleaving morphism associated to the $H$-cleft
extension $A_{H}\hookrightarrow A$.
}

\end{defin}

\begin{prop}
\label{propinversa}
Let $H$ be a weak Hopf algebra and $(A,\rho_{A})$ a right $H$-comodule algebra such that
$A_{H}\hookrightarrow A$ is a $H$-cleft extension with cleaving morphism $f$. Then the equality
\begin{equation}
\label{cleaving-lambda} \rho_A\co f^{-1}=(f^{-1}\ot \lambda_H)\co
c_{H,H}\co \delta_H
\end{equation}
holds.
\end{prop}

\begin{dem}
We define the  morphisms: $r=\rho_A\co f^{-1}$, $s=\rho_A\co f$ and $t=(f^{-1}\ot \lambda_H)\co c_{H,H}\co \delta_H$.

First of all, we show that $s\wedge r=s\wedge t$. Indeed:

\begin{itemize}

\item[ ]$\hspace{0.38cm} s\wedge r$

\item[ ]

\item[ ]$=\rho_A\co (f\wedge f^{-1})$

\item[ ]

\item[ ]$=(A\ot \Pi^{L}_{H})\co \rho_A\co (f\wedge f^{-1})$

\item[ ]

\item[ ]$=(\mu_A\ot (\Pi^{L}_{H}\co \mu_H\co (H\ot \Pi^{L}_{H})))\co (A\ot c_{H,A}\ot H)\co ((\rho_A\co f)\ot (\rho_A\co f^{-1}))\co \delta_H$

\item[ ]

\item[ ]$=(\mu_A\ot \Pi^{L}_{H})\co (A\ot c_{H,A})\co ((\rho_A\co f)\ot f^{-1})\co \delta_H$

\item[ ]

\item[ ]$=(\mu_A\ot H)\co (A\ot c_{H,A})\co (((f\ot (\mu_H\co (H\ot \lambda_H)\co \delta_H))\co \delta_{H})\ot f^{-1})\co \delta_H$

\item[ ]

\item[ ]$=s\wedge t.$

\end{itemize}

In the foregoing calculations, the first equality follows using that $A$ is a right $H$-comodule algebra; the second one because $A_{H}\hookrightarrow A$ is  $H$-cleft; in the third we use (\ref{PiLmuPiL}); the fourth relies on the equality
\begin{equation}
\label{rhoPiL} (\mu_A\ot \mu_H)\co (A\ot c_{H,A}\ot \Pi^{L}_{H})\co (\rho_A\ot \rho_A)=(\mu_A\ot H)\co (A\ot c_{H,A})\ot (\rho_A\ot A),
\end{equation}

the fifth one is a consequence of the definition of $\Pi^{L}_{H}$; finally, in the last one we use that $f$ is an integral.

Using similar techniques we obtain that $t\wedge s=r\wedge s$:

\begin{itemize}

\item[ ]

\item[ ]$\hspace{0.38cm} t\wedge s$

\item[ ]

\item[ ]$=(A\ot \mu_H)\co (c_{H,A}\ot H)\co (\lambda_H\ot (f^{-1}\wedge f)\ot H)\co (H\ot \delta_H)\co \delta_H$

\item[ ]

\item[ ]$=(A\ot \mu_H)\co (A\ot (\lambda_{H}\ot (\varepsilon_{H}\co \mu_{H}))
\co (\delta_{H}\ot H))\ot H)\co (c_{H,A}\ot H\ot H)\co (H\ot (\rho_A\co \eta_A)\ot H)\co \delta_H$

\item[ ]

\item[ ]$=(A\ot \mu_H)\co(A\ot (\lambda_{H}\co \mu_{H}\co (H\ot \overline{\Pi}_{H}^{L}))\ot H)\co (c_{H,A}\ot H\ot H)\co (H\ot (\rho_A\co \eta_A)\ot H)\co \delta_H$

\item[ ]

\item[ ]$=(A\ot \mu_H)\co (A\ot (\mu_{H}\co c_{H,H}\co (\lambda_{H}\ot \Pi_{H}^{L}))\ot H)\co (c_{H,A}\ot H\ot H)\co (H\ot (\rho_A\co \eta_A)\ot H)\co \delta_H$

\item[ ]

\item[ ]$=(A\ot \mu_H)\co ((\rho_A\co \eta_A)\co \Pi^{R}_{H})$

\item[ ]

\item[ ]$=(A\ot H\ot (\varepsilon_H\co \mu_H))\co (A\ot c_{H,H}\ot H)\co (c_{H,A}\ot \mu_H\ot H)\co (H\ot (\rho_A\co \eta_A)\ot (\delta_H\co \eta_H))$

\item[ ]

\item[ ]$=(A\ot H\ot (\varepsilon_H\co \mu_H))\co (A\ot c_{H,H}\ot H)\co (c_{H,A}\ot H\ot H)\co (H\ot ((\rho_A\ot H)\co (\rho_{A}\co \eta_A)))$

\item[ ]

\item[ ]$=\rho_A\co (f^{-1}\wedge f)$

\item[ ]

\item[ ]$=r\wedge s.$

\end{itemize}

In the previous calculations, the first equality follows because $f$ is an integral; the second one by coassociativity of $\delta_{H}$ and the naturality of $c$; in the third we use (\ref{PiLbarramu}). The fourth equality relies on the antimultiplicativity of the antipode and (\ref{composiciones}); the fifth and the sixth ones follow by the definition of $\Pi^{R}_{H}$; finally, in the seventh,  eight and ninth equalities we use that $A$ is a right $H$-comodule algebra.

On the other hand, under the conditions of this proposition we have that the equality
\begin{equation}
\label{f-1fmu}
(f^{-1}\wedge f)\co \mu_{H}=  ((\varepsilon_{H}\co \mu_{H})\ot (f^{-1}\wedge f))\co (H\ot \delta_{H})
\end{equation}
holds. Indeed:
\begin{itemize}

\item[ ]

\item[ ]$\hspace{0.38cm} (f^{-1}\wedge f)\co \mu_{H}$

\item[ ]

\item[ ]$=(A\ot (\varepsilon_{H}\co \mu_{H}\co (H\ot \mu_{H})))\co (c_{H,A}\ot H\ot H)\co
(H\ot c_{H,A}\ot H)\co (H\ot H\ot (\rho_{A}\co \eta_{A})) $

\item[ ]

\item[ ]$=(A\ot ((((\varepsilon_{H}\co \mu_{H})\ot (\varepsilon_{H}\co \mu_{H}))\co (H\ot \delta_{H}\ot H))))\co (c_{H,A}\ot H\ot H)\co
(H\ot c_{H,A}\ot H)\co (H\ot H\ot (\rho_{A}\co \eta_{A})) $

\item[ ]

\item[ ]$=(A\ot (\varepsilon_{H}\co \mu_{H}))\co ((\varepsilon_{H}\co \mu_{H})\ot c_{H,A}\ot H)\co
(H\ot \delta_{H}\ot \rho_{A}\co \eta_{A}))$

\item[ ]

\item[ ]$= ((\varepsilon_{H}\co \mu_{H})\ot (f^{-1}\wedge f))\co (H\ot \delta_{H}), $

\item[ ]

\end{itemize}

where the first and the last equalities follow because $f$ is an  integral, the second one by (a2) of Definition \ref{wha}
and the third one by the naturality of $c$.

Now we use the coassociativity of $\delta_{H}$, the naturality of $c$,
(\ref{f-1fmu}) and that $f$ is a convolution invertible integral to get
that $t\wedge s\wedge t=t$.

\begin{itemize}

\item[ ]

\item[ ]$\hspace{0.38cm} t\wedge s\wedge t$

\item[ ]

\item[ ]$=(\mu_A\ot H)\co (A\ot c_{H,A})\co (A\ot \mu_H\ot A)\co
(c_{H,A}\ot H\ot A)\co (\lambda_H\ot (((f^{-1}\wedge f)\ot
\Pi^{L}_{H})\co \delta_{H})\ot A)\co $

\item[ ]

\item[ ]$\hspace{0.38cm}
(\delta_H\ot f^{-1})\co \delta_H$

\item[ ]

\item[ ]$=(\mu_A\ot H)\co (A\ot c_{H,A})\co (A\ot \mu_H\ot A)\co
(c_{H,A}\ot H\ot A)\co $

\item[ ]

\item[ ]$\hspace{0.38cm} (\lambda_H\ot ((((f^{-1}\wedge f)\co
\mu_{H})\ot H)\co (H\ot c_{H,H})\co ((\delta_{H}\co \eta_{H})\ot
H))\ot A)\co (\delta_H\ot f^{-1})\co \delta_H$

\item[ ]

\item[ ]$=(\mu_A\ot H)\co (A\ot c_{H,A})\co (A\ot \mu_H\ot A)\co
(c_{H,A}\ot H\ot A)\co $

\item[ ]

\item[ ]$\hspace{0.38cm}  (\lambda_H\ot (((((\varepsilon_{H}\co
\mu_{H})\ot (f^{-1}\wedge f))\co (H\ot \delta_{H}))\ot H)\co (H\ot
c_{H,H})\co ((\delta_{H}\co \eta_{H})\ot H))\ot A)\co (\delta_H\ot
f^{-1})\co \delta_H$

\item[ ]

\item[ ]$=c_{H,A}\co ((\lambda_H\wedge \Pi^{L}_{H})\ot (f^{-1}\wedge f\wedge f^{-1}))\co \delta_H$

\item[ ]

\item[ ]$=t.$

\end{itemize}

Taking into account that $r\wedge s\wedge r=r$,
$$r=r\wedge s\wedge r=t\wedge s\wedge r=t\wedge s\wedge t=t$$
and we conclude the proof.

\end{dem}

\begin{prop}
\label{siempreCleft} Let $H$ be a cocommutative weak Hopf algebra
and let $(A, \rho_A)$ be a right $H$-comodule algebra. If there
exists   a convolution invertible integral $f:H\rightarrow A$, then
$A_{H}\hookrightarrow A$ is an $H$-cleft extension.

\end{prop}

\begin{dem} Let $f^{-1}$ be the convolution inverse of $f$. We have to show that $f\wedge f^{-1}$
factorizes through the equalizer $i_A$. Indeed:

\begin{itemize}

\item[ ]

\item[ ]$\hspace{0.38cm} \zeta_{A} \co (f\wedge f^{-1})$

\item[ ]

\item[ ]$=(A\ot \Pi^{L}_{H})\co \rho_A\co (A\ot (\varepsilon_H\co \mu_H))\co ((\rho_A\co \eta_A)\ot H)$

\item[ ]

\item[ ]$=(A\ot H\ot (\varepsilon_H\co \mu_H))\co (A\ot ((\Pi^{L}_{H}\ot H)\co\delta_H)\ot H)\co ((\rho_A\co \eta_A)\ot H)$

\item[ ]

\item[ ]$=(A\ot H\ot (\varepsilon_H\co \mu_H))\co (A\ot ((\Pi^{L}_{H}\ot H)\co c_{H,H}\co\delta_H)\ot H)\co ((\rho_A\co \eta_A)\ot H)$

\item[ ]

\item[ ]$=(A\ot (\varepsilon_H\co \mu_H)\ot H)\co (\rho_A\ot c_{H,H})\co  (((A\ot \Pi^{L}_{H})\co \rho_A\co
\eta_A)\ot H)$

\item[ ]

\item[ ]$=(A\ot (\varepsilon_H\co \mu_H)\ot H)\co (\rho_A\ot c_{H,H})\co ((\rho_A\co \eta_A)\ot H)$

\item[ ]

\item[ ]$=(A\ot H\ot (\varepsilon_H\co \mu_H))\co (A\ot (c_{H,H}\co \delta_H)\ot H)\co ((\rho_A\co \eta_A)\ot H)$

\item[ ]

\item[ ]$=(A\ot H\ot (\varepsilon_H\co \mu_H))\co (A\ot \delta_H\ot H)\co ((\rho_A\co \eta_A)\ot H)$

\item[ ]

\item[ ]$=\rho_A\co (f\wedge f^{-1}).$

\end{itemize}

In the foregoing calculations, the first and the last equalities
follow by (c2) of Definition \ref{integral}; the second, fourth and
 sixth ones use the condition of comodule for $A$; in the third
and seventh we use that $H$ is cocommutative; finally the fifth one
follows by (b3) of Definition \ref{H-comodalg}. $\square$

\end{dem}

\begin{apart}\label{relacioncleftdefil}
{\rm Let $H$ be a  weak Hopf algebra and let $(A, \rho_A)$ be a
right $H$-comodule algebra. We want to point out the relation
between the notion  of $H$-cleft extension and the notion of weak
$H$-cleft extension given by us in \cite{nmra1}. In \cite{nmra1} we
introduce the set $Reg^{WR}(H,A)$ as the one whose elements are the
morphisms $h:H\rightarrow A$ such that there exists a morphism
$h^{-1}:H\rightarrow A$, called the left weak inverse of $h$, such
that
\begin{equation}
\label{left-weak-inverse} h^{-1}\wedge h=e_{A}
\end{equation}
where $e_{A}$ is the morphism defined in (\ref{e-gamma}) for the
right-right weak entwining structure $\Gamma_{A}^{H}$ associated to
$(A,\rho_{A})$ (see (\ref{weak-Gamma})).

Then, following the Definition 1.9 of \cite{nmra1}, we say that
$A_{H}\hookrightarrow A$ is a weak $H$-cleft extension if there
exists a morphism $h:H\rightarrow A$ in $Reg^{WR}(H,A)$ of right
$H$-comodules such that
\begin{equation}
\label{weak-cleft} \Gamma_{A}^{H}\circ (H\otimes h^{-1})\circ \delta_{H}=
\zeta_{A}\circ (e_{A}\wedge h^{-1}).
\end{equation}

This definition of weak cleft extension is the one used in
\cite{ana1} where was proved that this kind of weak cleft extensions
induce weak crossed products with a left invertible cocycle (see
Definition 4.1 of \cite{ana1})

Also the definition of introduced in \cite{nmra1}   is a
generalization of the one given by Brzezi\'{n}ski \cite{BRZ} (see
\cite{D}, \cite{doi3} for the classical definitions in the Hopf
algebra setting) in the context of entwined structures. Note that,
while in the case of a cleft extension for an entwining structure
$h$ is required to be a comodule morphism and convolution
invertible, here both conditions are replaced by weaker ones. Also,
in \cite{NikaRamonAnaManel} we proved that weak $H$-cleft extensions
are exactly weak $H$-Galois extensions with normal basis (see
Definition 1.8 and Theorem 2.11).

Note that if $h$ is a morphism of right $H$-comodules
$$h\wedge e_{A}=\mu_{A}\co (id_{A}\ot e_{A})\co \rho_{A}\co h=h$$
and if $g=e_{A}\wedge h^{-1}$ we have
$$g\wedge h=(e_{A}\wedge h^{-1})\wedge h=e_{A}\wedge (h^{-1}\wedge
h)=e_{A}\wedge e_{A}=e_{A}$$ and
$$e_{A}\wedge g=e_{A}\wedge (e_{A}\wedge h^{-1})=(e_{A}\wedge e_{A})\wedge h^{-1}=e_{A}\wedge
h^{-1}=g.$$ Therefore, we can assume without loss of generality that
$e_{A}\wedge h^{-1}=h^{-1}$ and (\ref{weak-cleft}) can be expressed as
\begin{equation}
\label{new-weak-cleft} \Gamma_{A}^{H}\circ (H\otimes h^{-1})\circ
\delta_{H}= \zeta_{A}\circ  h^{-1}.
\end{equation}

Moreover, By Remark 1.10 of \cite{nmra1}, we know that if there
exists $h\in Reg^{WR}(H,A)$ of right $H$-comodules
\begin{equation}
\label{caracterizaciongamma} \Gamma_{A}^{H}= (\mu_{A}\ot H)\co (A\ot
(\rho_{A}\co\mu_{A}))\co (((h^{-1}\ot h)\co \delta_{H})\ot A)
\end{equation} and the extension is weak cleft, by Proposition 1.12 of
\cite{nmra1}, we obtain that
$$q_{A}=\mu_{A}\co (A\ot h^{-1})\co \rho_{A}:A\rightarrow A$$
factors through $i_{A}$. Therefore, there exists an unique morphism
$p_{A}:A\rightarrow A_{H}$ such that $q_{A}=i_{A}\co p_{A}$. Then,
$h\wedge h^{-1}=q_{A}\co h$ and, as a consequence, $h\wedge h^{-1}$
admits a factorization through $i_{A}$. }
\end{apart}

\begin{teo}
\label{equivalence-cleft-weak-cleft} Let $H$ be a  weak Hopf algebra
and let $(A, \rho_A)$ be a right $H$-comodule algebra. If there
exists $h\in Reg^{WR}(H,A)$ of right $H$-comodules such that $e_{A}\wedge
h^{-1}=h^{-1}$, the following assertions are equivalent:
\begin{itemize}

\item[(i)] The morphism $h\wedge h^{-1}$ factorizes through the equalizer $i_A$
and $h^{-1}$ satisfies (\ref{cleaving-lambda}).

\item[(ii)] The equality (\ref{new-weak-cleft}) holds.

\end{itemize}

\end{teo}

\begin{dem} If (ii) holds, $A_{H}\hookrightarrow A$ is a weak $H$-cleft
extension and then $h\wedge h^{-1}$ admits a factorization through
$i_{A}$. The equality (\ref{cleaving-lambda}) follows in a similar
way to the  proof given in Proposition \ref{propinversa} using
that $e_{A}\wedge h^{-1}=h^{-1}$.

Conversely, if we assume that (i) holds we have the following:

\begin{itemize}

\item[ ]

\item[ ]$\hspace{0.38cm}\Gamma_{A}^{H} \co (H\ot h^{-1})\co \delta_{H} $

\item[ ]

\item[ ]$=(\mu_{A}\ot H)\co (A\ot (\rho_{A}\co \mu_{A}))\co
(((h^{-1}\ot h)\co \delta_{H})\ot h^{-1})\co \delta_{H} $

\item[ ]

\item[ ]$=(\mu_{A}\ot H)\co  (h^{-1}\ot (\rho_{A}\co (h\wedge h^{-1})))\co \delta_{H} $

\item[ ]

\item[ ]$=(\mu_{A}\ot \Pi_{H}^{L})\co  (h^{-1}\ot (\rho_{A}\co (h\wedge h^{-1})))\co \delta_{H} $

\item[ ]

\item[ ]$=((\mu_{A}\co (A\ot \mu_{A}))\ot (\Pi_{H}^{L}\co \mu_{H}))\co
(h^{-1}\ot h\ot c_{H,A}\ot H)\co (H\ot \delta_{H}\ot (\rho_{A}\co
h^{-1}))\co (H\ot \delta_{H})\co \delta_{H} $

\item[ ]

\item[ ]$=(A\ot \Pi_{H}^{L})\co \mu_{A\ot H}\co (e_{A}\ot ((H\ot ((h^{-1}\ot \lambda_{H})
\co c_{H,H}\co \delta_{H}))\co \delta_{H}))\co \delta_{H} $

\item[ ]

\item[ ]$=(\mu_A\ot \Pi_{H}^{L})\co (e_{A}\ot c_{H,A})\co (\delta_{H}\ot h^{-1})\co \delta_{H} $

\item[ ]

\item[ ]$=(((A\ot \varepsilon_{H})\co \Gamma_{A}^{H} )\ot \Pi_{H}^{L})\co
(H\ot c_{H,A})\co (\delta_{H}\ot h^{-1})\co \delta_{H} $

\item[ ]

\item[ ]$=(A\ot (\varepsilon_{H}\co \mu_{H})\ot  \Pi_{H}^{L}) \co (c_{H,A}\ot c_{H,H})\co
(H\ot c_{H,A}\ot H)\co (\delta_{H}\ot (\rho_{A}\co h^{-1}))\co
\delta_{H} $

\item[ ]

\item[ ]$=(A\ot ( \Pi_{H}^{L}\co \mu_{H}\co (H\ot \Pi_{H}^{L})))\co (c_{H,A}\ot H)\co
 (H\ot (\rho_{A}\co h^{-1}))\co \delta_{H}$

\item[ ]

\item[ ]$=(A\ot ( \Pi_{H}^{L}\co \mu_{H}))\co (c_{H,A}\ot H)\co
 (H\ot (\rho_{A}\co h^{-1}))\co \delta_{H} $

\item[ ]

\item[ ]$= (A\ot ( \Pi_{H}^{L}\co \mu_{H}))\co (c_{H,A}\ot H)\co
 (H\ot ((h^{-1}\ot \lambda_H)\co
c_{H,H}\co \delta_H))\co \delta_{H}$

\item[ ]

\item[ ]$=(h^{-1}\ot \Pi_{H}^{L})\co
c_{H,H}\co \delta_H $

\end{itemize}

\end{dem}

where the first equality follows by (\ref{caracterizaciongamma}),
the second one by the coassociativity of $\delta_{H}$, the third one
by the factorization of $h\wedge h^{-1}$ through $i_{A}$ and the
fourth one because $A$ is a weak entwined module, $\delta_{H}$ is
coassociative and $h$ is a morphism of right $H$-comodules. In the
fifth equality we used the coassociativity of $\delta_{H}$ and
(\ref{cleaving-lambda}) and the sixth one is a consequence of the
coassociativity of $\delta_{H}$ and the properties of $\Pi_{H}^{L}$.
The seventh one follows by (\ref{d)}) and the eight one relies in
the definition of $\Gamma_{A}^{H}$ and the naturality of $c$. Using the
naturality of $c$ and (\ref{PiLmu}) we obtain the ninth equality and
the tenth one follows by (\ref{PiLmuPiL}). Finally, the eleventh one
follows by (\ref{cleaving-lambda}) and the last one by the
naturality of $c$ and  the properties of $\Pi_{H}^{L}$.

On the other hand,
$$\zeta_{A}\co h^{-1}=(A\ot \overline{\Pi}_{H}^{R})\co \rho_{A}\co
h^{-1}=(h^{-1}\ot (\overline{\Pi}_{H}^{R}\co\lambda_{H}))\co
c_{H,H}\co \delta_{H}=(h^{-1}\ot {\Pi}_{H}^{L})\co c_{H,H}\co
\delta_{H},$$ and the proof is complete.

As a corollary we have,

\begin{cor}
\label{prin-cor} Let $H$ be a  weak Hopf algebra and let $(A,
\rho_A)$ be a right $H$-comodule algebra. If $A_{H}\hookrightarrow
A$ is an $H$-cleft extension then it is a weak $H$-cleft extension.
\end{cor}

\begin{dem} If $A_{H}\hookrightarrow
A$ is an $H$-cleft extension, there exists an integral
$f:H\rightarrow A$ convolution invertible and such that the morphism
$f\wedge f^{-1}$ factorizes through the equalizer $i_A$, being $f^{-1}$
the convolution inverse of $f$. Then, $f\in Reg^{WR}(H,A)$, $e_{A}\wedge
f^{-1}=f^{-1}$ and by Proposition \ref{propinversa}, the equality
(\ref{cleaving-lambda}) holds. Therefore, by the previous Theorem we
obtain that $A_{H}\hookrightarrow A$ is a weak $H$-cleft extension.

\end{dem}

\begin{rem}
\label{rem-simplify} {\rm As a consequence of Corollary
\ref{prin-cor}, the results proved in \cite{nmra1} and
\cite{NikaRamonAnaManel} for weak $H$-cleft extensions can be
applied for $H$-cleft extensions. For example,
if   $A_{H}\hookrightarrow A$ is an $H$-cleft extension with cleaving morphism $f$, the morphism $q_A=\mu_A\co (A\ot f^{-1})\co\rho_A$ factors through the equalizer $i_A$, i. e., there exists a morphism $p_A:A\rightarrow A_H$ such that $i_A\co p_A=q_A$. Also, by Lemmas 3.9 and 3.11 of \cite{nmra4} we have the equalities:
\begin{equation}
\label{equcross1}
\mu_{A}\circ (A\otimes e_{A})\circ
\rho_{A}=id_{A},
\end{equation}
\begin{equation}
\label{equcross2} \mu_{A}\circ (q_{A}\otimes f)\circ
\rho_{A}=id_{A},\end{equation}
\begin{equation}
\label{equcross3}\rho_{A}\circ \mu_{A}=(\mu_{A}\otimes H)\circ
(q_{A}\otimes (\rho_{A}\circ \mu_{A}\circ (f\otimes A)))\circ
(\rho_{A}\otimes A),
\end{equation}
\begin{equation}
\label{equcross4}\mu_{A}\circ (i_{A}\otimes f)=\mu_{A}\circ
(q_{A}\otimes A)\circ (\mu_{A}\otimes f)\circ (i_{A}\otimes
(\rho_{A}\circ f)),
\end{equation}
\begin{equation}
\label{equcross5} p_{A}\circ \mu_A\circ (i_{A}\otimes
A)=\mu_{A_{H}}\circ (A_{H}\otimes p_{A}).
\end{equation}
}
\end{rem}

\begin{defin}\label{cleftequivalentes}

{\rm  Let $H$ be a  weak Hopf algebra. Two $H$-cleft extensions
$A_{H}\hookrightarrow A$ and $B_{H}\hookrightarrow B$ are equivalent
(written $A_{H}\hookrightarrow A \sim B_{H}\hookrightarrow B$) if $A_H=B_H$ and there is a
morphism of right $H$-comodule algebras $T:A\rightarrow B$ such that $T\co i_A=i_B$.

Note that, if the $H$-cleft extensions $A_{H}\hookrightarrow A$ and $A_{H}\hookrightarrow B$
are equivalent and $f$ is a cleaving morphism for $A_{H}\hookrightarrow A$, it is easy
to show that
$g=T\co f$ is  a cleaving morphism for $A_{H}\hookrightarrow B$ with $g^{-1}=T\co f^{-1}$. Under
these conditions, $T$ is an isomorphism. If $f$ is the cleaving morphism associated
to $A_{H}\hookrightarrow A$, we define the morphisms
$$\gamma_{A}=(p_{A}\ot H)\co \rho_{A}:A\rightarrow A_{H}\ot H,$$
$$\chi_{A}=\mu_{A}\co (i_{A}\ot f):A_{H}\ot H\rightarrow A$$
and
$$\gamma_{B}=(p_{B}\ot H)\co \rho_{B}:B\rightarrow A_{H}\ot H,$$
$$\chi_{B}=\mu_{B}\co (i_{B}\ot g):A_{H}\ot H\rightarrow B$$
where $p_{A}$ and $p_{B}$ are the factorizations of $q_A=\mu_A\co (A\ot f^{-1})\co\rho_A$, $q_B=\mu_B\co (A\ot g^{-1})\co\rho_B$ and $i_{A}$, $i_{B}$ the corresponding equalizer morphisms.  Then,
$$\chi_{B}\co \gamma_{A}= \mu_{B}\co ((i_{B}\co p_{A})\ot g)\co \rho_{A}=\mu_{B}\co
((T\co q_{A})\ot (T\co f))\co \rho_{A}=T\co \mu_{A}\co (A\ot (f^{-1}\wedge f))\co \rho_{A}$$
$$=T\co \mu_{A}\co (A\ot e_{A})\co \rho_{A}=T,$$
and
$$i_{B}\co p_{B}\co T=\mu_{B}\co (B\ot (T\co f^{-1}))\co \rho_{B}\co T=T\co q_{A}
=T\co i_{A}\co p_{A}=i_{B}\co p_{A}.$$
If we define $T^{-1}:B\rightarrow A$ by $T^{-1}=\chi_{A}\co \gamma_{B}$ we have:

$$T\co T^{-1}= T\co \mu_{A}\co ((i_{A}\co p_{B})\ot f)\co
\rho_{B}=\mu_{B}\co ((T\co i_A\co p_B)\ot (T\co f))\ot \rho_B$$
$$=\mu_{B}\co (B\ot e_{B})\co \rho_{B}=id_{B}$$
and, by a similar proof $T^{-1}\co T=id_{A}$. Therefore, $T$ is an isomorphism.

Obviously, $" \sim "$ is an equivalence relation, and we will denote
by $[B_{H}\hookrightarrow B]$
the isomorphisms class of the $H$-cleft extension $B_{H}\hookrightarrow B$.
}
\end{defin}

It is a well known fact in the Hopf algebra setting that, if $A_{H}\hookrightarrow A$ is
an $H$-cleft extension with convolution invertible integral $f$, the morphism $h=\mu_A\co (f\ot (f^{-1}\co\eta_H)))$
is a total integral convolution invertible. There is a similar result for weak Hopf algebras
although we want to point out that the proof is very different, because in our case
$\delta_H\co \eta_H\neq\eta_H\ot\eta_H$. Actually, in order to give the convolution inverse
of this morphism we assume that the weak Hopf algebra is cocommutative. This hypothesis can be removed in the classical case because for Hopf algebras the morphisms $\Pi^{L}_{H}$,
$\Pi^{R}_{H}$ $\overline{\Pi}^{L}_{H}$ and $\overline{\Pi}^{R}_{H}$ trivialize.

\begin{prop}
\label{existenciaintegraltotal}

Let $H$ be a weak Hopf algebra with invertible antipode.  If $A_{H}\hookrightarrow A$ is an $H$-cleft
extension with cleaving morphism $f$,  then  $h=\mu_A\co (f\ot (f^{-1}\co\eta_H)))$
is a total integral. Moreover, if $H$ is cocommutative $h$ is convolution invertible.

\end{prop}

\begin{dem}
We define $h=\mu_A\co (f\ot (f^{-1}\co\eta_H)))$, where $f^{-1}$ is the convolution
inverse of $f$. Using Proposition \ref{propinversa}, the properties of the antipode
and (\ref{composiciones}),
\begin{itemize}

\item[ ]$\hspace{0.38cm} (A\ot (\varepsilon_H\co\mu_H))\co (c_{H,A}\ot H)\co (H\ot (\rho_A\co f^{-1}\co \eta_H))$

\item[ ]

\item[ ]$=(A\ot (\varepsilon_H\co\mu_H))\co (c_{H,A}\ot H)\co (H\ot (( f^{-1}\ot \lambda_H)\co c_{H,H}\co\delta_H\co \eta_H))$

\item[ ]

\item[ ]$=((\varepsilon_H\co\mu_H\co c_{H,H})\ot f^{-1})\co (\lambda_H^{-1}\ot (\delta_H\co \eta_H))$

\item[ ]

\item[ ]$=f^{-1}\co \Pi^{L}_{H}\co \lambda_H^{-1}$

\item[ ]

\item[ ]$=f^{-1}\co \overline{\Pi}^{R}_{H}.$

\end{itemize}

As a consequence, using that $f$ is an integral and $A$ a comodule algebra,
we obtain that
\begin{equation}
\label{equalityforh}
h=f\wedge(f^{-1}\co \overline{\Pi}^{R}_{H}),
\end{equation}
and then
we get that $h$ is total because, by (\ref{equalityforh}) and (\ref{deltaPIRbarra}),
$$h\co \eta_{H}= (f\wedge(f^{-1}\co \overline{\Pi}^{R}_{H}))\co \eta_{H}=
(f\wedge f^{-1})\co \eta_{H}=\eta_{A}.$$

Moreover, $h$ is an integral. Indeed:

\begin{itemize}

\item[ ]$\hspace{0.38cm}\rho_A\co h$

\item[ ]

\item[ ]$=\mu_{A\ot H}\co ((\rho_A\co f)\ot (\rho_A\co f^{-1}\co \eta_H))$

\item[ ]

\item[ ]$=(A\ot (\lambda_H\co\lambda_H^{-1}))\co\mu_{A\ot H}\co(((f\ot H)\co\delta_H)\ot ((f^{-1}\ot\lambda_H)\co c_{H,H}\co\delta_H\co \eta_H))$

\item[ ]

\item[ ]$=(\mu_A\ot H)\co (f\ot (c_{H,A}\co (\lambda_{H}\ot f^{-1})\co (\mu_{H}\ot H)\co (H\ot c_{H,H})\co ((\delta_{H}\co \eta_{H})\ot \lambda_{H}^{-1})))\co \delta_{H}$

\item[ ]

\item[ ]$=(\mu_A\ot H)\co (f\ot (c_{H,A}\co (\lambda_{H}\ot f^{-1})\co (H\ot \Pi_{H}^{L})\co \delta_{H}\co \lambda_{H}^{-1}))\co \delta_{H}$

\item[ ]

\item[ ]$=(\mu_A\ot H)\co (f \ot ((f^{-1}\co  \Pi_{H}^{L}\co \lambda_{H}^{-1})\ot H)\co \delta_{H})\co \delta_{H}$

\item[ ]

\item[ ]$=(h\ot H)\co \delta_H$

\end{itemize}

The first equality follows because $A$ is a right $H$-comodule algebra, the second one uses that $\lambda_{H}$ is an isomorphism, the third and the fifth ones are a consequence of the antimultiplicative property for $\lambda_{H}$,  $\lambda_{H}^{-1}$ and the naturality of $c$, the fourth one follows by  (\ref{deltaPIL}) and, finally, the last one by the coassociativity of $\delta_{H}$ and (\ref{composiciones}).

Now we assume that $H$ is cocommutative. We define $h^{-1}=\mu_A\co
((f\co\eta_H)\ot f^{-1})$. Following a similar way to the one
developed for $h$, it is not difficult to prove the equalities:
\begin{equation}
\label{expresioninversah}
h^{-1}=(f\co \Pi^{R}_{H})\wedge f^{-1}
\end{equation}

and
\begin{equation}
\label{expresionauxiliar}
\mu_A\co (f^{-1}\ot (f\co \eta_H))=\mu_A\co (f^{-1}\ot (f\co \Pi^{R}_{H}\co \lambda_H))\co c_{H,H}\co\delta_H.
\end{equation}

As a consequence of the last equation, taking into account that $H$ is cocommutative we obtain that
\begin{equation}
\label{expresioncruzada}
\mu_A\co ((f^{-1}\co \eta_H)\ot (f\co \eta_H))=\eta_H.
\end{equation}

We conclude the proof showing that $h^{-1}$ is the convolution inverse of $h$. Condition (c2) follows because, by (\ref{expresioncruzada}),
$$h\wedge h^{-1}=f\wedge f^{-1}=(A\ot(\varepsilon_H\co\mu_H))\co ((\rho_A\co \eta_A)\ot H).$$

Moreover

\begin{itemize}

\item[ ]$\hspace{0.38cm} h^{-1}\wedge h$

\item[ ]

\item[ ]$=\mu_A\co (\mu_A\ot A)\co ((f\co \eta_H)\ot (f^{-1}\wedge f)\ot (f^{-1}\co \eta_H))$

\item[ ]

\item[ ]$=\mu_A\co (\mu_A\ot (\varepsilon_H\co \mu_H)\ot A)\co ((f\co \eta_H)\ot c_{H,A}\ot H\ot (f^{-1}\co \eta_H))\co (H\ot (\rho_A\co\eta_A))$

\item[ ]

\item[ ]$=\mu_A\co (A\ot (\varepsilon_H\co\mu_H\co c_{H,H}\co (\overline{\Pi}^{R}_{H}\ot H))\ot A)\co ((\rho_A\co f\co \eta_H)\ot H\ot (f^{-1}\co \eta_H))$

\item[ ]

\item[ ]$=\mu_A\co((f\co \Pi^{R}_{H})\ot (f^{-1}\co\eta_H))$

\item[ ]

\item[ ]$=h\co \Pi^{R}_{H}$

\item[ ]

\item[ ]$=f^{-1}\wedge f$

\item[ ]

\item[ ]$=(A\ot(\varepsilon_H\co\mu_H))\co(c_{H,A}\ot H)\co (H\ot(\rho_A\co \eta_A)).$

\end{itemize}

In the foregoing calculations, the first equality follows by the definition of $h$ and $h^{-1}$; the second and the last ones  because $f$ is convolution invertible; the third because $A$ is an $H$-comodule algebra; the fourth uses that $f$ is an integral and (\ref{deltaPIRbarra}); finally,
the fifth and the sixth ones are consequence of the definition of the total integral $h$.

The proof for the condition (c3) for $h$ follows a similar pattern and we leave the
details to the reader.
\end{dem}

\begin{rem}
{\rm
As a consequence of the previous proposition, in the cocommutative setting we can assume that the integral is total.

}
\end{rem}

In the following definition we recall the notion of weak left $H$-module algebra introduced in \cite{NikaRamon6}.

\begin{defin}
\label{weak-H-mod} {\rm Let $H$ be a weak Hopf algebra. We will say
that $A$ is a weak left $H$-module algebra if there exists a
morphism $\varphi_{A}:H\ot A\rightarrow A$ satisfying:

\begin{itemize}

\item[(d1)] $\varphi_{A}\co (\eta_{H}\ot A)=id_{A}$.

\item[(d2)] $\varphi_{A}\co (H\ot \mu_{A})=\mu_{A}\co
(\varphi_{A}\ot \varphi_{A})\co (H\ot c_{H,A}\ot A)\co
(\delta_{H}\ot A\ot A).$

\item[(d3)] $\varphi_{A}\co (\mu_H\ot \eta_{A})=\varphi_{A}\co (H\ot
(\varphi_{A}\co (H\ot \eta_{A}))).$

\end{itemize}

and any of the following equivalent conditions hold:

\begin{itemize}

\item[(d4)] $\varphi_{A}\co (\Pi^{L}_{H}\ot A)=\mu_{A}\co ((\varphi_{A}\co (H\ot
\eta_{A})\ot A).$

\item[(d5)] $\varphi_{A}\co (\overline{\Pi}^{L}_{H}\ot A)=\mu_{A}\co c_{A,A}\co
((\varphi_{A}\co (H\ot \eta_{A})\ot A).$

\item[(d6)] $\varphi_{A}\co (\Pi^{L}_{H}\ot \eta_{A})=\varphi_{A}\co (H\ot \eta_{A}).$

\item[(d7)] $\varphi_{A}\co (\overline{\Pi}^{L}_{H}\ot\eta_{A})=
\varphi_{A}\co (H\ot \eta_{A}).$

\item[(d8)] $\varphi_{A}\co (H\ot (\varphi_{A}\co (H\ot
\eta_{A})))=((\varphi_{A}\co (H\ot \eta_{A}))\ot (\varepsilon_{H}\co
\mu_{H}))\co (\delta_{H}\ot H).$

\item[(d9)] $\varphi_{A}\co (H\ot (\varphi_{A}\co (H\ot
\eta_{A})))=((\varepsilon_{H}\co \mu_{H})\ot (\varphi_{A}\co (H\ot \eta_{A})))\co (H\ot c_{H,H})\co (\delta_H\ot H).$

\end{itemize}

If we replace (d3) by

\begin{itemize}

\item[(d3-1)] $\varphi_A\co (\mu_H\ot A)=\varphi_A\co (H\ot \varphi_A)$

\end{itemize}

we will say that $(A, \varphi_A)$ is a left $H$-module algebra.

}
\end{defin}

\begin{rem}
\label{measuring}
{\rm
Note that by (d4) and (d5) if the weak Hopf algebra is cocommutative the morphism $\varphi_A\co (H\ot \eta_A)$ factors through the center of $A$. Moreover, if $H$ is a Hopf algebra and $(A, \varphi_A)$ is a weak left $H$-module algebra, conditions (d4)-(d9) imply that $\varepsilon_H\ot \eta_A=\varphi_A\co (H\ot \eta_A)$. As a consequence the equality (d3) is always true and $\varphi_A$ is a weak action of $H$ on $A$ (see \cite{BCM}).
}
\end{rem}

\begin{prop}
\label{estructuramoduloinducida}

Let $H$ be a cocommutative weak Hopf algebra. If $A_{H}\hookrightarrow A$ is an $H$-cleft extension with cleaving morphism  $f$, the pair $(A_H, \varphi_{A_H})$ is a  weak left $H$-module algebra, being
$\varphi_{A_H}$ the factorization of the morphism
$$\varphi_A=\mu_A\co (A\ot (\mu_A\co c_{A,A}))\co (((f\ot f^{-1})\co \delta_{H})\ot i_{A} )$$ through the equalizer $i_A$.

\end{prop}

\begin{dem} If $A_{H}\hookrightarrow A$ is
a $H$-cleft extension, by Corollary \ref{prin-cor}, we have that $A_{H}\hookrightarrow A$ is a weak $H$-cleft extension, and, by  Proposition 1.15 of \cite{nmra1},  we know that $\varphi_{A_H}$ factors through the equalizer $i_A$ and satisfies (b2). Moreover it is not difficult to see that

 \begin{equation}
\label{alternativamodulo}
\varphi_{A_H}=p_A\co \mu_A\co (f\ot i_A),
\end{equation}

and then (b1) holds.

As far as (b3),

\begin{itemize}

\item[ ]$\hspace{0.38cm} \varphi_{A_{H}}\co (H\ot(\varphi_{A_{H}}\co (H\ot \eta_{A_{H}})))$

\item[ ]

\item[ ]$=p_A\co\mu_A\co(f\ot(q_A\co f))$

\item[ ]

\item[ ]$=p_A\co\mu_A\co(f\ot(f\wedge f^{-1}))$

\item[ ]

\item[ ]$=((p_A\co \mu_A)\ot (\varepsilon_H\co \mu_H))\co (f\ot (\rho_A\co \eta_A)\co H)$

\item[ ]

\item[ ]$=(p_A \ot (\varepsilon_H\co \mu_H))\co (A\ot \overline{\Pi}^{R}_{H} \ot H)\co ((\rho_A\co f)\ot H)$

\item[ ]

\item[ ]$=((p_A\co q_A\co f)\ot (\varepsilon_H\co \mu_H))\co (\delta_H\ot H)$

\item[ ]

\item[ ]$=((p_A\co (f\wedge f^{-1}))\ot (\varepsilon_H\co \mu_H))\co (\delta_H\ot H)$

\item[ ]

\item[ ]$=(p_A\ot(\varepsilon_H\co \mu_H)\ot (\varepsilon_H\co \mu_H))\co ((\rho_A\co \eta_A)\ot \delta_H\ot H)$

\item[ ]

\item[ ]$=(p_A\ot(\varepsilon_H\co \mu_H))\co ((\rho_A\co \eta_A)\ot\mu_H)$

\item[ ]

\item[ ]$=(p_A\co(f\wedge f^{-1}))\co \mu_H$

\item[ ]

\item[ ]$=\varphi_{A_H}\co (\mu_{H}\ot \eta_{A_H}),$

\end{itemize}

where the first equality follows by (\ref{alternativamodulo}); the second because $q_A\co f=f\wedge f^{-1}$; the third, sixth and eighth ones by (c2); in the fourth one we use that $A$ is a right $H$-comodule algebra; in the fifth (\ref{PiRbarramu}) and that $f$ is an integral; the seventh equality follows because $H$ is a weak Hopf algebra; finally, in the last one we use that $p_A\co (f\wedge f^{-1})=\varphi_{A_H}\co (H\ot \eta_{A_H})$.

It only remains to show one of the equivalent conditions (d4)-(d9). We get (d6):

\begin{itemize}

\item[ ]

\item[ ]$\hspace{0.38cm} \varphi_{A_H}\co (\Pi^{L}_{H}\ot \eta_{A_H})$

\item[ ]

\item[ ]$=p_A\co (f\wedge f^{-1})\co \Pi^{L}_{H}$

\item[ ]

\item[ ]$=p_A\co(A\ot(\varepsilon_H\co \mu_H))\co ((\rho_A\co \eta_A)\ot\Pi^{L}_{H})$

\item[ ]

\item[ ]$=p_A\co (f\wedge f^{-1})$

\item[ ]

\item[ ]$=\varphi_{A_H}\co (H\ot \eta_{A_H}).$

\end{itemize}

\end{dem}

As a consequence of the previous result we have the following corollary.

\begin{cor}
\label{casotwisted}

In the conditions of Proposition \ref{estructuramoduloinducida},

\begin{itemize}

\item[(i)] $\varphi_{A_H}=\varphi_{A_H}\co (\overline{\Pi}_{H}^{L}\ot A_H)$ if and only if $\mu_A\co (f\ot i_A)=\mu_A\co c_{A,A}\co (f\ot i_A)$.

\item[(ii)] $\varphi_{A_H}=\varphi_{A_H}\co (\Pi_{H}^{L}\ot A_H)$ if and only if $\mu_A\co (f^{-1}\ot i_A)=\mu_A\co c_{A,A}\co (f^{-1}\ot i_A)$.

\end{itemize}

\end{cor}

\begin{dem}

We will show (i). Part (ii) is similar and we leave the details to
the reader. Assume that $\varphi_{A_H}=\varphi_{A_H}\co
(\overline{\Pi}_{H}^{L}\ot A_H)$. Taking into account (d5) and
Proposition \ref{estructuramoduloinducida}  we obtain that
$\varphi_{A}=\mu_A\co c_{A,A}\co ((\varphi_{A}\co (H\ot\eta_A))\ot
i_A)$. The result follows composing in this equality with $(H\ot
A_H\ot f)\co (H\ot c_{H,A_H})\co (\delta_H\ot A_H)$ on the right and
with $\mu_A$ on the left. Indeed, using the definition of
$\varphi_A$, that $f$ is a convolution invertible integral,(b2), the
properties of the equalizer $i_A$ and that $(A, \rho_A)$ is a right
$H$-comodule algebra,

\begin{itemize}

\item[ ]

\item[ ]$\hspace{0.38cm} \mu_A\co (\varphi_A\ot A)\co (H\ot ((i_A\ot f)\co c_{H,A_H}))\co (\delta_H\ot A_H)$

\item[ ]

\item[ ]$=\mu_A\co (A\ot \mu_A)\co (f\ot i_A\ot f^{-1}\wedge f))\co (H\ot c_{H,A_H})\co (\delta_H\ot A_H)$

\item[ ]

\item[ ]$=(\mu_A\ot (\varepsilon_H\co \mu_H))\co (f\ot c_{H,A}\ot H)\co (\delta_H\ot ((\mu_A\ot H)\co (i_A\ot (\rho_A\co \eta_A)))$

\item[ ]

\item[ ]$=(A\ot \varepsilon_H)\co \mu_{A\ot H}\co ((\rho_A\co f)\ot ((A\ot \overline{\Pi}_{H}^{R})\co \rho_A\co i_A))$

\item[ ]

\item[ ]$=(A\ot \varepsilon_H)\co \rho_A\co \mu_A\co (f\ot i_A)$

\item[ ]

\item[ ]$=\mu_A\co (f\ot i_A).$

\end{itemize}

On the other hand, by similar arguments,

\begin{itemize}

\item[ ]

\item[ ]$\hspace{0.38cm} \mu_A\co ((\mu_A\co c_{A,A}\co ((\varphi_{A}\co (H\ot\eta_A))\ot i_A)\ot f)\co (H\ot c_{H,A_H})\co (\delta_H\ot A_H)$

\item[ ]

\item[ ]$=\mu_A\co (i_A\ot f\wedge f^{-1}\wedge f)\co c_{H,A_H}$

\item[ ]

\item[ ]$=\mu_A\co (i_A\ot f)\co c_{H,A_H}$

\item[ ]

\item[ ]$=\mu_A\co c_{A,A}\co (f\ot i_{A}).$

\end{itemize}

Conversely, by the hypothesis,(d5) and the equality $f\wedge f^{-1}=\varphi_A\co (H\ot \eta_A)$,

\begin{itemize}

\item[ ]

\item[ ]$\hspace{0.38cm} i_A\co \varphi_{A_H}$

\item[ ]

\item[ ]$=\mu_A\co ((\mu_A\co (f\ot i_A))\ot f^{-1})\co (H\ot c_{H,A_H})\co (\delta_H\ot A_H)$

\item[ ]

\item[ ]$=\mu_A\co (i_A\ot f\wedge f^{-1})\co c_{H,A_H}$

\item[ ]

\item[ ]$=\mu_A\co c_{A,A}\co ((\varphi_A\co (H\ot \eta_A))\co i_A)$

\item[ ]

\item[ ]$=\varphi_{A}\co (\overline{\Pi}_{H}^{L}\ot i_A)$

\item[ ]

\item[ ]$=i_A\co \varphi_{A_H}\co (\overline{\Pi}_{H}^{L}\ot A_H),$

\end{itemize}

and then $\varphi_{A_H}=\varphi_{A_H}\co (\overline{\Pi}_{H}^{L}\ot A_H)$.

\end{dem}

\section{crossed systems for weak Hopf algebras}

In the first part of this section we generalize the theory of crossed systems over a Hopf algebra given by Doi in \cite{doi1} to the weak setting. Taking into account the theory developed in the previous section, being $H$ a cocommutative weak Hopf algebra, we will obtain a bijective correspondence between the isomorphisms classes of $H$-cleft extensions  $[A_{H}\hookrightarrow A]$ and the equivalence classes of crossed systems for $H$ over $A_H$.

By Propositions 1.4, 1.6 and 1.8 of \cite{NikaRamon6}, if $H$ is a weak Hopf algebra,  the morphisms
\begin{equation}
\label{omega-Lnew} \Omega_{H}^{L}= ((\varepsilon_{H}\co
\mu_{H})\ot H\ot H)\co \delta_{H\ot H}:H\ot H\rightarrow H\ot H
\end{equation}
\begin{equation}
\label{omega-Rnew} \Omega_{H}^{R}= (H\ot H\ot (\varepsilon_{H}\co
\mu_{H}))\co \delta_{H\ot H}:H\ot H\rightarrow H\ot H
\end{equation}
are idempotent and satisfy that $\mu_H=\mu_H\co \Omega_{H}^{L}=\mu_H\co \Omega_{H}^{R}$. Also, $\Omega_{H}^{L}$ is a morphism of left and right $H$-modules for the usual
regular actions because
\begin{equation}
\label{omega-L-1} \Omega_{H}^{L}=((\mu_{H}\co (H\ot
\Pi_{H}^{L}))\ot H)\co(H\ot \delta_{H})= (H\ot (\mu_{H}\co
(\overline{\Pi}_{H}^{R}\ot H)))\co ((c_{H,H}\co \delta_{H})\ot H).
\end{equation}

Moreover, if $H$ is cocommutative, $\Omega_{H}^{L}=\Omega_{H}^{R}$ and we denote it by $\Omega_{H}^{2}$. Then, we have
\begin{equation}
\label{Omegadelta}
(\Omega_{H}^{2}\ot H\ot H)\co \delta_{H\ot H}=(H\ot H\ot \Omega_{H}^{2})\co \delta_{H\ot H}=\delta_{H\ot H}\co \Omega_{H}^{2}.
\end{equation}

Following Definition 1.18 of \cite{NikaRamon6} we have:

\begin{defin}
 \label{regunoydos}
 {\rm  Let $H$ be a cocommutative weak Hopf algebra and
$(A,\varphi_{A})$ be a weak left $H$-module algebra. We define
$Reg_{\varphi_{A}}(H,A)$, as the set of morphisms $h:H\rightarrow
A$ such that  there exists a morphism $h^{-1}:H\rightarrow A$ (the
convolution inverse of $h$) satisfying the following equalities:
\begin{itemize}
\item[(e1)] $h\wedge h^{-1}= h^{-1}\wedge
h=u_{1},$
\item[(e2)] $h\wedge h^{-1}\wedge h=h,$
\item[(e3)] $h^{-1}\wedge h\wedge h^{-1}=h^{-1},$
\end{itemize}
where $u_{1}=\varphi_{A}\co (H\ot \eta_{A})$.

In a similar way, $Reg_{\varphi_{A}}(H\ot H,A)$ is the set of morphisms
$\sigma:H\ot H\rightarrow A$ such that there exists a morphism
$\sigma^{-1}:H\ot H\rightarrow A$ satisfying:
\begin{itemize}
\item[(f1)] $\sigma\wedge \sigma^{-1}= \sigma^{-1}\wedge
\sigma=u_{2}$
\item[(f2)] $\sigma\wedge \sigma^{-1}\wedge \sigma=\sigma.$
\item[(f3)] $\sigma^{-1}\wedge \sigma\wedge \sigma^{-1}=\sigma^{-1}.$
\end{itemize}
where $u_{2}=\varphi_A\co(H\ot u_{1})$.

Note that, by (d3) of Definition \ref{weak-H-mod},
\begin{equation}
\label{u1-u2}
u_{2}=u_{1}\co \mu_{H}
\end{equation}
and if $\sigma\in Reg_{\varphi_{A}}(H\ot H,A)$, in \cite{NikaRamon6}, we prove that
\begin{equation}
\label{sigmaOmega}
\sigma=\sigma\co \Omega_{H}^{2}
\end{equation}
Analogously, $\sigma^{-1}=\sigma^{-1}\co \Omega_{H}^{2}.$

Also, by Proposition 1.19 of \cite{NikaRamon6} we know that, if $H$ is a weak Hopf algebra and $(A,\varphi_{A})$
 a weak left $H$-module algebra such that there exists $h:H\rightarrow A$
satisfying that:
$$ h\wedge h^{-1}=
h^{-1}\wedge h=u_{1}, \;\; h\wedge h^{-1}\wedge h=h,
\;\;h^{-1}\wedge h\wedge h^{-1}=h^{-1},$$ the following equalities
are equivalent
\begin{equation}
\label{heta}
h\co \eta_{H}=\eta_{A},
\end{equation}
\begin{equation}
\label{hPiL} h\co \Pi_{H}^{L}=u_{1},
\end{equation}
\begin{equation}
\label{hbarPiL} h\co \overline{\Pi}_{H}^{L}=u_{1}.
\end{equation}
If $\lambda_{H}$ is an isomorphism they are equivalent to
\begin{equation}
\label{hPiR} h\co \Pi_{H}^{R}=u_{1}\co \lambda_{H},
\end{equation}
and
\begin{equation}
\label{hbarPiR} h\co \overline{\Pi}_{H}^{R}=u_{1}\co \lambda_{H}^{-1}.
\end{equation}

In a similar way, it is possible to see that, if
$\sigma:H\ot H\rightarrow A$ is a morphism such that

$$\sigma\wedge \sigma^{-1}= \sigma^{-1}\wedge
\sigma=u_{2}, \;\;\sigma\wedge \sigma^{-1}\wedge \sigma=\sigma, \;\;\sigma^{-1}\wedge \sigma\wedge \sigma^{-1}=\sigma^{-1},$$
  the following equalities are
equivalent:
\begin{equation}
\label{sigma-eta-l}
\sigma\co (\eta_{H}\ot H)= u_{1},
\end{equation}
\begin{equation}
\label{sigma-PiL-l}\sigma\co (\Pi_{H}^{L}\ot H)\co \delta_{H}=u_{1},
\end{equation}
\begin{equation}
\label{sigma-PiL-c}
\sigma\co c_{H,H}\co(H\ot \overline{\Pi}_{H}^{L})\co \delta_{H}=u_{1}.
\end{equation}
If  the antipode $\lambda_{H}$ is an isomorphism (\ref{sigma-eta-l})-(\ref{sigma-PiL-c}) are
equivalent to
\begin{equation}
\label{sigma-PiR-l}\sigma\co (\Pi_{H}^{R}\ot \lambda_{H})\co \delta_{H}=u_{1}\co \lambda_{H},
\end{equation}
and
\begin{equation}
\label{sigma-PiR-c} \sigma\co c_{H,H}\co(\lambda_{H}^{-1}\ot \overline{\Pi}_{H}^{R})\co \delta_{H}=
u_{1}\co \lambda_{H}^{-1}.
\end{equation}

Finally,  the following assertions are
equivalent:
\begin{equation}
\label{sigma-eta-r} \sigma\co (H\ot \eta_{H})= u_{1},
\end{equation}
\begin{equation}
\label{sigma-PiR-r}\sigma\co (H\ot \Pi_{H}^{R})\co \delta_{H}=u_{1},
\end{equation}
\begin{equation}
\label{sigma-PiR-r-c}\sigma\co c_{H,H}\co( \overline{\Pi}_{H}^{R}\ot H)\co \delta_{H}=u_{1}.
\end{equation}

If the antipode $\lambda_{H}$ is an isomorphism (\ref{sigma-eta-r})-(\ref{sigma-PiR-r-c}) are
equivalent to
\begin{equation}
\label{sigma-PiL-r}\sigma\co (\lambda_{H}\ot\Pi_{H}^{L})\co \delta_{H}=u_{1}\co \lambda_{H},
\end{equation}
and
\begin{equation}
\label{sigma-PiL-r-c} \sigma\co c_{H,H}\co( \overline{\Pi}_{H}^{L}\ot \lambda_{H}^{-1})\co \delta_{H}=
u_{1}\co \lambda_{H}^{-1}.
\end{equation}
}
\end{defin}

\begin{prop}
\label{f-f-1-eta} Let $H$ be a weak Hopf algebra and
$(A,\varphi_{A})$ be a weak left $H$-module algebra. If there exists
$h:H\rightarrow A$ satisfying the following equalities:
$$ h\wedge h^{-1}=
h^{-1}\wedge h=u_{1}, \;\; h\wedge h^{-1}\wedge h=h,
\;\;h^{-1}\wedge h\wedge h^{-1}=h^{-1},$$ we have that $h\co
\eta_{H}=\eta_{A}$ if and only if $h^{-1}\co \eta_{H}=\eta_{A}$.
\end{prop}

\begin{dem}
If $h\co \eta_{H}=\eta_{A}$, by (\ref{hPiL}) and (\ref{deltaPIL}),
we have
$$h^{-1}\co
\eta_{H}=(h^{-1}\wedge u_{1})\co \eta_{H}=(h^{-1}\wedge (h\co
\Pi_{H}^{L}))\co \eta_{H}=u_{1}\co \eta_{H}=\eta_{A}.
$$
Conversely, if $h^{-1}\co \eta_{H}=\eta_{A}$, by similar arguments,
$$h\co
\eta_{H}=(h\wedge u_{1})\co \eta_{H}=(h\wedge (h^{-1}\co
\Pi_{H}^{L}))\co \eta_{H}=u_{1}\co \eta_{H}=\eta_{A}.
$$

\end{dem}

\begin{defin}
\label{stmacruzado}
{\rm Let $H$ be a weak Hopf algebra, $(A, \varphi_A)$ a weak left $H$-module
algebra and $\sigma:H\ot H\rightarrow A$ a morphism satisfying (f1)-(f3). We say that
$(\varphi_A, \sigma)$ is a crossed system for $H$ over $A$ if the following conditions hold:

\begin{itemize}
\item[(g1)] $\mu_A\co (A\ot \varphi_A)\co (\sigma\ot \mu_H\ot A)\co (\delta_{H\ot H}\ot A)$
\item[ ]$\hspace{0.38cm} =\mu_A\co ((\varphi_A\co (H\ot \varphi_A))\ot A)\co (H\ot H\ot c_{A,A})\co (H\ot H\ot \sigma\ot A)\co (\delta_{H\ot H}\ot A).$
\item[(g2)] $((\varphi_A\co(H\ot \sigma)\co (\Omega_{H}^{L}\ot H)))\wedge (\sigma\co (H\ot \mu_H))=
((\sigma\ot \varepsilon_H)\co (H\ot \Omega_{H}^{L})\wedge (\sigma\co (\mu_H\ot H)).$
\item[(g3)] $\sigma\co (H\ot \eta_H)=\sigma\co (\eta_H\ot H)=\varphi_A\co (H\ot \eta_A).$
\end{itemize}

Note that equalities (g1) and (g3) are the same that the ones given by Doi in \cite{doi1},
while (g2) is slightly different on account of it includes the idempotent morphism $\Omega_{H}^{2}$.
Anyway, when $H$ is a Hopf algebra, $\Omega_{H}^{2}=id_{H\ot H}$ and our definition is a
generalization to the weak Hopf algebra setting of the one introduced in \cite{doi1}.

Also, it is clear that our condition (f1) over $\sigma$ implies that
it is left invertible in the sense of Definition 4.1 of \cite{ana1}.
In any case, to obtain the main results of this paper and a good
cohomological interpretation, we need the right invertibility, that
is (f1) and (f2), (f3).

Moreover, if the weak Hopf algebra $H$ is cocommutative the morphism $\sigma$ is in $Reg_{\varphi_{A}}(H\ot H,A)$ and we can remove the morphism $\Omega_{H}^{2}$ in condition (g2). Indeed, by (\ref{sigmaOmega}) and using that $\Omega_{H}^{L}=\Omega_{H}^{2}$ is a morphism of left and right $H$-modules we obtain that (g2) is equivalent to
\begin{equation}
\label{new-two-cocycle}
(\varphi_A\co(H\ot \sigma))\wedge (\sigma\co (H\ot \mu_H))=
((\sigma\ot \varepsilon_H))\wedge (\sigma\co (\mu_H\ot H)).
\end{equation}
Moreover, (\ref{new-two-cocycle}) is equivalent to
\begin{equation}
\label{new-two-cocycle-2}
\mu_{A}\co ((\varphi_{A}\co (H\ot \sigma^{-1}))\ot \sigma)\co (H\ot H\ot c_{H,H}\ot H)\co (H\ot H\ot H\ot c_{H,H})\co (\delta_{H\ot H}\ot H)
\end{equation}
$$=\mu_{A}\co (\sigma\ot \sigma^{-1})\co (H\ot \mu_{H}\ot \mu_{H}\ot H)\co
\delta_{H\ot H\ot H}$$
and to
\begin{equation}
\label{new-two-cocycle-3}
\mu_{A}\co (\sigma^{-1}\ot (\varphi_{A}\co (H\ot \sigma)))\co (\delta_{H\ot H}\ot H)
=(\sigma\co (\mu_{H}\ot H))\wedge (\sigma^{-1}\co (H\ot \mu_{H}))
\end{equation}

Two crossed systems for $H$ over $A$, $(\varphi_A, \sigma)$ and
$(\phi_A, \tau)$ are said to be equivalent, denoted by
$$(\varphi_A,
\sigma) \approx (\phi_A, \tau),$$  if $\varphi_A\co (H\ot
\eta_A)=\phi_A\co (H\ot \eta_A)$ and there exists $h$ in
$Reg_{\varphi_{A}}(H,A)\cap Reg_{\phi_{A}}(H,A)$ with $h\co
\eta_H=\eta_A$ and such that
\begin{equation}\label{relacionfis}
\varphi_A=\mu_A\co (\mu_A\ot A)\co (h\ot \phi_A\ot h^{-1})\co (\delta_H\ot c_{H, A})\co (\delta_H\ot A),
\end{equation}
\begin{equation}\label{relacionsigmas}
\sigma=\mu_A\co (\mu_A\ot h^{-1})\co (\mu_A\ot \tau \ot \mu_H)\co (h\ot \phi_A\ot \delta_{H\ot H})\co (\delta_H\ot h\ot H\ot H)\co \delta_{H\ot H}.
\end{equation}

}
\end{defin}

\begin{prop}
\label{relacionequivalencia}
Let $H$ be a cocommutative weak Hopf algebra. Then $\approx$ is an equivalence relation.

\end{prop}

\begin{dem}
Let $(\varphi_A, \sigma)$ be a crossed system. The morphism $u_{1}$
is in $Reg_{\varphi_{A}}(H,A)$ with inverse $u_{1}^{-1}=u_{1}$ and
satisfies that $u_{1}\co \eta_H=\eta_A$. Moreover, using that $(A,
\varphi_A)$ is a weak left $H$-module algebra,

\begin{itemize}

\item[ ]$\hspace{0.38cm} \mu_A\co (\mu_A\ot A)\co (u_1\ot \varphi_A\ot u_{1}^{-1})\co
(\delta_H\ot c_{H, A})\co (\delta_H\ot A)$

\item[ ]

\item[ ]$=\mu_A\co (\mu_A\ot A)\co ((\varphi_A\co (H\ot \eta_A))\ot \varphi_A\ot (\varphi_A\co (H\ot \eta_A)))\co (\delta_H\ot c_{H,A})\co (\delta_H\ot A)$

\item[ ]

\item[ ]$=\mu_A\co (\varphi_A\ot(\varphi_A\co (H\ot \eta_A)))\co (\delta_H\ot c_{H,A})\co (\delta_H\ot A)$

\item[ ]

\item[ ]$=\varphi_A,$

\end{itemize}

and we get (\ref{relacionfis}).

As far as (\ref{relacionsigmas}), using that $(A, \varphi_A)$ is a
weak left $H$-module algebra and taking into account that $\sigma$
is in $Reg_{\varphi_{A}}(H\ot H,A)$,

\begin{itemize}

\item[ ]$\hspace{0.38cm} \mu_A\co (\mu_A\ot u_{1}^{-1})\co (\mu_A\ot \sigma \ot \mu_H)\co (u_1\ot
\varphi_A\ot \delta_{H\ot H})\co (\delta_H\ot u_1\ot H\ot H)\co
\delta_{H\ot H}$

\item[ ]

\item[ ]$=\mu_A\co ((\mu_A\co (\varphi_A\co (H\ot \eta_A))\ot \varphi_A)\ot (\sigma \wedge (\varphi_A\co (\mu_H\ot \eta_A))))\co (\delta_H\ot (\varphi_A\co (H\ot \eta_A))\ot H\ot H\co \delta_{H\ot H})$

\item[ ]

\item[ ]$=\mu_A\co ((\varphi_A\co (\mu_H\ot \eta_A))\ot (\sigma \wedge \sigma^{-1}\wedge \sigma))\co \delta_{H\ot H})$

\item[ ]

\item[ ]$=\sigma\wedge\sigma^{-1}\wedge\sigma$

\item[ ]

\item[ ]$=\sigma,$

\end{itemize}

and the relation is reflexive.

In order to get that $\approx$ is symmetrical, assume that
$(\varphi_A, \sigma) \approx (\phi_A, \tau)$. Let $h$ be the
morphism in $Reg_{\varphi_{A}}(H,A)\cap Reg_{\phi_{A}}(H,A)$
satisfying (\ref{relacionfis}) and (\ref{relacionsigmas}) and such
that $h\co \eta_H=\eta_A$. Then the inverse $h^{-1}$ is in
$Reg_{\varphi_{A}}(H,A)\cap Reg_{\phi_{A}}(H,A)$ and by Proposition
\ref{f-f-1-eta} we obtain that $h^{-1}\co \eta_H=\eta_A$. Moreover,

\begin{itemize}

\item[ ]$\hspace{0.38cm} \mu_A\co (\mu_A\ot A)\co (h^{-1}\ot \varphi_A\ot h)\co
(\delta_H\ot c_{H, A})\co (\delta_H\ot A)$

\item[ ]

\item[ ]$=\mu_A\co (\mu_A\ot A)\co ((h^{-1}\wedge h)\ot \phi_A\ot (h^{-1}\wedge h))
\co (\delta_H\ot c_{H, A})\co (\delta_H\ot A)$

\item[ ]

\item[ ]$=\mu_A\co (\mu_A\ot A)\co ((\phi_A\co (H\ot \eta_A))\ot \phi_A\ot
(\phi_A\co (H\ot \eta_A)))\co (\delta_H\ot c_{H, A})\co (\delta_H\ot
A)$

\item[ ]

\item[ ]$=\mu_A\co (\phi_A\ot (\phi_A\co (H\ot \eta_A)))\co (H\ot c_{H, A})\co
(\delta_H\ot A)$

\item[ ]

\item[ ]$=\phi_A$

\end{itemize}

using that $(\varphi_A, \sigma) \approx (\psi_A, \tau)$, (e1) and
that $(A, \psi_A)$ is a weak left $H$-module algebra.

In a similar way we obtain (\ref{relacionsigmas}) and the relation is symmetrical.

Finally we show the transitivity. Assume that $(\varphi_A, \sigma)
\approx (\phi_A, \tau)$ and $(\phi_A, \tau) \approx (\chi_A,
\gamma)$ with morphisms $h$ in $Reg_{\varphi_{A}}(H,A)\cap
Reg_{\phi_{A}}(H,A)$ and $g$ in $Reg_{\phi_{A}}(H,A)\cap
Reg_{\chi_{A}}(H,A)$, respectively. Then,  the convolution product $h\wedge
g$ is in $Reg_{\varphi_{A}}(H,A)\cap Reg_{\chi_A}(H,A)$ and
$(h\wedge g)\co \eta_H=\eta_A$. Indeed:

\begin{itemize}

\item[ ]$\hspace{0.38cm} (h\wedge g)\co \eta_H$

\item[ ]

\item[ ]$=(h\wedge (g\co \Pi^{L}_{H}))\co \eta_H$

\item[ ]

\item[ ]$=(h\wedge g^{-1}\wedge g)\co \eta_H$

\item[ ]

\item[ ]$=(h\wedge h^{-1}\wedge h)\co \eta_H$

\item[ ]

\item[ ]$=h\co \eta_H$

\item[ ]

\item[ ]$=\eta_A$

\end{itemize}

In the previous calculations, the first equality follows by
(\ref{deltaPIL}), the second one by (\ref{hPiL}); in the third we
use that $g^{-1}\wedge g=h^{-1}\wedge h$; the fourth relies on (e2)
and the last one follows by the definition of $h$.

The proof for the conditions  (\ref{relacionfis}) and
(\ref{relacionsigmas}) follows a similar pattern to the well-know
in the classical case, and we leave the details to the reader.

\end{dem}

\begin{rem}\label{comparacioncasoclasico}

{\rm  We have given the detailed calculus for the above
Proposition in order to illustrate the differences when working
with weak Hopf algebras. Note that the proof  is trivial in the classical case:
if $H$ is a Hopf algebra, the relation is reflexive using the morphism
$h=\varepsilon_H\ot \eta_A$, and it is easy to get that it is symmetrical because
$h\wedge h^{-1}=h^{-1}\wedge h=\varepsilon_H\ot \eta_A$. Obviously, these equalities
are not true for weak Hopf algebras.
 }

\end{rem}

\begin{prop}
\label{sigma-sigma-1}
Let $H$ be a cocommutative weak Hopf algebra,
$(A,\varphi_{A})$ a weak left $H$-module algebra and
$\sigma \in Reg_{\varphi_{A}}(H\ot H,A)$. The following assertions hold.

\begin{itemize}

\item[(i)] $\sigma\co (\eta_{H}\ot H)=u_{1}\;\; \Leftrightarrow \;\;
\sigma^{-1}\co (\eta_{H}\ot H)=u_{1}.$

\item[(ii)] $\sigma\co (H\ot \eta_{H})=u_{1}\;\; \Leftrightarrow \;\;
\sigma^{-1}\co (H\ot \eta_{H})=u_{1}.$

\end{itemize}

\end{prop}

\begin{dem} We prove (i). The prove for (ii) is similar and we leave the details to the reader.

\begin{itemize}

\item[ ]

\item[ ]$\hspace{0.38cm} \sigma^{-1}\co (\eta_{H}\ot H)$

\item[ ]

\item[ ]$=(\sigma^{-1}\wedge u_{2})\co (\eta_{H}\ot H)    $

\item[ ]

\item[ ]$=\mu_{A}\co (\sigma^{-1}\ot (u_{1}\co \mu_{H}))\co \delta_{H\ot H}\co (\eta_{H}\ot H)  $

\item[ ]

\item[ ]$=\mu_{A}\co ((\sigma^{-1}\co c_{H,H})\ot (u_{1}\co \mu_{H}))\co (H\ot (\delta_{H}\co \eta_{H})\ot H)\co \delta_{H} $

\item[ ]

\item[ ]$=\mu_{A}\co ((\sigma^{-1}\co c_{H,H})\ot u_{1}))\co
 (((\overline{\Pi}_{H}^{L}\ot H)\co \delta_{H})\ot H)\co \delta_{H}   $

\item[ ]

\item[ ]$=\mu_{A}\co ((\sigma^{-1}\co c_{H,H}\co (H\ot \overline{\Pi}_{H}^{L} )\co \delta_{H})\ot (\sigma\co (\Pi_{H}^{L}\ot H)\co \delta_{H}))\co \delta_{H}$

\item[ ]

\item[ ]$=\mu_{A}\co ((\sigma^{-1}\co c_{H,H})\ot \sigma)\co (H\ot ((( \overline{\Pi}_{H}^{L}\co \mu_{H})\ot H)\co (H\ot c_{H,H})\co ((\delta_{H}\co \eta_{H})\ot H))\ot H)\co (H\ot \delta_{H})\co \delta_{H} $

\item[ ]

\item[ ]$=\mu_{A}\co ((\sigma^{-1}\co c_{H,H})\ot \sigma)\co (H\ot ((H\ot (((\varepsilon_{H}\co \mu_{H})\ot H)\co (H\ot c_{H,H})\co (\delta_{H}\ot H)))$

\item[ ]

\item[ ]$\hspace{0.38cm}\co ((\delta_{H}\co \eta_{H}) \ot (\delta_{H}\co \eta_{H})\ot H))\ot H)\co (H\ot \delta_{H})\co \delta_{H} $

\item[ ]

\item[ ]$=\mu_{A}\co ((\sigma^{-1}\co c_{H,H})\ot \sigma)\co (H\ot ((H\ot (((\varepsilon_{H}\co \mu_{H})\ot H)\co (H\ot c_{H,H})\co ((c_{H,H}\co \delta_{H})\ot H)))$

\item[ ]

\item[ ]$\hspace{0.38cm}\co ((\delta_{H}\co \eta_{H}) \ot (\delta_{H}\co \eta_{H})\ot H))\ot H)\co (H\ot (c_{H,H}\co \delta_{H}))\co \delta_{H} $

\item[ ]

\item[ ]$=((\mu_{A}\co (\sigma^{-1}\ot \sigma))\ot (\varepsilon_{H}\co \mu_{H}))\co
(H\ot c_{H,H}\ot c_{H,H}\ot H)\co (\delta_{H}\ot c_{H,H}\ot H)\co ((\delta_{H}\co \eta_{H})\ot ((\delta_{H}\ot H)\co \delta_{H})   $

\item[ ]

\item[ ]$=((\sigma^{-1}\wedge \sigma))\ot (\varepsilon_{H}\co \mu_{H}))\co \delta_{H\ot H}\co (\eta_{H}\ot H) $

\item[ ]

\item[ ]$=((u_{1}\co \mu_{H}))\ot (\varepsilon_{H}\co \mu_{H}))\co \delta_{H\ot H}\co (\eta_{H}\ot H)   $

\item[ ]

\item[ ]$=(u_{1}\ot \varepsilon_{H})\co \delta_{H}\co \mu_{H}\co  (\eta_{H}\ot H)  $

\item[ ]

\item[ ]$=u_{1},    $

\end{itemize}

 where the first and the twelfth equalities follow by the properties of $\sigma$, the second one by the definition of $u_{2}$, the third, tenth and eleventh ones by the naturality of $c$, the fourth one by (\ref{deltaPILbarra}),  the fifth one by the coassociativity of $\delta_{H}$ and  (\ref{sigma-PiL-l}), the sixth one by the coassociativity of $\delta_{H}$ and  (\ref{deltaPIL}), the seventh one by the definition of $\overline{\Pi}_{H}^{L}$ and the associativity of $\mu_{H}$, the eighth one by (a3) of Definition \ref{wha}, the ninth one by the cocommutativity of $H$, the thirteenth one by (a1) of Definition \ref{wha} and the last one by the unit-counit properties.

The proof for the converse is the same  changing $\sigma$ by $\sigma^{-1}$.

\end{dem}

\begin{apart}
\label{weak-crossed-products-exposition}
{\rm  The equalities (g1), (g2) and (g3) of Definition
\ref{stmacruzado}  have a clear meaning in the theory of weak
crossed products introduced in \cite{nmra4} and \cite{mra-preunit}.
The full details can be found in Section 2 of \cite{NikaRamon6}. In
this point we give a brief resume adapted to our setting, i.e.,
with some changes in the notation.

Let $H$ be a  weak Hopf algebra, $(A,\varphi_{A})$  a weak left
$H$-module algebra and $\sigma:H\ot H\rightarrow A$ a morphism. We
define the morphisms
$$\psi_{H}^{A}:H\ot A\rightarrow A\ot H, \;\;\;\sigma_{H}^{A}:H\ot H\rightarrow A\ot
H, $$ by
\begin{equation}
\label{psiAH} \psi_{H}^{A}=(\varphi_{A}\ot H)\co (H\ot c_{H,A})\co
(\delta_{H}\ot A)
\end{equation}
and
\begin{equation}
\label{sigmaAH} \sigma_{H}^{A}=(\sigma\ot \mu_{H})\co \delta_{H\ot H}.
\end{equation}

Then, $\psi_{H}^{A}$ satisfies
\begin{equation}\label{wmeas-wcp}
(\mu_A\ot H)\co (A\ot \psi_{H}^{A})\co (\psi_{H}^{A}\ot A) =
\psi_{H}^{A}\co (H\ot \mu_A).
\end{equation}

 As a consequence of (\ref{wmeas-wcp}), the morphism $\nabla_{A\ot H}:A\ot H\rightarrow
A\ot H$ defined by
\begin{equation}\label{idem-wcp}
\nabla_{A\ot H} = (\mu_A\ot H)\co(A\ot \psi_{H}^{A})\co (A\ot H\ot
\eta_A)
\end{equation}
is  idempotent. Moreover, $\nabla_{A\otimes H}$ satisfies that
$$\nabla_{A\otimes H}\co (\mu_A\ot H) = (\mu_A\ot H)\co
(A\ot \nabla_{A\otimes H}),$$ that is, $\nabla_{A\otimes H}$ is a left
$A$-module morphism (see Lemma 3.1 of \cite{mra-preunit}) for the
regular action  $\varphi_{A\ot H}=\mu_{A}\ot H$. With $A\times H$,
$i_{A\ot H}:A\times H\rightarrow A\ot H$ and $p_{A\otimes H}:A\ot
H\rightarrow A\times H$ we denote the object, the injection and the
projection associated to the factorization of $\nabla_{A\otimes H}$.
Moreover, if $\psi_{H}^{A}$ satisfies (\ref{wmeas-wcp}), the
following identities hold
\begin{equation}\label{fi-nab}
(\mu_{A}\ot H)\co (A\ot \psi_{H}^{A})\co (\nabla_{A\otimes H}\ot A)=
(\mu_{A}\ot H)\co (A\ot \psi_{H}^{A})=\nabla_{A\otimes H}\co(\mu_{A}\ot
H)\co (A\ot \psi_{H}^{A})
\end{equation}
\begin{equation}\label{fi-nab-2}
 \nabla_{A\otimes H}\co \psi_{H}^{A}= \psi_{H}^{A}
\end{equation}
and, by the naturality of $c$, we have
\begin{equation}\label{nabla-u1}
\nabla_{A\otimes H}=((\mu_{A}\co (A\ot u_{1}))\ot H)\co (A\ot \delta_{H})
\end{equation}
and this implies that $\nabla_{A\otimes H}$ is a morphism of right $H$-comodules for
$\rho_{A\ot H}=A\ot \delta_{H}$.

Also, in Propositions 2.7 and 2.8 of \cite{NikaRamon6} we can find the proof of the
equalities:
\begin{equation}
\label{nabla-nabla} \nabla_{A\otimes H}=((\mu_{A}\co (A\ot u_{1}))\ot
H)\co (A\ot \delta_{H}),
\end{equation}
\begin{equation}
\label{nabla-fi} \mu_{A}\co (u_{1}\ot \varphi_{A})\co (\delta_{H}\ot
A)=\varphi_{A},
\end{equation}
\begin{equation}
\label{nabla-fiAH} (\mu_{A}\ot H)\co (u_{1}\ot \psi_{H}^{A})\co
(\delta_{H}\ot A)=\psi_{H}^{A},
\end{equation}
\begin{equation}
\label{eta-psi-varep} (A\ot \varepsilon_{H})\co \psi_{H}^{A}\co
(H\ot \eta_{A})=u_{1},
\end{equation}
\begin{equation}
\label{eta-psi-complex} (\mu_{A}\ot H)\co (u_{1}\ot c_{H,A})\co
(\delta_{H}\ot A)=(\mu_{A}\ot H)\co (A\ot c_{H,A})\co
((\psi_{H}^{A}\co (H\ot \eta_{A}))\ot A),
\end{equation}
\begin{equation}
\label{nabla-varep} (A\ot \varepsilon_{H})\co \nabla_{A\otimes H}=\mu_{A}\co (A\ot u_{1}).
\end{equation}
\begin{equation}
\label{nabla-delta} (A\ot \delta_{H})\co \nabla_{A\otimes H}=(\nabla_{A\otimes H}\ot H)\co (H\ot \delta_{H}).
\end{equation}
\begin{equation}
\label{delta-sigmaHA} (A\ot \delta_{H})\co
\sigma_{H}^{A}=(\sigma_{H}^{A}\ot \mu_{H})\co \delta_{H\ot H}.
\end{equation}

Moreover, if  $H$ is a cocommutative  and $\sigma\in
Reg_{\varphi_{A}}(H\ot H,A)$,  the morphism $\sigma_{H}^{A}$
 satisfies the following identities:
\begin{equation}
\label{sigmaOmegaH} \sigma_{H}^{A}\co \Omega_{H}^{2}=\sigma_{H}^{A},
\end{equation}
\begin{equation}
\label{Nablasigma}\nabla_{A\otimes H}\co\sigma_{H}^{A}=\sigma_{H}^{A},
\end{equation}
\begin{equation}
\label{varepsilonsigma}(A\ot \varepsilon_{H})\co
\sigma_{H}^{A}=\sigma.
\end{equation}

If we consider the quadruple ${\Bbb A}_{H}=(A,H, \psi_{H}^{A},
\sigma_{H}^{A})$, where  $H$ a cocommutative weak Hopf algebra,
$(A,\varphi_{A})$ is a weak left $H$-module algebra and $\sigma\in
Reg_{\varphi_{A}}(H\ot H,A)$, we say that  ${\Bbb A}_{H}$ satisfies
the twisted condition if
\begin{equation}\label{twis-wcp}
(\mu_A\ot H)\co (A\ot \psi_{H}^{A})\co (\sigma_{H}^{A}\ot A) =
(\mu_A\ot H)\co (A\ot \sigma_{H}^{A})\co (\psi_{H}^{A}\ot H)\co
(H\ot \psi_{H}^{A})
\end{equation}
and   the  cocycle condition holds if
\begin{equation}\label{cocy2-wcp}
(\mu_A\ot H)\co (A\ot \sigma_{H}^{A}) \co (\sigma_{H}^{A}\ot H) =
(\mu_A\ot H)\co (A\ot \sigma_{H}^{A})\co (\psi_{H}^{A}\ot H)\co
(H\ot\sigma_{H}^{A}).
\end{equation}

Note that, if ${\Bbb A}_{H}=(A, H, \psi_{H}^{A}, \sigma_{H}^{A})$
satisfies the twisted condition, in Proposition 3.4 of
\cite{mra-preunit}, we proved that the following equalities hold:
\begin{equation}\label{c1}
(\mu_A\otimes H)\circ (A\otimes \sigma_{H}^{A})\circ
(\psi_{H}^{A}\otimes H)\circ (H\otimes \nabla_{A\otimes H}) =
\nabla_{A\otimes H}\circ (\mu_A\otimes H)\circ (A\otimes
\sigma_{H}^{A})\circ (\psi_{H}^{A}\otimes H),
\end{equation}
\begin{equation}\label{aw}
\nabla_{A\otimes H}\circ (\mu_A\otimes H)\circ
(A\otimes\sigma_{H}^{A})\circ (\nabla_{A\otimes H}\otimes H) =
\nabla_{A\otimes H}\circ (\mu_A\otimes H)\circ
(A\otimes\sigma_{H}^{A}).
\end{equation}

Then, by (\ref{Nablasigma}), we obtain
\begin{equation}\label{c11}
(\mu_A\otimes H)\circ (A\otimes \sigma_{H}^{A})\circ
(\psi_{H}^{A}\otimes H)\circ (H\otimes \nabla_{A\otimes H}) =
 (\mu_A\otimes H)\circ (A\otimes
\sigma_{H}^{A})\circ (\psi_{H}^{A}\otimes H),
\end{equation}
\begin{equation}\label{aw1}
 (\mu_A\otimes H)\circ
(A\otimes\sigma_{H}^{A})\circ (\nabla_{A\otimes H}\otimes H) =
(\mu_A\otimes H)\circ (A\otimes\sigma_{H}^{A}).
\end{equation}

For the  product defined by
\begin{equation}\label{prod-todo-wcp}
\mu_{A\ot_{\sigma}  H} = (\mu_A\ot H)\co (\mu_A\ot \sigma_{H}^{A})\co (A\ot
\psi_{H}^{A}\ot H),
\end{equation}
if the twisted and the cocycle conditions hold, we obtain that it is
associative and normalized with respect to $\nabla_{A\otimes H}$ (i.e.
$\nabla_{A\otimes H}\co \mu_{A\ot_{\sigma} H}=\mu_{A\ot_{\sigma} H}=
\mu_{A\ot_{\sigma}H}\co (\nabla_{A\otimes H}\ot \nabla_{A\otimes H}$)).
Due to the normality condition,
\begin{equation}
\label{prod-wcp} \mu_{A\times_{\sigma} H} =
p_{A\otimes H}\co\mu_{A\ot H}\co
(i_{A\otimes H}\ot i_{A\otimes H}),
\end{equation}
is associative as
well (Proposition 3.7 of \cite{mra-preunit}). Hence if ${\Bbb
A}_{H}=(A,H, \psi_{H}^{A}, \sigma_{H}^{A})$  satisfies
(\ref{twis-wcp}) and (\ref{cocy2-wcp}) we say that $A\ot_{\sigma}H=(A\ot H,
\mu_{A\ot_{\sigma} H})$ is a weak crossed product.

If ${\Bbb A}_{H}$ satisfies
\begin{equation}\label{pre1-wcp}
    (\mu_A\otimes H)\circ (A\otimes \sigma_{H}^{A})\circ
    (\psi_{H}^{A}\otimes H)\circ (H\otimes \nu) =
    \nabla_{A\otimes H}\circ
    (\eta_A\otimes H),
    \end{equation}
\begin{equation}\label{pre2-wcp}
    (\mu_A\otimes H)\circ (A\otimes \sigma_{H}^{A})\circ
    (\nu\otimes H) = \nabla_{A\otimes H}\circ (\eta_A\otimes H),
    \end{equation}
\begin{equation}\label{pre3-wcp}
(\mu_A\otimes H)\circ (A\otimes \psi_{H}^{A})\circ (\nu\otimes A) =
\beta_{\nu},
\end{equation}
for $\nu:K\rightarrow A\ot H$ and \begin{equation} \label{beta-nu}
\beta_{\nu}=(\mu_A\otimes
H)\circ (A\otimes \nu):A\rightarrow A\otimes H,\;
\end{equation}
by Theorem 2.2 of \cite{NikaRamon6},  we obtain that $\nu$ is a
preunit (see (73) of \cite{NikaRamon6} for the definition in this
setting) for the product $\mu_{A\ot_{\sigma} H}$ defined in
(\ref{prod-todo-wcp}). Therefore $A\times H$ is an algebra with
the product defined in (\ref{prod-wcp}) and unit $\eta_{A\times_{\sigma}
H}=p_{A\otimes H}\circ\nu$. In what follows we denote this algebra by
$A\times_{\sigma}H$.

Following Definition 2.11 of \cite{NikaRamon6}, we say that $\sigma$ satisfies the
twisted condition if
\begin{equation}
\label{twisted-sigma} \mu_{A}\co ((\varphi_{A}\co (H\ot
\varphi_{A}))\ot A)\co (H\ot H\ot c_{A,A})\co (((H\ot H\ot
\sigma)\co \delta_{H\ot H}) \ot A)=\mu_{A}\co (A\ot \varphi_{A})\co
(\sigma_{H}^{A}\ot A),
\end{equation}
and if
\begin{equation}
\label{2-cocycle-sigma}
\mu_{A}\co (A\ot \sigma)\co (\sigma_{H}^{A}\ot H)=
\mu_{A}\co (A\ot \sigma)\co (\psi_{H}^{A}\ot H)\co (H\ot
\sigma_{H}^{A})
\end{equation}
holds, we will say that $\sigma$ satisfies the
2-cocycle condition.

It is a trivial calculus to prove that (g1) of Definition \ref{stmacruzado} is equivalent to
(\ref{twisted-sigma})
and also (g2) is equivalent to (\ref{2-cocycle-sigma}). Moreover, by Theorem 2.13 of
\cite{NikaRamon6}, we know that $\sigma$ satisfies the
twisted condition (\ref{twisted-sigma}) if and only if ${\Bbb
A}_{H}$ satisfies the twisted condition (\ref{twis-wcp}), and, by Theorem 2.14 \cite{NikaRamon6}, we also
know that $\sigma$ satisfies the
2-cocycle condition (\ref{2-cocycle-sigma}) if and only if ${\Bbb
A}_{H}$ satisfies the cocycle condition (\ref{cocy2-wcp}).

On the other hand, (g3) of Definition \ref{stmacruzado} is exactly the normal condition
introduced in Definition 2.16 of \cite{NikaRamon6} and also, by Theorem 2.18 and Corollary 2.19 of \cite{NikaRamon6},
we have that $\nu=\nabla_{A\otimes H}\co (\eta_{A}\ot\eta_{H})$ is a preunit for the weak crossed
product associated to ${\Bbb A}_{H}$ if and only if (g3) holds.

Therefore, if $(\varphi_A, \sigma)$ is a crossed system for $H$ over $A$, we
have that $A\ot_{\sigma} H$ is a weak
crossed product, with preunit $\nu=\nabla_{A\otimes H}\co
(\eta_{A}\ot\eta_{H})$.
Conversely, if the pair $(\varphi_A, \sigma)$ is such that $A\ot_{\sigma} H$ is a weak
crossed product, with preunit $\nu=\nabla_{A\otimes H}\co
(\eta_{A}\ot\eta_{H})$ and normalized with respect to $\nabla_{A\otimes H}$, we obtain that
$(\varphi_A, \sigma)$ is a crossed system for $H$ over $A$ (see Corollary 2.20 of
\cite{NikaRamon6}).

Note that, if $(\varphi_A, \sigma)$ and
$(\phi_A, \tau)$ are  equivalent crossed systems for $H$ over $A$, the equality
 $$u_{1}=\varphi_A\co (H\ot \eta_A)=\phi_A\co (H\ot \eta_A)=u_{1}^{\prime}$$ holds and then the associated
idempotent morphisms are the same because, by (\ref{nabla-u1}),
\begin{equation}
\label{nablaprima-nabla}
\nabla_{A\ot H}^{\prime}=((\mu_{A}\co (A\ot u_{1}^{\prime}))\ot H)\co (A\ot \delta_{H})=
((\mu_{A}\co (A\ot u_{1}))\ot H)\co (A\ot \delta_{H})=\nabla_{A\ot H}.
\end{equation}
Therefore, they define an algebra structure over the same object $A\times H$.

}
\end{apart}

In the following result we characterize the crossed products where $\varphi_A$ is an $H$-module structure.

\begin{teo}
\label{cocycle-center} Let $H$ be a cocommutative weak Hopf algebra,
$(A,\varphi_{A})$  a weak left $H$-module algebra and $\sigma\in
Reg_{\varphi_{A}}(H\ot H,A)$ satisfying the twisted condition
(\ref{twisted-sigma}). The following assertions are equivalent:
\begin{itemize}
\item[(i)] $(A, \varphi_{A})$ is a left $H$-module algebra.
\item[(ii)] The morphism $\sigma$ factors through the center of $A$.
\end{itemize}
\end{teo}

\begin{dem}
Let $(A, \varphi_{A})$ be a left $H$-module algebra. We define
$\gamma_{\sigma}:A\ot H\ot H\rightarrow A$ as
$$\gamma_{\sigma}=\mu_{A}\co ((\mu_{A}\co c_{A,A})\ot A)\co (A\ot ((\sigma\ot \sigma^{-1})\co \delta_{H\ot H})).$$
Then,
\begin{equation}
\label{gamma-u2} \gamma_{\sigma}=\mu_{A}\co (A\ot u_{2})
\end{equation}
because

\begin{itemize}

\item[ ]$\hspace{0.38cm}\gamma_{\sigma}$

\item[ ]

\item [ ]$= \mu_{A}\co ((\mu_{A}\co c_{A,A})\ot A)\co (A\ot (((\sigma\wedge
u_{2})\ot \sigma^{-1})\co \delta_{H\ot H})) $

\item[ ]

\item [ ]$=\mu_{A}\co ((\mu_{A}\co  (A\ot \mu_{A})\co (A\ot (c_{A,A}\co (A\ot u_{1}))))\ot A)\co
(c_{A,A}\ot H\ot A)\co (A\ot \sigma_{H}^{A}\ot \sigma^{-1})\co (A\ot
\delta_{H\ot H}) $

\item[ ]

\item [ ]$=\mu_{A}\co ((\mu_{A}\co (A\ot (\varphi_{A}\co
(\Pi_{H}^{L}\ot A)\co c_{A,H}))\co (c_{A,A}\ot H))\ot A)\co (A\ot
\sigma_{H}^{A}\ot \sigma^{-1})\co (A\ot \delta_{H\ot H})    $

\item[ ]

\item [ ]$=\mu_{A}\co ((\mu_{A}\co (A\ot (A\ot (\varphi_{A}\co c_{A,H}\co
(A\ot (\mu_{H}\co (H\ot \lambda_{H})\co (\mu_{H}\ot \mu_{H})\co
\delta_{H\ot H}))))))\ot A)$

\item[ ]

\item[ ]$\hspace{0.38cm}\co (c_{A,A}\ot H\ot H\ot A)\co (A\ot \sigma\ot H\ot H\ot A)\co (A\ot
\delta_{H\ot H}\ot \sigma^{-1})\co (A\ot  \delta_{H\ot H})$

\item[ ]

\item [ ]$=\mu_{A}\co ((\mu_{A}\co (A\ot \varphi_{A})\co (\sigma_{H}^{A}\ot A)\co
(H\ot H\ot (\varphi_{A}\co c_{A,H}\co (A\ot \lambda_{H})))\co (H\ot
c_{A,H}\ot \mu_{H})$

\item[ ]

\item[ ]$\hspace{0.38cm}\co(c_{A,H}\ot H\ot H\ot H)\co (A\ot\delta_{H\ot H}))\ot
\sigma^{-1})\co (A\ot  \delta_{H\ot H})$

\item[ ]

\item [ ]$= \mu_{A}\co ((\mu_{A}\co ((\varphi_{A}\co (H\ot \varphi_{A}))\ot A)\co
(H\ot H\ot c_{A,A})\co (((H\ot H\ot \sigma)\co \delta_{H\ot H}) \ot
A)$

\item[ ]

\item[ ]$\hspace{0.38cm}\co(H\ot H\ot (\varphi_{A}\co c_{A,H}\co (A\ot \lambda_{H})))\co
(H\ot c_{A,H}\ot \mu_{H})\co (c_{A,H}\ot H\ot H\ot H)\co
(A\ot\delta_{H\ot H}))\ot \sigma^{-1})\co $

\item[ ]

\item[ ]$\hspace{0.38cm}(A\ot  \delta_{H\ot H})$

\item[ ]

\item [ ]$=\mu_{A}\co (((\varphi_A\co (\mu_H\ot A)\co (H\ot \mu_H\ot A))\co (H\ot H\ot c_{A,H}))\ot \mu_{A})\co (H\ot c_{A,H}
\ot (c_{A,H}\co (\sigma\ot (\lambda_{H}\co \mu_{H}))$

\item[ ]

\item[ ]$\hspace{0.38cm}\co (A\ot
c_{H,H}\ot H)\co ((c_{H,H}\co \delta_{H})\ot (c_{H,H}\co
\delta_{H})))\ot A) \co (c_{A,H}\ot H\ot H\ot H\ot \sigma^{-1})$

\item[ ]

\item[ ]$\hspace{0.38cm}\co
(H\ot \delta_{H\ot H}\ot H\ot H)\co (A\ot  \delta_{H\ot H})$

\item[ ]

\item [ ]$= \mu_{A}\co (((\varphi_{A}\co (\mu_H\ot A))\co (H\ot  c_{A,H}))\ot A)\co (c_{A,H}\ot \lambda_{H}\ot A)\co
(A\ot ((\mu_{H}\ot \mu_{H})\co \delta_{H\ot H})\ot (\sigma\wedge
\sigma^{-1}))$

\item[ ]

\item[ ]$\hspace{0.38cm}\co  (A\ot  \delta_{H\ot H}) $

\item[ ]

\item [ ]$=\mu_{A}\co ((\varphi_{A}\co
c_{A,H}\co (A\ot (\Pi_{H}^{L}\co \mu_{H})))\ot u_{2})\co (A\ot
\delta_{H\ot H})   $

\item[ ]

\item [ ]$=\mu_{A}\co ((\varphi_{A}\co (\overline{\Pi}_{H}^{L}\ot A)\co c_{A,H})\ot u_{1})\co
 (A\ot (\delta_{H}\co \mu_{H})) $

\item[ ]

\item [ ]$=\mu_{A}\co (( \mu_{A}\co c_{A,A}\co (u_{1}\ot A)\co c_{A,H})\ot u_{1})\co
(A\ot (\delta_{H}\co \mu_{H}))  $

\item[ ]

\item [ ]$=\mu_{A}\co (A\ot (u_{1}\wedge u_{1}))\co (A\ot \mu_{H})$

\item[ ]

\item [ ]$=\mu_{A}\co (A\ot u_{2}) $

\item[ ]

\end{itemize}

where the first and the ninth equalities follow by
(f1)-(f3), the second one by (\ref{u1-u2}) and the naturality of $c$, the third one by
the associativity of $\mu_{A}$ and by (d4) of Definition
\ref{weak-H-mod}, the fourth one by the definition of $\Pi_{H}^{L}$,
the fifth and the eight  ones by the naturality of $c$, the
coassociativity of $\delta_{H}$ and (d3-1) of Definition
\ref{weak-H-mod}. The sixth one is a consequence of the twisted
condition, the seventh follows by  the naturality of $c$ and the
coassociativity and cocommutativity of $\delta_{H}$, the tenth one
by the cocommutativity of $\delta_{H}$, the eleventh one by (d5) of
Definition \ref{weak-H-mod}, the twelfth one by the naturality of
$c$ and finally, the last one, by (d2)  of Definition
\ref{weak-H-mod}.

Therefore, $\sigma$ factors through the center of $A$ because

\begin{itemize}

\item[ ]

\item[ ]$\hspace{0.38cm}\mu_{A}\co (A\ot \sigma) $

\item[ ]

\item [ ]$=\mu_{A}\co (A\ot (u_{2}\wedge\sigma))  $

\item[ ]

\item [ ]$= \mu_{A}\co  ((\mu_{A}\co (A\ot u_{2}))\ot \sigma)\co (A\ot \delta_{H\ot H})   $

\item[ ]

\item [ ]$= \mu_{A}\co  (\gamma_{\sigma}\ot \sigma)\co (A\ot \delta_{H\ot H}) $

\item[ ]

\item [ ]$=\mu_{A}\co  ((\mu_{A}\co ((\mu_{A}\co c_{A,A})\ot A)\co (A\ot ((\sigma\ot \sigma^{-1})\co \delta_{H\ot H})))
\ot \sigma)\co (A\ot \delta_{H\ot H}) $

\item[ ]

\item [ ]$= \mu_{A}\co ((\mu_{A}\co c_{A,A})\ot A)\co (A\ot
((\sigma\ot (\sigma^{-1}\wedge \sigma))\co \delta_{H\ot H}))$

\item[ ]

\item [ ]$= \mu_{A}\co ((\mu_{A}\co c_{A,A})\ot A)\co (A\ot
((\sigma\ot u_{2})\co \delta_{H\ot H})) $

\item[ ]

\item [ ]$=\mu_{A}\co ( A\ot (\varphi_{A}\co c_{A,H}\co (A\ot \overline{\Pi}_{H}^{L})))\co (c_{A,A}\ot H)
\co (A\ot \sigma_{H}^{A})  $

\item[ ]

\item [ ]$=\mu_{A}\co ( A\ot (\varphi_{A}\co (\Pi_{H}^{L}\ot A) \co c_{A,H}))\co (c_{A,A}\ot H)
\co (A\ot \sigma_{H}^{A})    $

\item[ ]

\item [ ]$=\mu_{A}\co ( A\ot (\mu_{A}\co (u_{1}\ot A)\co c_{A,H}))\co (c_{A,A}\ot H)
\co (A\ot \sigma_{H}^{A})  $

\item[ ]

\item [ ]$=\mu_{A}\co ((\sigma\wedge u_{2})\ot A)\co (H\ot c_{A,H})\co (c_{A,H}\ot H)  $

\item[ ]

\item [ ]$=\mu_{A}\co c_{A,A}\co (A\ot \sigma) $

\item[ ]

\end{itemize}

where the first, the sixth and the eleventh equalities follow by (f1)-(f3), the second one by the associativity of
$\mu_{A}$, the third one by (\ref{gamma-u2}), the fourth one by the
definition of $\gamma_{\sigma}$, the fifth one by the naturality of
$c$, the coassociativity of $\delta_{H}$ and the associativity of
$\mu_{A}$, the seventh and the ninth ones by (d5) of Definition \ref{weak-H-mod},
the eight one by the cocommutativity of $\delta_{H}$,  the tenth one by the
naturality of $c$ and the associativity of $\mu_{A}$.

Conversely, assume that the morphism $\sigma$ factors through the center
of $A$. Then, $(A,\varphi_{A})$ is a left $H$-module algebra because

\begin{itemize}

\item[ ]

\item[ ]$\hspace{0.38cm}\varphi_{A}\co (H\ot \varphi_{A}) $

\item[ ]

\item [ ]$= \varphi_{A}\co
(H\ot (\mu_{A}\co (\varphi_{A}\ot \varphi_{A})\co (H\ot c_{H,A}\ot
A)\co (\delta_{H}\ot \eta_{A}\ot A)))  $

\item[ ]

\item [ ]$=\mu_{A}\co (\varphi_{A}\ot \varphi_{A})\co (H\ot c_{H,A}\ot A)\co
(\delta_{H}\ot ((\varphi_{A}\ot \varphi_{A})\co (H\ot c_{H,A}\ot
A)\co (\delta_{H}\ot \eta_{A}\ot A)))  $

\item[ ]

\item [ ]$= \mu_{A}\co (u_{2}\ot (\varphi_{A}\co (H\ot \varphi_{A})))\co (\delta_{H\ot H}\ot A) $

\item[ ]

\item [ ]$= \mu_{A}\co ((\sigma^{-1}\wedge \sigma)\ot (\varphi_{A}\co (H\ot \varphi_{A})))\co (\delta_{H\ot H}\ot A)$

\item[ ]

\item [ ]$=\mu_{A}\co ( \sigma^{-1} \ot
(\mu_{A}\co (\sigma \ot (\varphi_{A}\co (H\ot \varphi_{A})))\co (\delta_{H\ot H}\ot
A)))\co (\delta_{H\ot H}\ot A) $

\item[ ]

\item [ ]$=\mu_{A}\co ( \sigma^{-1} \ot
(\mu_{A}\co c_{A,A}\co (\sigma \ot (\varphi_{A}\co (H\ot \varphi_{A})))\co
(\delta_{H\ot H}\ot A)))\co (\delta_{H\ot H}\ot A)  $

\item[ ]

\item [ ]$=\mu_{A}\co ( \sigma^{-1} \ot
(\mu_{A}\co ((\varphi_{A}\co (H\ot \varphi_{A}))\ot A)\co (H\ot H\ot
c_{A,A})\co (((H\ot H\ot \sigma)\co \delta_{H\ot H}) \ot A))\co
(\delta_{H\ot H}\ot A)  $

\item[ ]

\item [ ]$=\mu_{A}\co ( \sigma^{-1} \ot
(\mu_{A}\co (A\ot \varphi_{A})\co (\sigma_{H}^{A}\ot A))\co
(\delta_{H\ot H}\ot A)   $

\item[ ]

\item [ ]$= \mu_{A}\co ( (\sigma^{-1}\wedge \sigma) \ot (\varphi_{A}\co (\mu_{H}\ot A)))\co
(\delta_{H\ot H}\ot A)   $

\item[ ]

\item [ ]$= \mu_{A}\co (u_{2}\ot (\varphi_{A}\co (\mu_{H}\ot A)))\co
( \delta_{H\ot H}\ot A)$

\item[ ]

\item [ ]$=\mu_{A}\co (\varphi_{A}\ot \varphi_{A})\co (H\ot c_{H,A}\ot A)\co
((\delta_{H}\co \mu_{H})\ot \eta_{A}\ot A))) $

\item[ ]

\item [ ]$=\varphi_{A}\co (\mu_{H}\ot A)  $

\item[ ]

\end{itemize}

where the first, the second and the twelfth equalities follows by
(d2) of Definition \ref{weak-H-mod}, the third one by the naturality
of $c$ and (d3) of Definition \ref{weak-H-mod}, the fourth  and
the tenth ones by (f1)-(f3), the fifth and the
ninth  ones by the naturality of $c$ and the coassociativity of
$\delta_{H}$, the sixth one by the cocommutativity of $\delta_{H}$
an the factorization of $\sigma$ through  $Z(A)$, the
seventh one by the naturality of $c$ and the cocommutativity of $H$, the eight one by the the
twisted condition (\ref{twisted-sigma}) and, finally, the eleventh
one by (a1) of Definition \ref{wha}.

\end{dem}

\begin{cor}
\label{sigmayu2}
Let $H$ be a cocommutative weak Hopf algebra,
$(A,\varphi_{A})$  a weak left $H$-module algebra. The following assertions are equivalent:
\begin{itemize}
\item[(i)] $(A, \varphi_{A})$ is a left $H$-module algebra.
\item[(ii)] $(\varphi_{A}, u_2)$ is a crossed system for $H$ over $A$.
\end{itemize}

\end{cor}

\begin{dem}
$(i)\Rightarrow (ii)$ Trivially, the morphism $u_2$ is in $Reg_{\varphi_{A}}(H\ot H,A)$ and satisfies (g2) and (g3). In order to get (g1), using the definition of $u_2$, that $H$ is a weak Hopf algebra, and that $(A, \varphi_{A})$ is a left $H$-module algebra,

\begin{itemize}

\item[ ]$\hspace{0.38cm} \mu_A\co (A\ot \varphi_A)\co (u_2\ot \mu_H\ot A)\co (\delta_{H\ot H}\ot A)$

\item[ ]

\item[ ]$=\mu_A\co ((\varphi_A\co (H\ot \eta_A))\ot \varphi_A)\co ((\delta_{H}\co \mu_H)\ot A)$

\item[ ]

\item[ ]$=\varphi_A\co (\mu_H\ot A)$

\item[ ]

\item[ ]$=\mu_A\co (\varphi_A\ot A)\co (H\ot c_{A,A})\co (H\ot (\varphi_A\co (H\ot \eta_A))\ot A)\co ((\delta_{H}\co \mu_H)\ot A)$

\item[ ]

\item[ ]$=\mu_A\co ((\varphi_A\co (H\ot \varphi_A))\ot A)\co (H\ot H\ot c_{A,A})\co (H\ot H\ot u_2\ot A)\co (\delta_{H\ot H}\ot A).$

\end{itemize}

As far as $(ii)\Rightarrow (i)$, note that by cocommutativity of $H$ we have that $\Pi_{H}^{L}=\overline{\Pi}_{H}^{L}$ and $u_2$ factors through the center of $A$. Applying Theorem \ref{cocycle-center} we get that $(A, \varphi_{A})$ is a left $H$-module algebra.

\end{dem}

\begin{rem}
{\rm  In the conditions of Corollary \ref{sigmayu2}, if $(A,\varphi_{A})$ is a left $H$-module algebra, we have that for the
crossed system
$(\varphi_{A}, u_2)$
\begin{equation}
\label{smash-condition}
\sigma_{H}^{A}=(u_{1}\ot H)\co \delta_{H}\co \mu_{H}.
\end{equation}
Then, the associated crossed product defined in (\ref{prod-todo-wcp}) is
\begin{equation}
\label{productosmash}
\mu_{A\ot_{u_{2}} H}=\nabla_{A\otimes H}\co (\mu_A\ot \mu_H)\co (A\ot \psi_{H}^{A} \ot H).
\end{equation}
and therefore
\begin{equation}
\label{productosmashaspa}
\mu_{A\times_{u_{2}} H}=p_{A\otimes H}\co (\mu_A\ot \mu_H)\co (A\ot \psi_{H}^{A} \ot H)\co (i_{A\otimes H}\ot i_{A\otimes H}).
\end{equation}
In this case we say that the weak crossed product is smash.

On the other hand, if  for a weak left $H$-module algebra the equality
\begin{equation}
\label{productotwistedaspa}
\varphi_A=\varphi_A\co (\Pi_{H}^{L}\ot A),
\end{equation}
holds,  by (d4) of Definition \ref{weak-H-mod}, we obtain that

\begin{itemize}

\item[ ]$\hspace{0.38cm} \mu_{A\times_{\sigma} H}$

\item[ ]

\item[ ]$=p_{A\otimes H}\co (\mu_A\ot H)\co (\mu_A\ot \sigma_{H}^{A})\co (A\ot c_{H,A}\ot H)\co ((\nabla_{A\otimes H}\co i_{A\otimes H})\ot i_{A\otimes H})$

\item[ ]

\item[ ]$=p_{A\otimes H}\co (\mu_A\ot H)\co (\mu_A\ot \sigma_{H}^{A})\co (A\ot c_{H,A}\ot H)\co (i_{A\otimes H}\ot i_{A\otimes H})$.

\end{itemize}

In this case the weak crossed product is called twisted.

}
\end{rem}

\begin{prop}
\label{H-comod-alg} Let $H$ be a cocommutative weak Hopf algebra and
 $(\varphi_A, \sigma)$  a crossed system for $H$ over $A$. Then, the algebra
 $A\times_{\sigma} H$ is
a right $H$-comodule algebra for the coaction
$$\rho_{A\times_{\sigma} H}=
(p_{A\otimes H}\ot H)\co (A\ot \delta_{H})\co i_{A\otimes H}.$$

Moreover, $(A\times_{\sigma} H)_{H}=A$.

\end{prop}

\begin{dem} By Proposition 3.2 of \cite{NikaRamon6}, we obtain
that $A\times_{\sigma} H$ is
a right $H$-comodule algebra for the coaction
$$\rho_{A\times_{\sigma} H}=
(p_{A\otimes H}\ot H)\co (A\ot \delta_{H})\co i_{A\otimes H}.$$

Moreover,
$$
\setlength{\unitlength}{3mm}
\begin{picture}(30,4)
\put(2,2){\vector(1,0){6}}
\put(14,2.5){\vector(1,0){10}}
\put(14,1.5){\vector(1,0){10}}
\put(1,2){\makebox(0,0){$A$}}
\put(11,2){\makebox(0,0){$A\times H$}}
\put(28,2){\makebox(0,0){$A\times H\ot H$}}
\put(5.5,3){\makebox(0,0){${i_{A\times_{\sigma} H}}$}}
\put(18,4){\makebox(0,0){$\rho_{A\times_{\sigma} H}$}}
\put(18,0.0){\makebox(0,0){$(A\times H\ot \Pi_{H}^{L})\co
\rho_{A\times_{\sigma} H}$}}
\end{picture}
$$

is an equalizer diagram, where
$$i_{A\times_{\sigma} H}=p_{A\otimes H}\co (A\ot \eta_H).$$

Indeed, using that $\nabla_{A\otimes H}$ is a morphism of
right $H$-comodules (see (\ref{nabla-u1})) and (\ref{deltaPIL}) we have:

\begin{itemize}

\item[ ]$\hspace{0.38cm} \rho_{A\times_{\sigma} H}\co i_{A\times_{\sigma} H}$

\item[ ]

\item[ ]$=(p_{A\otimes H}\ot H)\co (A\ot \delta_H)\co \nabla_{A\otimes H}\co (A\ot \eta_H)$

\item[ ]

\item[ ]$=(p_{A\otimes H}\ot \Pi^{L}_{H})\co (A\ot (\delta_H\co \eta_H))$

\item[ ]

\item[ ]$=(p_{A\otimes H}\ot (\Pi^{L}_{H}\co \Pi^{L}_{H}))\co (A\ot (\delta_H\co \eta_H))$

\item[ ]

\item[ ]$=(A\times H\ot \Pi^{L}_{H})\co \rho_{A\times_{\sigma} H}
\co i_{A\times_{\sigma} H}.$

\end{itemize}

Moreover, if $g:Q\rightarrow A\times H$ is a morphism such
that $\rho_{A\times_{\sigma} H}\co g=(A\times H\ot \Pi^{L}_{H})\co
\rho_{A\times_{\sigma} H}\co g$ we obtain

\begin{itemize}

\item[ ]$\hspace{0.38cm} (A\ot \delta_H)\co i_{A\otimes H}\co g$

\item[ ]

\item[ ]$=(\nabla_{A\otimes H}\ot H)\co (A\ot \delta_H)\co i_{A\otimes H}\co g$

\item[ ]

\item[ ]$=(i_{A\otimes H}\ot H)\co \rho_{A\times_{\sigma} H}\co g$

\item[ ]

\item[ ]$=(i_{A\otimes H}\ot \Pi^{L}_{H})\co \rho_{A\times_{\sigma} H}\co g$

\item[ ]

\item[ ]$=(A\ot H\ot \Pi^{L}_{H})\co (A\ot \delta_H)\co i_{A\otimes H}\co g,$

\end{itemize}

and then
\begin{equation}\label{propiedadfactorizacion}
i_{A\otimes H}\co g=(A\ot  \Pi^{L}_{H})\co i_{A\otimes H}\co g.
\end{equation}

Now we will show that $h=(A\ot \varepsilon_H)\co i_{A\otimes H}\co g$
is the unique morphism such that $i_{A\otimes H}\co h=g$.

First note that, by (\ref{nabla-u1}), we have
$$h=(A\ot \varepsilon_{H})\co i_{A\otimes H}\co g=(A\ot \varepsilon_{H})\co
\nabla_{A\otimes H}\co (A\ot  \Pi^{L}_{H})\co i_{A\otimes H}\co g$$
$$=\mu_{A}\co (A\ot u_{1})\co i_{A\otimes H}\co g, $$
and the equality
\begin{equation}\label{newnabla-2}
\nabla_{A\otimes H}\co ((\mu_{A}\co (A\ot u_{1}))\ot \eta_{H})=
(A\ot \Pi_{H}^{L})\co \nabla_{A\otimes H}
\end{equation}
holds because
\begin{itemize}

\item[ ]$\hspace{0.38cm} \nabla_{A\otimes H}\co ((\mu_{A}\co (A\ot u_{1}))\ot \eta_{H})$

\item[ ]

\item[ ]$=(\mu_{A}\ot H)\co (A\ot (\nabla_{A\otimes H}\co (u_{1}\ot \eta_{H})))$

\item[ ]

\item[ ]$=(\mu_{A}\ot H)\co (A\ot (\psi_{H}^{A}\co (\eta_{H}\ot u_{1})))$

\item[ ]

\item[ ]$=(\mu_{A}\ot H)\co (A\ot (((u_{1}\co \mu_{H})\ot H)\co
(H\ot c_{H,H})\co ((\delta_{H}\co \eta_{H})\ot H)))$

\item[ ]

\item[ ]$=(\mu_{A}\ot H)\co (A\ot ((u_{1}\ot \Pi_{H}^{L})\co \delta_{H}) $

\item[ ]

\item[ ]$=(A\ot \Pi_{H}^{L})\co \nabla_{A\otimes H}$

\end{itemize}

In the foregoing calculations, the first equality follows using that $\nabla_{A\otimes H}$ is a
morphism of left $A$-modules;
the second one is a consequence of the definition of $\nabla_{A\otimes H}$; the third
follows by (d3) of the definition of weak left $H$-module algebra; the fourth one
by (\ref{deltaPIL}), and finally, the last one
by (\ref{nabla-u1}).

Therefore,
$$i_{A\times_{\sigma} H}\co h=p_{A\otimes H}\co ((\mu_{A}\co (A\ot u_{1})
\co i_{A\otimes H}\co g)\ot \eta_{H})=p_{A\otimes H}\co (A\ot  \Pi^{L}_{H})\co i_{A\otimes H}\co g=g $$

Finally, if $r:Q\rightarrow A$ is a morphism such that $i_{A\times H}\co r=g$, we have
$$r=(A\ot \varepsilon_H)\co \nabla_{A\otimes H}\co (r\ot \eta_H)=(A\ot \varepsilon_H)
\co i_{A\otimes H}\co g=h$$
and we conclude the proof.

\end{dem}

In the following proposition we establish the relation between crossed
systems and $H$-cleft extensions.

\begin{prop}
\label{sistemacruzadocleft}Let $H$ be a cocommutative weak Hopf algebra and
 $(\varphi_A, \sigma)$  a crossed system for $H$ over $A$.
 Then $A\hookrightarrow A\times_{\sigma}H$ is an $H$-cleft extension.

\end{prop}

\begin{dem}

The morphism $f=p_{A\otimes H}\co (\eta_A\ot H):H\rightarrow A\times_{\sigma} H$ is a total integral.
Obviously, $f\co\eta_H=\eta_{A\times_{\sigma} H}$. Moreover, using that
$\nabla_{A\otimes H}$ is a morphism of right $H$-comodules
we get that $f$ is an integral because:

\begin{itemize}

\item[ ]$\hspace{0.38cm} \rho_{A\times_{\sigma} H}\co f$

\item[ ]

\item[ ]$=(p_{A\otimes H}\ot H)\co (A\ot \delta_H)\co\nabla_{A\otimes H}\co (\eta_A\ot H)$

\item[ ]

\item[ ]$=(f\ot H)\co \delta_H.$

\end{itemize}

We define $f^{-1}=p_{A\otimes H}\co(\sigma^{-1}\ot H)\co (H\ot c_{H,H})\co
((\delta_H\co \lambda_H)\ot H)\co \delta_H$.
We will show that $f^{-1}$ is the convolution inverse of $f$. First note that (c1) holds:

\begin{itemize}

\item[ ]$\hspace{0.38cm} f^{-1}\wedge f$

\item[ ]

\item[ ]$=p_{A\otimes H}\co (\mu_A\ot \sigma_{H}^{A})\co (A\ot ((u_{1}\ot H)\co \delta_{H})
\ot H)\co (((\sigma^{-1}\ot H)\co (H\ot c_{H,H})\co ((\delta_{H}\co
\lambda_{H})\ot H)\co \delta_{H})\ot H)\co \delta_{H}$

\item[ ]

\item[ ]$=p_{A\otimes H}\co ((\mu_{A}\co (A\ot (u_{2}\wedge \sigma)))\ot \mu_{H})\co
(\sigma^{-1}\ot \delta_{H\ot H})\co ((( (H\ot c_{H,H})\co ((\delta_{H}\co \lambda_{H})\ot H)\co
\delta_{H})\ot H)\co \delta_{H}$

\item[ ]

\item[ ]$=p_{A\otimes H}\co (\mu_{A}\ot H)\co (\sigma^{-1}\ot \sigma_{H}^{A})\co
((( (H\ot c_{H,H})\co ((\delta_{H}\co \lambda_{H})\ot H)\co
\delta_{H})\ot H)\co \delta_{H}$

\item[ ]

\item[ ]$=p_{A\otimes H}\co ((\sigma^{-1}\wedge \sigma)\ot \mu_H)\co \delta_{H\ot H}\co
(\lambda_H\ot H)\co \delta_H$

\item[ ]

\item[ ]$=p_{A\otimes H}\co
 ((u_{1}\co \mu_H)\ot \mu_H)\co \delta_{H\ot H}\co (\lambda_H\ot H)\co \delta_H$

\item[ ]

\item[ ]$=p_{A\otimes H}\co (u_{1}\ot H)\co\delta_H\co \Pi^{R}_{H}$

\item[ ]

\item[ ]$=p_{A\otimes H}\co \nabla_{A\otimes H}\co (\eta_A\ot \Pi^{R}_{H})$

\item[ ]

\item[ ]$=f\co \Pi^{R}_{H}$

\item[ ]

\item[ ]$=(A\times H\ot(\varepsilon_H\co\mu_H))\co
(c_{H,A\times H}\ot H)\co (H\ot (\rho_{A\times_{\sigma} H}\co \eta_{A\times_{\sigma} H})).$

\end{itemize}

In the previous calculations, the first equality follows by the
normalized condition for the product $\mu_{A\ot_{\sigma}H}$; the second one
uses that $(A, \varphi_A)$ is a weak left $H$-module, the coassociativity of $\delta_{H}$
and the naturality of $c$; the third one relies on  the coassociativity of $\delta_{H}$, the naturality of $c$ and
because $\sigma$ is in $Reg_{\varphi_{A}}(H\ot H,A)$; the fourth one is a consequence of
the coassociativity of $\delta_{H}$ and the naturality of $c$. The fifth uses that
$\sigma\in Reg_{\varphi_{A}}(H\ot H,A)$. The sixth follows by the definition
of $\Pi^{R}_{H}$ and that $H$ is a weak Hopf algebra; the seventh one  by (\ref{nabla-u1}),
the eigth is a consequence
of the equality $p_{A\otimes H}\co \nabla_{A\otimes H}=p_{A\otimes H}$; finally,
in the last one we use that $f$ is a total integral.

In order to give condition (c2), we will show the equality
\begin{equation}\label{igualdadparacociclo}
((\sigma\co (H\ot\mu_H))\wedge (\sigma^{-1}\co (\mu_H\ot H)))\co
(H\ot \lambda_H\ot H)\co (H\ot \delta_H)\co \delta_H=u_{1}.
\end{equation}

Indeed: Using the anticomultiplicativity of the antipode, the coassocitivity
of $\delta_{H}$, the naturality of $c$,
 the equalities (\ref{deltaPIL}) and (\ref{deltaPIR}) as well as
$\sigma\in Reg_{\varphi_{A}}(H\ot H,A)$ and (\ref{new-two-cocycle-3}), we get
(\ref{igualdadparacociclo}), because

\begin{itemize}

\item[ ]$\hspace{0.38cm} ((\sigma\co (H\ot\mu_H))\wedge (\sigma^{-1}\co (\mu_H\ot H)))\co (H\ot \lambda_H\ot H)\co (H\ot \delta_H)\co \delta_H$

\item[ ]

\item[ ]$=\mu_A\co (\sigma\ot \sigma^{-1})\co (H\ot c_{H,H}\ot H)\co (((H\ot \Pi^{L}_{H})\co \delta_H)\ot ((\Pi^{R}_{H}\ot H)\co \delta_H))\co \delta_H$

\item[ ]

\item[ ]$=((\sigma\co (\mu_H\ot H))\wedge (\sigma^{-1}\co (H\ot \mu_H)))\co (\eta_H\ot H\ot \eta_H)$

\item[ ]

\item[ ]$=\mu_A\co (\sigma^{-1}\ot \varphi_A)\co (H\ot c_{H,H}\ot (\sigma\co (H\ot \eta_H)))\co ((\delta_H\co\eta_H)\ot \delta_H)$

\item[ ]

\item[ ]$=(\sigma^{-1}\wedge u_{2}))\co (\eta_H\ot H)$

\item[ ]

\item[ ]$=\sigma^{-1}\co (\eta_H\ot H)$

\item[ ]

\item[ ]$=u_{1}.$

\end{itemize}

Now we prove (c2):

\begin{itemize}

\item[ ]$\hspace{0.38cm} f\wedge f^{-1}$

\item[ ]

\item[ ]$=p_{A\otimes H}\co (\mu_{A}\ot H)\co (A\ot \sigma_{H}^{A})\co (\psi_{H}^{A}\ot H)\co (H
\ot ((\sigma^{-1}\ot H)\co (H\ot c_{H,H})\co ((\delta_{H}\co \lambda_{H})\ot H)\co \delta_{H}))
\co \delta_{H}$

\item[ ]

\item[ ]$=p_{A\otimes H}\co ((\mu_{A}\co ((\varphi_{A}\co (H\ot \sigma^{-1}))\ot \sigma)\co (H\ot H\ot c_{H,H}\ot H)\co(H\ot c_{H,H}\ot c_{H,H}\ot H)$

\item[ ]

\item[ ]$\hspace{0.38cm}\co (\delta_{H}\ot \delta_{H}\ot H))\ot H)\co (H\ot H\ot c_{H,H})\co (H\ot ((H\ot \mu_{H})\co (c_{H,H}\ot H)\co (H\ot
(\delta_{H}\co\lambda_{H})))\ot H)$

\item[ ]

\item[ ]$\hspace{0.38cm} \co (\delta_{H}\ot \delta_{H})\co \delta_{H}$

\item[ ]

\item[ ]$=p_{A\otimes H}\co ((\mu_{A}\co (\sigma\ot \sigma^{-1})\co (H\ot \mu_{H}\ot \mu_{H}\ot H)\co
\delta_{H\ot H\ot H})\ot H)\co  (H\ot H\ot c_{H,H})$

\item[ ]

\item[ ]$\hspace{0.38cm}  \co(H\ot ((H\ot \mu_{H})\co (c_{H,H}\ot H)\co (H\ot
(\delta_{H}\co\lambda_{H})))\ot H)\co (\delta_{H}\ot \delta_{H})\co \delta_{H}$

\item[ ]

\item[ ]$=p_{A\otimes H}\co ((\mu_{A}\co (\sigma\ot \sigma^{-1})\co (H\ot \mu_{H}\ot \mu_{H}\ot H)\co
\delta_{H\ot H\ot H})\ot H)\co (H\ot H\ot c_{H,H}) $

\item[ ]

\item[ ]$\hspace{0.38cm}  \co (H\ot ((H\ot \mu_{H})\co (c_{H,H}\ot H)\co (H\ot
((\lambda_{H}\ot \lambda_{H})\co c_{H,H}\co \delta_{H}))\ot H)\co
(\delta_{H}\ot \delta_{H})\co \delta_{H}$

\item[ ]

\item[ ]$=p_{A\otimes H}\co ((\mu_{A}\co (\sigma\ot \sigma^{-1})\co (H\ot \mu_{H}\ot \mu_{H}\ot H)\co
\delta_{H\ot H\ot H})\ot H)\co (H\ot H\ot c_{H,H})$

\item[ ]

\item[ ]$\hspace{0.38cm} \co (H\ot (c_{H,H}\co (\Pi_{H}^{L}\ot \lambda_{H}))
\ot H)\co (\delta_{H}\ot \delta_{H})\co \delta_{H}$

\item[ ]

\item[ ]$=p_{A\otimes H}\co ((\mu_{A}\co (\sigma\ot \sigma^{-1})\co (H\ot \mu_{H}\ot \mu_{H}\ot H)\co
\delta_{H\ot H\ot H})\ot H)$

\item[ ]

\item[ ]$\hspace{0.38cm} \co (H\ot ((((\lambda_{H}\ot H)\co \delta_{H})\ot \Pi_{H}^{L})\co
c_{H,H}\co \delta_{H}))\co \delta_{H}$

\item[ ]

\item[ ]$=p_{A\otimes H}\co ((((\sigma\co (H\ot\mu_H))\wedge (\sigma^{-1}\co (\mu_H\ot H)))\co
(H\ot \lambda_H\ot H)\co (H\ot \delta_H)\co \delta_H)\ot \Pi_{H}^{L})\co \delta_{H}$

\item[ ]

\item[ ]$=p_{A\otimes H}\co (u_{1}\ot\Pi^{L}_{H})\co \delta_H$

\item[ ]

\item[ ]$=p_{A\otimes H}\co (A\ot \Pi^{L}_{H})\co \nabla_{A\otimes H}\co (\eta_A\ot H)$

\item[ ]

\item[ ]$=p_{A\otimes H}\co  \nabla_{A\otimes H}\co(\eta_A\ot \overline{\Pi}^{L}_{H})$

\item[ ]

\item[ ]$=f\co \overline{\Pi}^{L}_{H}$

\item[ ]

\item[ ]$=(A\times H\ot (\varepsilon_H\co\mu_H))\co
((\rho_{A\times_{\sigma} H}\co \eta_{A\times_{\sigma} H})\ot H).$

\end{itemize}

The first equality follows by (\ref{fi-nab}), the second one by
the naturality of $c$ and the coassociativity of $\delta_{H}$, the third one by
(\ref{new-two-cocycle-2}), the fourth one by the antimultiplicative property of $\lambda_{H}$, the fifth one by
the naturality of $c$ and the definition of $\Pi_{H}^{L}$. In the sixth one we use the naturality of $c$ and the seventh
follows by the cocommutativity of $H$ and the coassociativity of $\delta_{H}$. The eight one is a
consequence of (\ref{igualdadparacociclo}) and the ninth of (\ref{nabla-u1}). In the tenth
equality we use that if $H$
is cocommutative $\Pi^{L}_{H}=\overline{\Pi}^{L}_{H}$ and
$(A\ot \overline{\Pi}^{L}_{H})\co \nabla_{A\otimes H}=
\nabla_{A\otimes H}\co (A\ot \overline{\Pi}^{L}_{H})$. The eleventh one relies on
the properties of $\nabla_{A\ot H}$ and
finally the last one follows because $f$ is a total integral.

To finish the proof we only need to show that $f^{-1}\wedge f\wedge f^{-1}=f^{-1}$.
First of all, using that $H$ is cocommutative it is not difficult to see
that $(f^{-1}\ot \lambda_H)\co c_{H,H}\co \delta_H=\rho_{A\times_{\sigma} H}\co f^{-1}$.
Then by this equality, the fact that $\lambda_H\co \lambda_H=id_H$ (which follows
because $H$ is cocommutative) and (\ref{composiciones}),

\begin{itemize}

\item[ ]$\hspace{0.38cm} f^{-1}\wedge f\wedge f^{-1}$

\item[ ]

\item[ ]$=\mu_{A\times_{\sigma} H}\co (f^{-1}\ot ((f\wedge f^{-1})\co\lambda_H\co \lambda_H))
\co \delta_H$

\item[ ]

\item[ ]$=\mu_{A\times_{\sigma} H}\co (f^{-1}\ot
(f\co\overline{\Pi}^{L}_{H}\co\lambda_H\co \lambda_H))\co c_{H,H}\co \delta_H$

\item[ ]

\item[ ]$=\mu_{A\times_{\sigma} H}\co(A\times H\ot (f\co \Pi^{R}_{H}))\co
\rho_{A\times H}\co f^{-1}$

\item[ ]

\item[ ]$=\mu_{A\times_{\sigma} H}\co(A\times H\ot (f^{-1}\wedge f))\co
\rho_{A\times H}\co f^{-1}$

\item[ ]

\item[ ]$=f^{-1},$

\end{itemize}

and we conclude the proof.

\end{dem}

\begin{prop}
\label{first-step}
Let $H$ be a cocommutative weak Hopf algebra and let $A$ be an algebra.
If $(\varphi_A, \alpha)$ and $(\phi_A, \beta)$ are two equivalent crossed systems, so are the associated $H$-cleft extensions $A\hookrightarrow A\times_{\alpha}H$ and $A\hookrightarrow A\times_{\beta}H$.
\end{prop}

\begin{dem}We will begin showing that this correspondence is well defined: Let $h$ be the morphism in
$Reg_{\varphi_{A}}(H,A)\cap Reg_{\phi_{A}}(H,A)$ satisfying
conditions (\ref{relacionfis}) and (\ref{relacionsigmas}). We denote by
$A\hookrightarrow A\times_{\alpha}H$ and $A\hookrightarrow A\times_{\beta}H$
the $H$-cleft extensions defined by $(\varphi_A, \alpha)$ and $(\phi_A, \beta)$,
respectively. We will show that
\begin{equation}
\label{isomorfismocleft}
T=p_{A\otimes H}\co (\mu_A\ot H)\co (A\ot h\ot H)\co (A\ot \delta_H)\co i_{A\ot H}
\end{equation}
is a morphism of $H$-comodule algebras such that
$ T\co i_{A\times_{\alpha} H}=i_{A\times_{\beta} H}$.
Firstly of all, note that, by (\ref{nablaprima-nabla})
the idempotent morphisms  defined by the two crossed
systems coincide. We denote it by $\nabla_{A\ot H}$. Moreover,

\begin{itemize}

\item[ ]$\hspace{0.38cm} T\co \eta_{A\times_{\alpha} H}$

\item[ ]

\item[ ]$=p_{A\otimes H}\co (\mu_A\ot H)\co (A\ot h\ot H)\co (A\ot \delta_H)\co \nabla_{A\ot H}\co (\eta_A\ot \eta_H)$

\item[ ]

\item[ ]$=p_{A\otimes H}\co ((u_{1}\wedge h)\ot H)\co \delta_H\co \eta_H$

\item[ ]

\item[ ]$=p_{A\otimes H}\co (h\ot H)\co \delta_H\co \eta_H$

\item[ ]

\item[ ]$=p_{A\otimes H}\co ((h\co \overline{\Pi}^{L}_{H})\ot H)\co \delta_H\co \eta_H$

\item[ ]

\item[ ]$=p_{A\otimes H}\co ((\varphi_A\co (H\ot \eta_A))\ot H)\co \delta_H\co \eta_H$

\item[ ]

\item[ ]$=\eta_{A\times_{\beta} H}.$

\end{itemize}

In the foregoing calculations, the first and the last equalities follow by the definition
of $\nabla_{A\ot H}$; in the second one we use that $\nabla_{A\ot H}$ is a morphism of
left $A$-modules and right $H$-comodules; the third that $h$ is in $Reg_{\varphi_{A}}(H,A)$.
Finally, the fourth equality follows by (\ref{deltaPILbarra}) and the fifth one by
(\ref{hbarPiL}).

To prove the multiplicative condition for $T$ we need to fix a new notation and get two
auxiliary identities. First note that, by \ref{weak-crossed-products-exposition}, the crossed systems $(\varphi_A, \alpha)$ and
$(\phi_A, \beta)$ define two quadruples
$$(A,H,\psi_{H,\alpha}^{A}=(\varphi_{A}\ot H)\co (H\ot c_{H,A})\co
(\delta_{H}\ot A),\sigma_{H,\alpha}^{A}=(\alpha\ot \mu_{H})\co \delta_{H\ot H})$$
$$(A,H,\psi_{H,\beta}^{A}=(\phi_{A}\ot H)\co (H\ot c_{H,A})\co
(\delta_{H}\ot A),\sigma_{H,\beta}^{A}=(\beta\ot \mu_{H})\co \delta_{H\ot H})$$
that induce the corresponding weak crossed products. On the other hand,
the following equalities hold:
\begin{equation}
\label{prin-1-aux}
\mu_{A}\co (A\ot h)\co \sigma_{H,\alpha}^{A}=\mu_{A}\co ((\mu_{A}\co (h\ot \phi_{A})
\co (\delta_{H}\ot h))\ot \beta)\co \delta_{H\ot H},
\end{equation}
\begin{equation}
\label{prin-2-aux}
(\nabla_{A\ot H}\ot H)\co ((\mu_{A}\co (h\ot \phi_{A})
\co (\delta_{H}\ot A))\ot \delta_{H})\co (H\ot c_{H,A})\co (\delta_{H}\ot A)=
(\mu_{A}\ot \delta_{H})\co
(h\ot\psi_{H,\beta}^{A})\co (\delta_{H}\ot A).
\end{equation}

The proof for (\ref{prin-1-aux}) is as follows:
\begin{itemize}

\item[ ]$\hspace{0.38cm}\mu_{A}\co (A\ot h)\co \sigma_{H,\alpha}^{A} $

\item[ ]

\item[ ]$=\mu_{A}\co ((\mu_{A}\co (\mu_{A}\ot h^{-1})\co
((\mu_{A}\co (h\ot \phi_{A})\co (\delta_{H}\ot h))\ot \sigma_{H,\beta}^{A}))\ot
(h\co \mu_{H}))\co
(\delta_{H\ot H}\ot H\ot H)\co \delta_{H\ot H}  $

\item[ ]

\item[ ]$=\mu_{A}\co ((\mu_{A}\co (h\ot \phi_{A})\co (\delta_{H}\ot h))\ot (\mu_{A}\co (\beta\ot
(h^{-1}\wedge h))))\co
(\delta_{H\ot H}\ot \mu_H)\co \delta_{H\ot H} $

\item[ ]

\item[ ]$=\mu_{A}\co ((\mu_{A}\co (h\ot \phi_{A})\co (\delta_{H}\ot h))\ot
(\beta\wedge (\phi_{A}\co (\mu_{H}\ot \eta_{A}))))\co
\delta_{H\ot H} $

\item[ ]

\item[ ]$=\mu_{A}\co ((\mu_{A}\co (h\ot \phi_{A})
\co (\delta_{H}\ot h))\ot \beta)\co \delta_{H\ot H}$

\end{itemize}

where the first equality follows by the equivalence between $(\varphi_{A}, \alpha)$ and
$(\phi_{A}, \beta)$, the second one by the coassociativity of $\delta_{H\ot H}$ and (a1)
 of the definition of weak Hopf algebra, the third one by the coassociativity of $\delta_{H\ot H}$ and
 because $h\in Reg_{\varphi_{A}}(H,A)$ and finally, in the last one we use that
$\beta\in Reg_{\varphi_{A}}(H\ot H,A)$.

Taking into account that $\nabla_{A\ot H}$ is a morphism of left $A$-modules and right $H$comodules,
the coassociativity of $\delta_{H}$ and (\ref{fi-nab-2}), we obtain (\ref{prin-2-aux}) because:

\begin{itemize}

\item[ ]$\hspace{0.38cm}(\nabla_{A\ot H}\ot H)\co ((\mu_{A}\co (h\ot \phi_{A})
\co (\delta_{H}\ot A))\ot \delta_{H})\co (H\ot c_{H,A})\co (\delta_{H}\ot A) $

\item[ ]

\item[ ]$=(A\ot \delta_{H})\co \nabla_{A\ot H}\co (\mu_{A}\ot H)\co
(h\ot \psi_{H,\beta}^{A})\co
(\delta_{H}\ot A)$

\item[ ]

\item[ ]$=(\mu_{A}\ot \delta_{H})\co
(h\ot\psi_{H,\beta}^{A})\co (\delta_{H}\ot A).$

\end{itemize}

We are now in position to show that $T$ is a multiplicative morphism.

\begin{itemize}

\item[ ]

\item[ ]$\hspace{0.38cm} T\co \mu_{A\times_{\alpha} H}$

\item[ ]

\item[ ]$=p_{A\otimes H}\co ((\mu_{A}\co (A\ot h))\ot H)\co (A\ot \delta_{H})\co
\nabla_{A\ot H}\co (\mu_{A}\ot H)\co (\mu_{A}\ot \sigma_{H,\alpha}^{A})\co (A\ot
\psi_{H,\alpha}^{A}\ot H)$

\item[ ]

\item[ ]$\hspace{0.38cm}\co (i_{A\ot H}\ot i_{A\ot H})$

\item[ ]

\item[ ]$=p_{A\otimes H}\co (\mu_{A}\ot H)\co
(A\ot (\mu_{A}\co (A\ot h)\co \sigma_{H,\alpha}^{A})\ot \mu_{H})\co (\mu_{A}\ot \delta_{H\ot H})\co
(A\ot \psi_{H,\alpha}^{A}\ot H)\co  (i_{A\ot H}\ot i_{A\ot H}) $

\item[ ]

\item[ ]$=p_{A\otimes H}\co (\mu_{A}\ot H)\co (A\ot \mu_{A}\ot H)$

\item[ ]

\item[ ]$\hspace{0.38cm} \co(A\ot ((\mu_{A}\co ((\mu_{A}\co (h\ot A))\ot A)\co (H\ot
\phi_{A}\ot h^{-1})\co (\delta_{H}\ot c_{H,A})\co (\delta_{H}\ot
A))$

\item[ ]

\item[ ]$\hspace{0.38cm} \ot (\mu_{A}\co ((\mu_{A}\co (h\ot \phi_{A})
\co (\delta_{H}\ot h))\ot \beta))\co \delta_{H\ot H})\ot \mu_{H})\co
(A\ot H\ot A\ot \delta_{H\ot H})$

\item[ ]

\item[ ]$\hspace{0.38cm} \co (A\ot H\ot c_{H,A}\ot H)\co
(A\ot \delta_{H}\ot A\ot H)\co (i_{A\ot H}\ot i_{A\ot H})$

\item[ ]

\item[ ]$=p_{A\otimes H}\co (\mu_{A}\ot H)\co (\mu_{A}\ot A\ot H)\co ((\mu_{A}\co
(A\ot (\mu_{A}\co (h\ot \phi_{A})\co (\delta_{H}\ot A))))$

\item[ ]

\item[ ]$\hspace{0.38cm}\ot(\mu_{A}\co ((h^{-1}\wedge h)\ot \phi_{A})\co
(\delta_{H}\ot h))\ot \sigma_{H,\beta}^{A})\co(A\ot H\ot A\ot \delta_{H\ot H})\co (A\ot H\ot c_{H,A}\ot H)$

\item[ ]

\item[ ]$\hspace{0.38cm} \co
(A\ot \delta_{H}\ot A\ot H)\co (i_{A\ot H}\ot i_{A\ot H})$

\item[ ]

\item[ ]$=p_{A\otimes H}\co (\mu_{A}\ot H)\co (\mu_{A}\ot A\ot H)\co (A\ot (\phi_{A}\co (H\ot h)\ot \sigma_{H,\beta}^{A}))\co (\nabla_{A\ot H}\ot H\ot H\ot H)$

\item[ ]

\item[ ]$\hspace{0.38cm}\co (\mu_A\co (A\ot (\mu_A\co (h\ot \phi_{A})))\ot \delta_{H\ot H})\co (A\ot \delta_H\ot c_{H,A}\ot H)\co (A\ot \delta_H\ot A\ot H)\co (i_{A\ot H}\ot i_{A\ot H})$

\item[ ]

\item[ ]$=p_{A\otimes H}\co (\mu_{A}\ot H)\co (\mu_{A}\ot \sigma_{H,\beta}^{A})\co (\mu_{A}\ot \phi_{A}\ot H\ot H)\co
(A\ot A\ot H\ot c_{H,A}\ot H)$

\item[ ]

\item[ ]$\hspace{0.38cm} \co (A\ot ((\mu_{A}\ot \delta_{H})\co
(h\ot\psi_{H,\beta}^{A})\co (\delta_{H}\ot A))\ot ((h\ot H)\co \delta_{H}))\co
\co (i_{A\ot H}\ot i_{A\ot H})$

\item[ ]

\item[ ]$=p_{A\otimes H}\co (\mu_{A}\ot H)\co (\mu_{A}\ot \sigma_{H,\beta}^{A})\co (\mu_{A}\ot
((\mu_{A}\ot H)\co (A\ot \psi_{H,\beta}^{A})\co (\psi_{H,\beta}^{A}\ot A))\ot H)$

\item[ ]

\item[ ]$\hspace{0.38cm} \co (A\ot ((h\ot H)\co \delta_{H})\ot A\ot ((h\ot H)\co \delta_{H}))
\co (i_{A\ot H}\ot i_{A\ot H})$

\item[ ]

\item[ ]$=p_{A\otimes H}\co (\mu_{A}\ot H)\co (\mu_{A}\ot \sigma_{H,\beta}^{A})\co
(A\ot \psi_{H,\beta}^{A})\ot H)$

\item[ ]

\item[ ]$\hspace{0.38cm} \co ((((\mu_{A}\co (A\ot h))\ot H)\co (A\ot \delta_{H}))\ot
(((\mu_{A}\co (A\ot h))\ot H)\co (A\ot \delta_{H})) ) \co (i_{A\ot
H}\ot i_{A\ot H})$

\item[ ]

\item[ ]$=\mu_{A\times_{\beta} H}\co (T\ot T).$

\end{itemize}

In the previous calculations, the first equality follows by the definition, the second one because
 $\nabla_{A\ot H}$  is a morphism of left $A$-modules and  (\ref{Nablasigma}),
 (\ref{delta-sigmaHA}) hold, the third one
 is a consequence of the equivalence between $(\varphi_{A}, \alpha)$ and
$(\phi_{A}, \beta)$ and (\ref{prin-1-aux}),
the fourth one follows by the associativity of $\mu_{A}$, the coassociativity of
$\delta_{H}$ and the naturality of $c$, the fifth one by the associativity of $\mu_{A}$,
the sixth one by (\ref{prin-2-aux}), the seventh one by the definition of $\psi_{H,\beta}^{A}$,
the eight one by (\ref{wmeas-wcp}) and, finally, the ninth one by the normalized condition for
$\mu_{A\ot_{\beta}H}$.

Using that
\begin{equation}\label{otraigualdadparaT}
i_{A\ot H}\co T=(\mu_A\ot H)\co (A\ot h\ot H)\co (A\ot \delta_H)\co i_{A\ot H},
\end{equation}

it is not difficult to see that $T$ is a morphism of right $H$-comodules.
Moreover, by associativity of $\mu_{A}$, the coassociativity of $\delta_{H}$,
(\ref{deltaPILbarra}) and (\ref{hbarPiL}), we have

\begin{itemize}

\item[ ]$\hspace{0.38cm} T\co i_{A\times_{\alpha} H}$

\item[ ]

\item[ ]$=p_{A\otimes H}\co (\mu_A\ot H)\co (A\ot h\ot H)\co (A\ot \delta_H)\co
\nabla_{A\ot H}\co (A\ot \eta_H)$

\item[ ]

\item[ ]$=p_{A\otimes H}\co (\mu_A\ot H)\co (A\ot ((u_{1}\wedge h)
\ot H)\co (A\ot \delta_H)\co (A\ot \eta_H)$

\item[ ]

\item[ ]$=p_{A\otimes H}\co (\mu_A\ot H)\co (A\ot ( h\ot H))\co
(A\ot \delta_H)\co (A\ot \eta_H)$

\item[ ]

\item[ ]$=p_{A\otimes H}\co (\mu_A\ot H)\co (A\ot (h\co
\overline{\Pi}^{L}_{H})\ot H)\co (A\ot \delta_H)\co (A\ot \eta_H)$

\item[ ]

\item[ ]$=p_{A\otimes H}\co (\mu_A\ot H)\co (A\ot u_{1}\ot H)
\co (A\ot \delta_H)\co (A\ot \eta_H)$

\item[ ]

\item[ ]$=p_{A\otimes H}\co \nabla_{A\ot H}\co (A\ot \eta_H)$

\item[ ]

\item[ ]$=i_{A\times_{\beta} H}.$

\end{itemize}

Therefore, the associated $H$-cleft extensions $A\hookrightarrow A\times_{\alpha}H$ and $A\hookrightarrow A\times_{\beta}H$ are equivalent.

\end{dem}

\begin{rem}
\label{commentstofirststep}
{\rm In the conditions of the previous result, $T$ is an isomorphism  with inverse
$$T^{-1}=p_{A\otimes H}\co (\mu_A\ot H)\co (A\ot h^{-1}\ot H)\co
(A\ot \delta_H)\co i_{A\ot H}$$
because
\begin{itemize}

\item[ ]$\hspace{0.38cm} T^{-1}\co T$

\item[ ]

\item[ ]$=p_{A\otimes H}\co (\mu_A\ot H)\co (A\ot (h\wedge h^{-1})\ot H)\co
(A\ot \delta_H)\co i_{A\ot H}$

\item[ ]

\item[ ]$=p_{A\otimes H}\co (\mu_A\ot H)\co (A\ot u_{1}
\ot H)\co (A\ot \delta_H)\co i_{A\ot H}$

\item[ ]

\item[ ]$=p_{A\otimes H}\co \nabla_{A\ot H}\co i_{A\ot H}$

\item[ ]

\item[ ]$=id_{A\times_{\alpha} H}.$

\end{itemize}
}
\end{rem}

\begin{prop}
\label{construccioncociclo}
Let $H$ be a cocommutative weak Hopf algebra. If $A_{H}\hookrightarrow A$ is an
$H$-cleft extension, the morphism
$$\sigma_{A}:=(\mu_A\co (f\ot f))\wedge (f^{-1}\co \mu_H):H\ot H\rightarrow A, $$
where $f:H\rightarrow A$ is a convolution invertible total integral,
factors through the equalizer $i_A$. Moreover, if
$\varphi_{A_H}:H\ot A_H\rightarrow A_H$ is the weak left $H$-module
structure defined in Proposition \ref{estructuramoduloinducida}, the
factorization of $\sigma_{A}$ is a morphism in
$Reg_{\varphi_{A_{H}}}(H\ot H,A_{H})$ satisfying the normal
condition (g3) and  with convolution inverse the factorization
through the equalizer $i_A$ of the morphism
$$\sigma_{A}^{-1}:=(f\co \mu_H)\wedge (\mu_A\co c_{A,A}\co  (f^{-1}\ot f^{-1})).$$

\end{prop}

\begin{dem}
If $A_{H}\hookrightarrow A$ is
a $H$-cleft extension, by Corollary \ref{prin-cor}, we have that $A_{H}\hookrightarrow A$ is a weak $H$-cleft extension. Then, by Proposition 1.17 of \cite{nmra1}, we obtain that $\sigma_{A}$ factors through the equalizer $i_A$ and, if $\sigma_{A_{H}}$ is the factorization of $\sigma_{A}$, the equality
\begin{equation}
\label{sigma-f}
\sigma_{A_{H}}=p_{A}\co \mu_{A}\co (f\ot f)
\end{equation}
holds.

On the other hand, the morphism $\sigma_{A}^{-1}$ factors through the equalizer $i_A$ because
$$\rho_{A}\circ
\sigma^{-1}_{A}=(A\otimes \overline{\Pi}_{H}^{R})\circ
\rho_{A}\circ \sigma^{-1}_{A}.$$
Indeed:

\begin{itemize}

\item[ ]$\hspace{0.38cm}  \rho_{A}\circ
\sigma^{-1}_{A}$

\item[ ]

\item[ ]$=\mu_{A\otimes H}\circ (\rho_A\otimes (\rho_A\co \mu_A))\circ (A\otimes c_{A,A})\circ ((f\circ
\mu_{H})\otimes f^{-1}\otimes f^{-1})\circ \delta_{H\otimes H}$

\item[ ]

\item[ ]$=\mu_{A\otimes H}\circ (f\otimes H\otimes
\mu_{A\otimes H})\circ (\mu_{H}\otimes \mu_{H}\otimes (\rho_{A}\circ
 f^{-1})\otimes (\rho_{A}\circ f^{-1}))\circ
 (\delta_{H\otimes H}\otimes c_{H,H})\circ \delta_{H\otimes H}$

\item[ ]

\item[ ]$=(\mu_{A}\otimes \mu_{H})\circ (\mu_{A}\otimes A
\otimes \mu_{H}\otimes H)\circ (A\otimes A\otimes c_{H,A}\otimes
\mu_{H}\otimes H)\circ (f\otimes c_{H,A}\otimes c_{H,A}\otimes H\otimes H)$

\item[ ]

\item[ ]$\hspace{0.38cm}\circ
(\mu_{H}\otimes H\otimes c_{H,A}\otimes c_{H,A}\otimes H)\circ
(H\otimes c_{H,H}\otimes H\otimes (\rho_{A}\circ
 f^{-1})\otimes (\rho_{A}\circ f^{-1}))\co (\delta_{H}\otimes H\otimes
 \delta_{H}\otimes H)$

\item[ ]

\item[ ]$\hspace{0.38cm}\circ(H\otimes H\otimes c_{H,H})\circ
\delta_{H\otimes H}$

\item[ ]

\item[ ]$= (\mu_{A}\otimes \mu_{H})\circ (f\otimes
\mu_{A\otimes H}\otimes H)\circ (H\otimes c_{H,A}\otimes
c_{H,A}\otimes H)$

\item[ ]

\item[ ]$\hspace{0.38cm}\co(\mu_{H}\otimes H\otimes [\Gamma_{A}^{H}\circ (H\otimes f^{-1})\circ
\delta_{H}]\otimes (\rho_{A}\circ f^{-1}))\circ (H\otimes
c_{H,H}\otimes c_{H,H})\circ (\delta_{H}\otimes H\otimes H\otimes
H)\circ \delta_{H\otimes H} $

\item[ ]

\item[ ]$= (\mu_{A}\otimes \mu_{H})\circ (f\otimes
\mu_{A\otimes H}\otimes H)\circ (H\otimes c_{H,A}\otimes
c_{H,A}\otimes H) $

\item[ ]

\item[ ]$\hspace{0.38cm}\co (\mu_{H}\otimes H\otimes [(A\otimes
\overline{\Pi}_{H}^{R})\circ \rho_{A}\circ f^{-1}]\otimes
(\rho_{A}\circ f^{-1}))\circ (H\otimes c_{H,H}\otimes c_{H,H}) \circ
(\delta_{H}\otimes H\otimes H\otimes H)\circ \delta_{H\otimes H}$

\item[ ]

\item[ ]$= (\mu_{A}\otimes \mu_{H})\circ (A\otimes \mu_{A}\otimes
H\otimes (\mu_{H}\circ ( \overline{\Pi}_{H}^{R}\otimes H)))\circ
(f\otimes A\otimes c_{H,A}\otimes H\otimes H)$

\item[ ]

\item[ ]$\hspace{0.38cm} \co(\mu_{H}\otimes c_{H,A}\otimes
c_{H,A}\otimes H)\circ (H\otimes H\otimes H\otimes (\rho_{A}\circ
 f^{-1})\otimes (\rho_{A}\circ f^{-1}))\circ (H\otimes H\otimes
 H\otimes c_{H,H})$

 \item[ ]

\item[ ]$\hspace{0.38cm}\co (H\otimes H\otimes
\delta_{H}\otimes H)\circ \delta_{H\otimes H} $

\item[ ]

\item[ ]$=(\mu_{A}\otimes H)\circ (f\otimes [(A\otimes
\mu_{H})\circ (c_{H,A}\otimes (((\varepsilon_{H}\circ
\mu_{H})\otimes H)\circ (H\otimes\delta_{H})))\circ (H\otimes
\mu_{A}\otimes H\otimes H)$

\item[ ]

\item[ ]$\hspace{0.38cm}\co (H\otimes A\otimes c_{H,A}\otimes H)
\circ (H\otimes (\rho_{A}\circ
 f^{-1})\otimes (\rho_{A}\circ f^{-1}))\circ (H\otimes c_{H,H})\circ
 (\delta_{H}\otimes H)])$

 \item[ ]

\item[ ]$\hspace{0.38cm}\co (\mu_{H}\otimes H\otimes H)\circ
 \delta_{H\otimes H}$

\item[ ]

\item[ ]$=(\mu_{A}\otimes H)\circ (f\otimes [(A\otimes
\mu_{H})\circ  (c_{H,A}\otimes H)\circ (H\otimes \mu_{A}\otimes
H)\circ (H\otimes f^{-1}\otimes (\rho_{A}\circ f^{-1}))$

\item[ ]

\item[ ]$\hspace{0.38cm} \co(H\otimes c_{H,H})\circ
(\delta_{H}\otimes H) ])\circ (\mu_{H}\otimes H\otimes H)\circ
 \delta_{H\otimes H}$

\item[ ]

\item[ ]$=(\mu_{A}\otimes H)\co (f\ot [((\mu_A\co c_{A,A}\co (f^{-1}\ot f^{-1}))\ot \Pi_{H}^{L})\co (H\ot c_{H,H})\co ((c_{H,H}\ot \delta_H)\ot H)])$

\item[ ]

\item[ ]$\hspace{0.38cm}  \co (\mu_H\ot H\ot H)\co \delta_{H\ot H}$

\item[ ]

\item[ ]$=(\mu_{A}\otimes H)\co (f\ot [(A\otimes
\mu_{H})\circ  (c_{H,A}\otimes H)\circ (H\otimes \mu_{A}\otimes
H)\circ (H\otimes f^{-1}\otimes (\rho_{A}\circ f^{-1}))\co (H\otimes
c_{H,H})\circ (\delta_{H}\otimes H)])$

\item[ ]

\item[ ]$\hspace{0.38cm}  \co (\mu_{H}\otimes H\otimes H)\circ
 \delta_{H\otimes H}$

\item[ ]

\item[ ]$=(A\otimes \overline{\Pi}_{H}^{R})\circ
\rho_{A}\circ \sigma^{-1}_{A}.$

\end{itemize}

In the last computations, the first equality follows by the definitions,
the second one because $A$ is a right $H$-comodule algebra
and $H$ a weak Hopf algebra, the third one by the coassociativity  of
$\delta_{H}$, the fourth one by the associativity of $\mu_{A}$ and
$\mu_{H}$ and the fifth one by Theorem \ref{equivalence-cleft-weak-cleft}. In the
sixth equality  we use the associativity of $\mu_{H}$ and the seventh one
follows by (\ref{PiRbarramu}). To prove the eight equality we use the
definition of right $H$-comodule algebra and in the ninth one we
apply the  identity
\begin{equation}
\label{apoyo-dem} ((\mu_A\co c_{A,A}\co (f^{-1}\ot f^{-1}))\ot
\Pi_{H}^{L})\co (H\ot c_{H,H})\co ((c_{H,H}\co \delta_H)\ot H)
\end{equation}
$$=(A\otimes
\mu_{H})\circ  (c_{H,A}\otimes H)\circ (H\otimes \mu_{A}\otimes
H)\circ (H\otimes f^{-1}\otimes (\rho_{A}\circ f^{-1}))\co (H\otimes
c_{H,H})\circ (\delta_{H}\otimes H)$$

obtained in Proposition \ref{propinversa}. Finally, the tenth one is
obtained by (\ref{composiciones2}) and the idempotent character of $\overline{\Pi}_{H}^{R}$ and
the eleventh one by repetition of the previous computations but in
inverse order.

Let $\sigma_{A_{H}}^{-1}$ be the factorization of $\sigma_{A}^{-1}$.
We will finish the proof showing that $\sigma_{A_{H}}$ is a morphism
in $Reg_{\varphi_{A_{H}}}(H\ot H,A_{H})$ with inverse
$\sigma_{A_{H}}^{-1}$. First of all, note that $A_{H}\hookrightarrow
A$ is an $H$-cleft extension and then, by Proposition
\ref{siempreCleft}, the morphism $f\wedge f^{-1}$ factors through
the equalizer $i_A$. Now, using that $H$ is a weak Hopf algebra, $f$
an integral, $A$ an $H$-comodule algebra, $A_H$ a weak $H$-module
algebra and the equality (\ref{equcross2}) we obtain

\begin{itemize}

\item[ ]$\hspace{0.38cm} i_A\co (\sigma_{A_{H}}\wedge \sigma_{A_{H}}^{-1})$

\item[ ]

\item[ ]$=\mu_A\co (\mu_A\ot \mu_A)\co (((\mu_A\co (f\ot f))\ot((f^{-1}\wedge f)\co \mu_H)\ot ((f^{-1}\ot f^{-1})\co c_{H,H}))\co (H\ot H\ot \delta_{H\ot H})\co \delta_{H\ot H}$

\item[ ]

\item[ ]$=\mu_A\co (\mu_A\ot \mu_A)\co (A\ot (f^{-1}\wedge f)\ot ((f^{-1}\ot f^{-1})\co c_{H,H}))\co (\mu_{A\ot H}\ot H\ot H)$

\item[ ]

\item[ ]$\hspace{0.38cm}\co ((\rho_A\co f)\ot (\rho_A\co f)\ot H\ot H)\co \delta_{H\ot H}$

\item[ ]

\item[ ]$=\mu_A\co ((\mu_A\co (A\ot (f^{-1}\wedge f))\co \rho_A)\ot \mu_{A})\co ((\mu_{A}\co (f\ot f))\ot ((f^{-1}\ot f^{-1})\co c_{H,H}))\co \delta_{H\ot H}$

\item[ ]

\item[ ]$=q_A\co \mu_{A}\co(f\ot (f\wedge f^{-1}))$

\item[ ]

\item[ ]$=q_A\co \mu_{A}\co(f\ot (q_{A}\co \mu_{A}\co (f\ot (i_{A}\co\eta_{A_{H}}))))$

\item[ ]

\item[ ]$=i_A\co \varphi_{A_H}\co(H\ot (\varphi_{A_H}\co(H\ot \eta_{A_H}))$

\item[ ]

\item[ ]$=i_A\co \varphi_{A_H}\co(\mu_H\ot \eta_{A_H}),$

\end{itemize}

and then, using that $i_A$ is a monomorphism, $\sigma_{A_{H}}\wedge \sigma_{A_{H}}^{-1}=\varphi_{A_H}\co(\mu_H\ot \eta_{A_H})$.

To prove
$\sigma_{A_{H}}^{-1}\wedge\sigma_{A_{H}}=\varphi_{A_{H}}\co(\mu_H\co (H\ot \eta_{A_{H}}))$ we need to check the equality
\begin{equation}
\label{new-aux-2}
e_A\circ\mu_{H}\circ ( \Pi_{H}^{R}\otimes
H)=e_A\circ\mu_{H}
\end{equation}
where $e_A$ is the morphism defined in (\ref{e-gamma}). Indeed, using the equality (\ref{PiRmu}) and the definition of
 weak Hopf algebra we have:

\begin{itemize}

\item[ ]$\hspace{0.38cm}  e_A\circ\mu_{H}\circ
(\Pi_{H}^{R}\otimes H))$

\item[ ]

\item[ ]$=(A\otimes [((\varepsilon_{H}\circ \mu_{H})\otimes
 (\varepsilon_{H}\circ \mu_{H}))\circ (H\otimes
 (c_{H,H}\circ \delta_{H})\otimes H)]) \circ (c_{H,A}\otimes
 H\otimes H)\circ
 (H\otimes c_{H,A}\otimes H)$

 \item[ ]

\item[ ]$\hspace{0.38cm}\co (H\otimes H\otimes (\rho_{A}\circ \eta_{A}))$

\item[ ]

\item[ ]$= (A\otimes [\varepsilon_{H}\circ \mu_{H}\circ
(\mu_{H}\otimes H)]) \circ (c_{H,A}\otimes H\otimes H)\circ
 (H\otimes c_{H,A}\otimes H)\circ
 (H\otimes H\otimes (\rho_{A}\circ \eta_{A}))$

\item[ ]

\item[ ]$= e_A\circ\mu_{H}.$

\end{itemize}

Then, $\sigma^{-1}_{A_{H}} \wedge
\sigma_{A_{H}}=\varphi_{A_H}\co(\mu_H\ot \eta_{A_H})$, because by composing with
the equalizer $i_{A}$ we have

\begin{itemize}

\item[ ]$\hspace{0.38cm}i_{A}\co (\sigma^{-1}_{A_{H}} \wedge
\sigma_{A_{H}})$

\item[ ]

\item[ ]$=q_{A}\circ \mu_{A}\circ (A\otimes \mu_{A})\circ
(A\otimes\mu_{A}\otimes \mu_{A}) \circ (f\otimes f^{-1}\otimes
f^{-1}\otimes f\otimes f)\circ (\mu_{H}\otimes c_{H,H}\otimes
H\otimes H)$

\item[ ]

\item[ ]$\hspace{0.38cm}\co (\delta_{H\otimes H}\otimes H\otimes H)\circ
\delta_{H\otimes H} $

\item[ ]

\item[ ]$= q_{A}\circ \mu_{A}\circ (f\otimes \mu_{A}) \circ
(\mu_{H}\otimes(\mu_{A}\circ (f^{-1}\otimes f^{-1}))\otimes
(\mu_{A}\circ (f\otimes f)))\circ (H\otimes H\otimes
c_{H,H}\otimes H\otimes H)$

\item[ ]

\item[ ]$\hspace{0.38cm}\co (H\otimes c_{H,H}\otimes c_{H,H}\otimes H)\circ
(H\otimes \delta_{H\otimes H}\otimes H)\circ (\delta_{H}\otimes
\delta_{H}) $

\item[ ]

\item[ ]$=q_{A}\circ \mu_{A}\circ ((f\circ \mu_{H})\otimes
(\mu_{A}\circ (f^{-1}\otimes A)))\circ (H\otimes \delta_{H}\otimes
[\mu_{A}\circ (e_A\otimes f)])\circ \delta_{H\otimes H}$

\item[ ]

\item[ ]$=q_{A}\circ \mu_{A}\circ ((f\circ \mu_{H})\otimes
(\mu_{A}\circ (f^{-1}\otimes A)))\circ (H\otimes \delta_{H}
\otimes [ (A\otimes \varepsilon_{H})\circ \Gamma_{A}^{H}\circ (H\otimes
f)])\circ \delta_{H\otimes H}$

\item[ ]

\item[ ]$= q_{A}\circ \mu_{A}\circ ((f\circ \mu_{H})\otimes
(\mu_{A}\circ (f^{-1}\otimes f))\otimes (\varepsilon_{H}\circ
\mu_{H}))\circ (H\otimes \delta_{H}\otimes c_{H,H}\otimes H)$

\item[ ]

\item[ ]$\hspace{0.38cm}\co (H\otimes H\otimes H\otimes
\delta_{H})\circ \delta_{H\otimes H}$

\item[ ]

\item[ ]$= q_{A}\circ \mu_{A}\circ ((f\circ \mu_{H})\otimes
e_A\otimes (\varepsilon_{H}\circ \mu_{H}))\circ (H\otimes
H\otimes c_{H,H}\otimes H)\circ (H\otimes H\otimes H\otimes
\delta_{H})\circ \delta_{H\otimes H}$

\item[ ]

\item[ ]$= q_{A}\circ \mu_{A}\circ ((f\circ \mu_{H})\otimes
(e_A\circ\mu_{H}\circ (\Pi_{H}^{R}\otimes H)))\circ
\delta_{H\otimes H}$

\item[ ]

\item[ ]$=q_{A}\circ (f\wedge e_A)\circ \mu_{H} $

\item[ ]

\item[ ]$= q_{A}\circ  f\circ \mu_{H}.$

\item[ ]

\item[ ]$= (f\wedge f^{-1})\co  \mu_{H}.$

\item[ ]

\item[ ]$=i_A\co \varphi_{A_H}\co(\mu_H\ot \eta_{A_H}),$

\end{itemize}

The first equality follows from (\ref{sigma-f}) and (\ref{equcross5}), the second
one relies on the coassociativity of $\delta_{H}$ and the associativity of
$\mu_{A}$, the third one by the naturality of the braiding, the
fourth one follows from the properties of
$\Gamma_{A}^{H}$ and the fifth one by $(A\otimes \varepsilon_{H})\circ
\Gamma_{A}^{H}\circ (H\otimes f)=(f\otimes (\varepsilon_{H}\circ
\mu_{H}))\circ (c_{H,H}\otimes H)\circ (H\otimes \delta_{H})$. The
naturality of the braiding yields the sixth one and the equality
(\ref{PiRmu})
implies the seventh one. The eight one is a consequence of
(\ref{new-aux-2})
and the fact that $H$ is a weak Hopf algebra. Finally, in the ninth
one we use that $f\wedge e_A=f$, the tenth one follows because $f$ is a morphism of right
$H$-comodules and
in the eleventh one  we use the definition
of $\varphi_{A_{H}}$.

To prove (f2) we compose with the equalizer $i_A$

\begin{itemize}

\item[ ]$\hspace{0.38cm} i_A\co (\sigma_{A_{H}}\wedge \sigma_{A_{H}}^{-1}\wedge \sigma_{A_{H}})$

\item[ ]

\item[ ]$=i_A\co (\sigma_{A_{H}}\wedge (\varphi_{A_H}\co(\mu_{H}\ot \eta_{A_H}))$

\item[ ]

\item[ ]$=(\mu_A\co (f\ot f))\wedge (f^{-1}\co \mu_H)\wedge ((f\wedge f^{-1})\co \mu_H)$

\item[ ]

\item[ ]$=(\mu_A\co (f\ot f))\wedge ((f^{-1}\wedge f\wedge f^{-1})\co \mu_H)$

\item[ ]

\item[ ]$=i_A\co \sigma_{A_{H}},$

\end{itemize}

and then $\sigma_{A_{H}}\wedge \sigma_{A_{H}}^{-1}\wedge \sigma_{A_{H}}=\sigma_{A_{H}}$.
In a similar way, using that $f\wedge f^{-1}\wedge f=f$ we get (f3).

To finish the proof, we only need to show that $\sigma_{A_{H}}$ satisfies the normal condition.
Indeed: by the usual arguments
$$i_{A}\co \sigma_{A_{H}}\co (\eta_{H}\ot H)=q_{A}\co \mu_{A}\co (\eta_{A}\ot f)=q_{A}\co f=f\wedge f^{-1}=
i_{A}\co \varphi_{A_{H}}\co (H\ot \eta_{A_{H}})$$
and
$$i_{A}\co \sigma_{A_{H}}\co (H\ot \eta_{H})=q_{A}\co f=f\wedge f^{-1}=
i_{A}\co \varphi_{A_{H}}\co (H\ot \eta_{A_{H}}).$$ Therefore,
$\sigma_{A_{H}}\co (\eta_{H}\ot H)=\sigma_{A_{H}}\co (H\ot
\eta_{H})= \varphi_{A_{H}}\co (H\ot \eta_{A_{H}})$.

\end{dem}

\begin{cor}
\label{casosmash}

In the conditions of Proposition \ref{construccioncociclo}, the following assertions are equivalent:

\begin{itemize}

\item[(i)] $\sigma_{A_H}=\varphi_{A_H}\co (\mu_H\ot \eta_{A_H}).$

\item[(ii)] $\mu_A\co (f\ot f)=f\co \mu_H.$

\end{itemize}

\begin{dem}
{\rm
$(i)\Rightarrow (ii)$. Composing with $i_A$, and using the definition of $\sigma_{A}$

 \begin{itemize}

\item[ ]

\item[ ]$\hspace{0.38cm}(i_A\co \sigma_{A_H})\wedge (f\co \mu_H)$

\item[ ]

\item[ ]$=(\mu_A\co (f\ot f))\wedge (f^{-1}\co \mu_H) \wedge (f\co \mu_H)$

\item[ ]

\item[ ]$=(\mu_A\co (f\ot f))\wedge u_2$

\item[ ]

\item[ ]$=\mu_A\co (A\ot e_A)\co \rho_A\co \mu_A\co (f\ot f)$

\item[ ]

\item[ ]$=\mu_A\co (f\ot f)$

\end{itemize}

On the other hand, by the hypothesis, $(i_A\co \sigma_{A_H})\wedge (f\co \mu_H)=u_2\wedge (f\co \mu_H)=(u_1\wedge f)\co \mu_H=f\co \mu_H$ and we obtain (ii).

The converse  is an easy consequence of (\ref{sigma-f}). Indeed:

\begin{itemize}

\item[ ]

\item[ ]$\hspace{0.38cm}i_A\co \sigma_{A_H}$

\item[ ]

\item[ ]$=q_A\co \mu_A\co (f\ot f)$

\item[ ]

\item[ ]$=(f\wedge f^{-1})\co \mu_H$

\item[ ]

\item[ ]$=u_2$

\item[ ]

\item[ ]$=i_A\co \varphi_{A_H}\co (\mu_H\ot \eta_{A_H})$

\end{itemize}

and then $\sigma_{A_H}=\varphi_{A_H}\co (\mu_H\ot \eta_{A_H}).$
}
\end{dem}

\end{cor}

In the next Theorem we prove that each $H$-cleft extension
 determines an unique equivalence class of crossed systems for $H$ over $A$. First we need a fundamental result in the
 study of $H$-cleft extensions that generalizes Theorem 11 of \cite{doi3}.

\begin{teo}
\label{cleftimplicasistemacruz-previo}
Let $H$ be a cocommutative weak Hopf algebra and let
$A_{H}\hookrightarrow A$ be an
$H$-cleft extension  with
$f$ an associated convolution invertible total integral. Then, the pair $(\varphi_{A_H}, \sigma_{A_H})$ is a crossed system for $H$ over $A_{H}$,
where $\varphi_{A_H}$ is the weak $H$-module structure defined
in Proposition \ref{estructuramoduloinducida} and $\sigma_{A_H}$
the morphism obtained in Proposition \ref{construccioncociclo}. Moreover, the $H$-cleft extensions
$A_{H}\hookrightarrow A$ and $A_{H}\hookrightarrow A_{H}\times_{\sigma_{A_{H}}}H$ are equivalent.

\end{teo}

\begin{dem}
First note that in this case
$$\psi_{H}^{A_{H}}=(p_{A}\ot H)\co \rho_{A}\co \mu_{A}\co (f\ot i_{A})$$
and
$$\sigma_{H}^{A_{H}}=(p_{A}\ot H)\co \rho_{A}\co \mu_{A}\co (f\ot f)$$
and then, by Proposition 3.13 of \cite{nmra4}, we have that the quadruple
$(A_{H}, H, \psi_{H}^{A_{H}}, \sigma_{H}^{A_{H}})$ satisfies the twisted and the
cocycle conditions. Moreover, the normal
condition for $\sigma_{A_{H}}$ implies that there exists a preunit. Therefore, by the theory
exposed in \ref{weak-crossed-products-exposition}, we obtain that $(\varphi_{A_H}, \sigma_{A_H})$
is a crossed system for $H$ over $A_{H}$. Moreover, by Lemma 3.11 of \cite{nmra4}, we obtain that
\begin{equation}
\label{newexpressionfornabla}
\nabla_{A_{H}\ot H}=(p_{A}\ot H)\co \rho_{A}\co \mu_{A}\co (i_{A}\ot f).
\end{equation}

By Proposition \ref{sistemacruzadocleft}, we known that $A_{H}\hookrightarrow A_{H}\times_{\sigma_{A_{H}}}H$ is an $H$-cleft extension. Also, by 3.10 of \cite{nmra4}, there exists a right $H$-comodule algebra isomorphism
$$T=p_{A_{H}\ot H}\co (p_A\ot H)\co \rho_{A}:A\rightarrow A_{H}\times_{\sigma_{A_{H}}}H$$
such that
$$T^{-1}=\mu_{A}\co (i_{A}\ot f)\co i_{A_{H}\ot H}$$
and
$$T^{-1}\co i_{A_{H}\times_{\sigma_{A_{H}}}H}=\mu_{A}\co (i_{A}\ot f)\co \nabla_{A_{H}\ot H}\co (A_{H}\ot \eta_{H})=\mu_{A}\co (i_{A}\ot (f\co \eta_{H}))=i_{A}.$$

Therefore, we obtain that $A_{H}\hookrightarrow A$ and $A_{H}\hookrightarrow A_{H}\times_{\sigma_{A_{H}}}H$ are equivalent.

\end{dem}

\begin{prop}
\label{third-step}
Let $H$ be a cocommutative weak Hopf algebra and let
$(\varphi_{A}, \sigma)$ be a crossed system for $H$ over $A$. Let $A\hookrightarrow
A\times_{\sigma}H$ be
the $H$-cleft extension constructed in Proposition \ref{sistemacruzadocleft}. Then,
if $(\phi_{A}, \tau)$ is the crossed system associated to the $H$-cleft extension
 $A\hookrightarrow A\times_{\sigma}H$, we have that
$(\phi_{A}, \tau)=(\varphi_{A}, \sigma)$.
\end{prop}

\begin{dem} By Proposition \ref{H-comod-alg} and
Theorem \ref{cleftimplicasistemacruz-previo},
the convolution invertible  total integral $f=p_{A\otimes H}\co (\eta_A\ot H)$ determines a
crossed system $(\phi_A, \tau)$, where $\phi_A$ and $\tau$ are defined by
\begin{equation}
\label{weakHmoduloparah}
\phi_A=p_{A\times_{\sigma} H}\co \mu_{A\times_{\sigma} H}\co (f\ot i_{A\times_{\sigma} H})
\end{equation}
and
\begin{equation}
\label{cocicloparah}
\tau=p_{A\times_{\sigma} H}\co \mu_{A\times_{\sigma} H}\co (f\ot f),
\end{equation}

where $p_{A\times_{\sigma} H}$ is the factorization through the equalizer
$i_{A\times_{\sigma} H}=p_{A\otimes H}\co (A\ot \eta_H)$ of the morphism
$q_{A\times_{\sigma} H}=\mu_{A\times_{\sigma} H}\co (A\times_{\sigma} H\ot f^{-1})\co \rho_{A\times_{\sigma} H}$.

We will whow that the crossed systems $(\varphi_A, \sigma)$ and
$(\phi_A, \tau)$ coincide. First of all we prove the following
equality:
$$u_{1}$$
\begin{equation}
\label{ecuacion} =\mu_A\co ((\varphi_A\co (H\ot \sigma^{-1}))\ot
\sigma)\co (H\ot H\ot ((c_{H,H}\ot H)\co ( H\ot  c_{H,H})))\co
((\delta_{H\ot H}\co (H\ot \lambda_{H})\co \delta_{H})\ot H) \co
\delta_H.
\end{equation}

Indeed:

\begin{itemize}

\item[ ]

\item[ ]$\hspace{0.38cm}\mu_A\co ((\varphi_A\co (H\ot \sigma^{-1}))\ot \sigma)\co (H\ot H\ot
((c_{H,H}\ot H)\co ( H\ot  c_{H,H})))\co
((\delta_{H\ot H}\co (H\ot \lambda_{H})\co \delta_{H})\ot H) \co \delta_H$

\item[ ]

\item[ ]$=\mu_A\co ((\varphi_A\co (H\ot \sigma^{-1}))\ot \sigma)\co (H\ot H\ot c_{H,H}\ot (\lambda_H\wedge \Pi^{L}_{H}))\co (H\ot c_{H,H}\ot c_{H,H})$

\item[ ]

\item[ ]$\hspace{0.38cm} \co(\delta_H\ot \lambda_H\ot H\ot H)\co (H\ot (c_{H,H}\co \delta_H)\ot H)\co (H\ot \delta_H)\co \delta_H$

\item[ ]

\item[ ]$=\mu_A\co ((\varphi_A\co (H\ot \sigma^{-1}))\ot (\sigma\co (H\ot \mu_H)))\co(H\ot H\ot c_{H,H}\ot H\ot \Pi^{L}_{H})$

\item[ ]

\item[ ]$\hspace{0.38cm} \co (H\ot c_{H,H}\ot c_{H,H}\ot H)\co (\delta_H\ot (c_{H,H}\co \delta_H\co \lambda_H)\ot (c_{H,H}\co \delta_H))\co (H\ot \delta_H)\co \delta_H$

\item[ ]

\item[ ]$=((\varphi_A\co (H\ot \sigma^{-1}))\wedge (\sigma\co (H\ot (\mu_H\co (H\ot \Pi^{L}_{H}))))\co (H\ot \lambda_H\ot H)\co (H\ot \delta_H)\co \delta_H$

\item[ ]

\item[ ]$=((\sigma\co (H\ot \mu_H))\wedge (\sigma^{-1}\co (\mu_H\ot H)))\co (H\ot \lambda_H\ot H)\co (H\ot \delta_H)\co \delta_H$

\item[ ]

\item[ ]$=\mu_A\co (\sigma\ot \sigma^{-1})\co (H\ot (c_{H,H}\co (\Pi^{L}_{H}\ot \Pi^{R}_{H}))\ot H)\co (\delta_H\ot \delta_H)\co \delta_H$

\item[ ]

\item[ ]$=\mu_A\co ((\sigma\co (H\ot\Pi^{R}_{H})\co \delta_H)\ot (\sigma^{-1}\co (\Pi^{L}_{H}\ot H)\co \delta_H))\co \delta_H$

\item[ ]

\item[ ]$=u_{1}\wedge u_{1}$

\item[ ]

\item[ ]$=u_{1}$

\end{itemize}

In the previous computations, the first, third and fifth equalities follow by
the properties of the antipode $\lambda_H$ and (\ref{propiedadesgeneralespiLR});
in the  second and sixth ones we use that $H$ is cocommutative; the fourth follows
 by (g2); the seventh is a consequence of (\ref{sigma-PiR-r}) and (\ref{sigma-PiL-l}); finally, in the last
 equality we use that $A$ is a weak $H$-module algebra.

Now we can obtain a simple expression for the morphism
$q_{A\times_{\sigma} H}$:

\begin{equation}
\label{expresionparaq}
q_{A\times_{\sigma} H}=p_{A\otimes H}\co (A\ot \Pi^{L}_{H})\co i_{A\otimes H}.
\end{equation}

\begin{itemize}

\item[ ]

\item[ ]$\hspace{0.38cm} q_{A\times_{\sigma} H}$

\item[ ]

\item[ ]$=\mu_{A\times_{\sigma} H}\co (p_{A\otimes H}\ot f^{-1})\co (A\ot \delta_H)\co i_{A\ot H}$

\item[ ]

\item[ ]$=p_{A\otimes H}\co (\mu_A\ot H)\co (\mu_A\ot \sigma\ot \mu_H)\co (A\ot \varphi_A\ot \delta_{H\ot H})\co (A\ot H\ot c_{H,A}\ot H)$

\item[ ]

\item[ ]$\hspace{0.38cm} \co (A\ot \delta_H\ot \sigma^{-1}\ot H)\co (A\ot H\ot H\ot c_{H,H})\co (A\ot H\ot (\delta_H\co \lambda_H)\ot H)\co (A\ot H\ot \delta_H)\co (A\ot \delta_H)\co i_{A\ot H}$

\item[ ]

\item[ ]$=p_{A\otimes H}\co (\mu_A\ot \Pi^{L}_{H})\co (A\ot \mu_A\ot H)\co (A\ot ((H\ot \varphi_A)\co \sigma^{-1})\ot \sigma\ot H)\co (A\ot H\ot H\ot c_{H,H}\ot H\ot H)$

\item[ ]

\item[ ]$\hspace{0.38cm} \co (A\ot H\ot c_{H,H}\ot c_{H,H}\ot H)\co
(A\ot \delta_H\ot (\delta_H\co \lambda_H)\ot H)\co (A\ot H\ot
(c_{H,H}\co \delta_H))$

\item[ ]

\item[ ]$\hspace{0.38cm} \co (A\ot \delta_H)\co i_{A\ot H}$

\item[ ]

\item[ ]$=p_{A\otimes H}\co (\mu_A\ot \Pi^{L}_{H})\co (A\ot
(\mu_A\co ((\varphi_A\co (H\ot \sigma^{-1}))\ot \sigma)\co (H\ot
H\ot c_{H,H}\ot H)$

\item[ ]

\item[ ]$\hspace{0.38cm} \co (H\ot c_{H,H}\ot c_{H,H})\co (\delta_H\ot
(\delta_H\co \lambda_H)\ot H)\co (H\ot \delta_H)\co \delta_H))\co
i_{A\ot H}$

\item[ ]

\item[ ]$=p_{A\otimes H}\co (\mu_A\ot \Pi^{L}_{H})\co (A\ot u_{1}\ot H)\co (A\ot \delta_H)\co i_{A\ot H}$

\item[ ]

\item[ ]$=p_{A\otimes H}\co (A\ot \Pi^{L}_{H})\co \nabla_{A\ot H}\co i_{A\ot H}$

\item[ ]

\item[ ]$=p_{A\otimes H}\co (A\ot \Pi^{L}_{H})\co i_{A\ot H}.$

\end{itemize}

In the foregoing calculations, the first equality follows by the $H$-comodule structure for
$A\times_{\sigma} H$; in the second one we use that $\mu_{A\ot_{\sigma} H}\co
(\nabla_{A\ot H}\ot \nabla_{A\ot H})=\mu_{A\ot_{\sigma} H}$, the third relies
on the antimultiplicativity of the antipode; the fourth by cocommutativity of $H$.
The fifth one follows by (\ref{ecuacion}), the seventh ia a consequence of the definition
of $\nabla_{A\ot H}$; finally the last one uses that $\nabla_{A\ot H}\co i_{A\ot H}= i_{A\ot H}$.

On the other hand, $q_{A\times_{\sigma} H}=i_{A\times_{\sigma} H}\co p_{A\times_{\sigma} H}
=p_{A\otimes H}\co (p_{A\times_{\sigma} H}\ot \eta_H)$ and then

$$(A\ot \varepsilon_H)\co i_{A\ot H}\co q_{A\times_{\sigma} H}=
(A\ot \varepsilon_H)\co \nabla_{A\ot H}\co (p_{A\times_{\sigma} H}\ot \eta_H)=p_{A\times_{\sigma} H}$$

As a consequence,
\begin{equation}
\label{expresionparap}
p_{A\times_{\sigma} H}=(A\ot \varepsilon_H)\co i_{A\ot H}
\end{equation}
because
$$p_{A\times_{\sigma} H}=(A\ot \varepsilon_H)\co \nabla_{A\ot H}\co (A\ot \Pi^{L}_{H})
\co i_{A\ot H}=(\mu_A\ot \varepsilon_H)\ot (A\ot (\varphi_A\co (\Pi^{L}_{H}\ot \eta_A))\ot H)
\co i_{A\ot H}$$
$$=(A\ot \varepsilon_H)\co \nabla_{A\ot H}\co i_{A\ot H}=
(A\ot \varepsilon_H)\co i_{A\ot H}.$$

Using this equality it is not difficult to see that $(\varphi_A,
\sigma)=(\phi_A, \tau)$. Indeed:

\begin{itemize}

\item[ ]

\item[ ]$\hspace{0.38cm} \phi_A$

\item[ ]

\item[ ]$=p_{A\times_{\sigma} H}\co \mu_{A\times_{\sigma} H}\co (f\ot i_{A\times_{\sigma} H})$

\item[ ]

\item[ ]$=(\mu_A\ot \varepsilon_H)\co (\varphi_A\ot \sigma\ot \mu_H)\co (H\ot A\ot \delta_{H\ot H})\co (H\ot c_{H,A}\ot H)\co (\delta_H\ot A\ot \eta_H)$

\item[ ]

\item[ ]$=\mu_A\co (\varphi_A\ot (\sigma\co (H\ot\Pi_{H}^{R})\co \delta_H))\co (H\ot c_{H,A})\co (\delta_H\ot A)$

\item[ ]

\item[ ]$=\mu_A\co (\varphi_A\ot u_1)\co (H\ot c_{H,A})\co (\delta_H\ot A)$

\item[ ]

\item[ ]$=\varphi_A$,

\end{itemize}

and

\begin{itemize}

\item[ ]

\item[ ]$\hspace{0.38cm} \tau$

\item[ ]

\item[ ]$=p_{A\times_{\sigma} H}\co \mu_{A\times_{\sigma} H}\co (f\ot f)$

\item[ ]

\item[ ]$=(A\ot \varepsilon_H)\co  \nabla_{A\ot H}\co (\sigma\ot \mu_H)\co \delta_{H\ot H}$

\item[ ]

\item[ ]$=\sigma\wedge u_{2}$

\item[ ]

\item[ ]$=\sigma$.

\end{itemize}

\end{dem}

The following theorem is the weak version of Lemma 2.1 of \cite{doi1}.

\begin{prop}
\label{cleftimplicasistemacruz}
Let $H$ be a cocommutative weak Hopf algebra and let
$A_{H}\hookrightarrow A$ be an
$H$-cleft extension  with
$f$ an associated convolution invertible total integral. Assume that $g:H\rightarrow A$ is another convolution
invertible total integral  with associated crossed system $(\phi_{A_H}, \tau_{A_H})$. Then
the crossed systems $(\varphi_{A_H}, \sigma_{A_H})$ and $(\phi_{A_H}, \tau_{A_H})$
are equivalent.

\end{prop}

\begin{dem}
The morphism $\widetilde{h}=f\wedge g^{-1}$ factors through the equalizer $i_A$. Indeed, by (\ref{cleaving-lambda}), the coassociativity of $\delta_{H}$ and the naturality of $c$, we have

\begin{itemize}

\item[ ]$\hspace{0.38cm}\rho_A\co \widetilde{h}$

\item[ ]

\item[ ]$=(\rho_A\co f)\wedge (\rho_A\co g^{-1})$

\item[ ]

\item[ ]$=((f\ot H)\co \delta_H)\wedge ((g^{-1}\ot \lambda_H)\co c_{H,H}\co \delta_H)$

\item[ ]

\item[ ]$=((\mu_{A}\co (f\ot g^{-1}))\ot \Pi_{H}^{L})\co (H\ot (c_{H,H}\co \delta_H))\co \delta_H$

\item[ ]

\item[ ]$=((\mu_{A}\co (f\ot g^{-1}))\ot (\Pi_{H}^{L}\co \Pi_{H}^{L}))\co (H\ot (c_{H,H}\co \delta_H))\co \delta_H$

\item[ ]

\item[ ]$=(A\ot \Pi_{H}^{L})\co \rho_A\co \widetilde{h},$

\end{itemize}

and then there exists a morphism $h:H\rightarrow A_H$ such that $\widetilde{h}=i_A\co h$. Note that, in the conditions of this theorem, $f\wedge f^{-1}=g\wedge g^{-1}$,
and $f^{-1}\wedge f=g^{-1}\wedge g$. Then,
\begin{equation}
\label{interseccion}
\varphi_{A_{H}}\co (H\ot \eta_{A_{H}})=\phi_{A_{H}}\co
(H\ot \eta_{A_{H}})
\end{equation}
because
$$i_{A}\co \varphi_{A_{H}}\co (H\ot \eta_{A_{H}})=q_{A}\co f=f\wedge f^{-1}$$
and
$$i_{A}\co \phi_{A_{H}}\co (H\ot \eta_{A_{H}})=q_{A}^{\prime}\co g=g\wedge g^{-1}$$
where $q_{A}$ is the morphism defined in Remark \ref{rem-simplify} and $q_{A}^{\prime}$ the analogous for $g$.

On the other hand, using that $f$ and $g$ are convolution invertible total integrals, we have
\begin{itemize}

\item[ ]$\hspace{0.38cm} \widetilde{h}\co \eta_H$

\item[ ]

\item[ ]$=\mu_A\co (A\ot g^{-1})\co \rho_A\co f\co \eta_H$

\item[ ]

\item[ ]$=\mu_A\co (A\ot g^{-1})\co \rho_A\co g\co \eta_H$

\item[ ]

\item[ ]$=(g\wedge g^{-1})\co \eta_H$

\item[ ]

\item[ ]$=\eta_A,$

\end{itemize}

Therefore, taking into account that $\eta_A=i_A\co \eta_{A_H}$, we obtain that $h\co \eta_H=\eta_{A_{H}}$.

The morphism $\widetilde{h}^{-1}=g\wedge f^{-1}$ admits a factorization
through the equalizer $i_A$ (the proof is similar to the one developed for $\widetilde{h}$)
and the factorization $h^{-1}$ is the convolution inverse
of $h$. As a consequence $h$ is in $Reg_{\varphi_{A_{H}}}(H,A_H)\cap Reg_{\phi_{A_{H}}}(H,A_H)$. Indeed:
First note that
$$ i_{A}\co (h\wedge h^{-1})= \widetilde{h}\wedge \widetilde{h}^{-1}=f\wedge g^{-1}\wedge
g\wedge f^{-1}=f\wedge f^{-1}\wedge f\wedge f^{-1}=\varphi_{A}\co (H\ot \eta_{A})=
i_{A}\co \varphi_{A_{H}}\co (H\ot \eta_{A_{H}})$$
and, by (\ref{interseccion}),  $h\wedge h^{-1}=\varphi_{A_{H}}\co (H\ot \eta_{A_{H}})=
\phi_{A_{H}}\co (H\ot \eta_{A_{H}}).$
Similarly, $h^{-1}\wedge h=\varphi_{A_{H}}\co (H\ot \eta_{A_{H}})=\phi_{A_{H}}\co
(H\ot \eta_{A_{H}}).$ Moreover,
$$i_{A}\co (h\wedge h^{-1}\wedge h)=\widetilde{h}\wedge \widetilde{h}^{-1}\wedge \widetilde{h}=\widetilde{h}=i_{A}\co h$$
and
$$i_{A}\co (h^{-1}\wedge h\wedge h^{-1})=\widetilde{h}^{-1}\wedge \widetilde{h}\wedge \widetilde{h}^{-1}=\widetilde{h}^{-1}=i_{A}\co h^{-1}$$
Then $h\wedge h^{-1}\wedge h=h$ and $h^{-1}\wedge h\wedge h^{-1}=h^{-1}$.

The proof for (\ref{relacionfis}) follows by the definition of $h$ and $\phi_{A_{H}}$. Indeed:

\begin{itemize}

\item[ ]$\hspace{0.38cm} i_A\co \mu_{A_H}\co (\mu_{A_H}\ot A_H)\co
(h\ot \phi_{A_H}\ot h^{-1})\co (\delta_H\ot c_{H,A_{H}})\co (\delta_H\ot A_H)$

\item[ ]

\item[ ]$=\mu_A\co (\mu_A\ot (f\wedge g^{-1}\wedge g))\co ((g^{-1}\wedge g\wedge f^{-1})\ot c_{H,A})\co (\delta_H\ot i_A)$

\item[ ]

\item[ ]$=\mu_A\co (\mu_A\ot (f^{-1}\wedge f\wedge f^{-1}))\co ((f\wedge f^{-1}\wedge f)\ot c_{H,A})\co (\delta_H\ot i_A)$

\item[ ]

\item[ ]$=i_A\co\varphi_{A_H},$

\end{itemize}

and then $\varphi_{A_H}=\mu_{A_H}\co (\mu_{A_H}\ot A_H)\co (h\ot \phi_{A_H}\ot h^{-1})\co (\delta_H\ot c_{H,A_{H}})\co (\delta_H\ot A_H).$

In order to get (\ref{relacionsigmas}), we begin by showing the equality
\begin{equation}\label{fygjuntas}
\mu_A\co \mu_{A\ot A}\co (f\ot (g^{-1}\wedge g)\ot (f\wedge g^{-1})\ot g)\co (\delta_H\ot \delta_H)=\mu_A\co (f\ot f),
\end{equation}

which follows because $f$ and $g$ are convolution invertible integrals, $A$ a right $H$-comodule algebra, (\ref{PiLmu}) and the equality $f^{-1}\wedge f=g^{-1}\wedge g$. Indeed:

\begin{itemize}

\item[ ]$\hspace{0.38cm} \mu_A\co \mu_{A\ot A}\co (f\ot (g^{-1}\wedge g)\ot (f\wedge g^{-1})\ot g)\co (\delta_H\ot \delta_H)$

\item[ ]

\item[ ]$=(A\ot \varepsilon_{H})\co \mu_{A\ot H}\co (\mu_{A}\ot H\ot ((\mu_{A}\ot H)\co (A\ot c_{H,A})\co ((\rho_{A}\co \eta_{A})\ot A)))\co (A\ot c_{H,A}\ot A)$

\item[ ]

\item[ ]$\hspace{0.38cm} \co((\mu_{A}\co (f\ot f))\ot H\ot ((g^{-1}\ot g)\co \delta_H))\co \delta_{H\ot H} $

\item[ ]

\item[ ]$=(A\ot \varepsilon_{H})\co \mu_{A\ot H}\co ((\mu_{A}\co (f\ot f))\ot H\ot ((\mu_{A}\ot \Pi_{H}^{L})\co (A\ot \rho_{A})\co (g^{-1}\ot g)\co \delta_H))\co \delta_{H\ot H} $

\item[ ]

\item[ ]$=(A\ot \varepsilon_{H})\co \mu_{A\ot H}\co ((\mu_{A}\co (f\ot f))\ot H\ot ((g^{-1}\wedge g)\ot \Pi_{H}^{L})\co \delta_{H})\co \delta_{H\ot H} $

\item[ ]

\item[ ]$=(A\ot \varepsilon_{H})\co \mu_{A\ot H}\co
(((f\ot H)\co \delta_{H})\ot (((f\wedge f^{-1}\wedge f)\ot H)\co \delta_{H}))$

\item[ ]

\item[ ]$=(A\ot \varepsilon_{H})\co \mu_{A\ot H}\co (((f\ot H)\co \delta_{H})\ot ((f\ot H)\co \delta_{H}))$

\item[ ]

\item[ ]$=\mu_A\co (f\ot f).$

\end{itemize}

Using this equality and similar arguments to the ones developed above we will finish the proof showing (\ref{relacionsigmas}):

\begin{itemize}

\item[ ]$\hspace{0.38cm} i_A\co\mu_{A_H}\co (\mu_{A_H}\ot h^{-1})\co (\mu_{A_H}\ot \tau_{A_{H}} \ot \mu_H)\co (h\ot \phi_{A_H}\ot \delta_{H\ot H})\co (\delta_H\ot h\ot H\ot H)\co \delta_{H\ot H}$

\item[ ]

\item[ ]$=\mu_A\co  ((\mu_{A}\co (\mu_{A}\ot g^{-1})\co ((f\wedge f^{-1}\wedge f)\ot
(f\wedge g^{-1})\ot H)\co (H\ot c_{H,H})\co (\delta_{H}\ot H))$

\item[ ]

\item[ ]$\hspace{0.38cm} \ot(\mu_{A}\co (A\ot
(f^{-1}\wedge f\wedge f^{-1}))\co \rho_{A}\co \mu_{A}\co (g\ot g)))\co \delta_{H\ot H}$

\item[ ]

\item[ ]$=(\mu_A\co \mu_{A\ot A}\co (f\ot (g^{-1}\wedge g)\ot (f\wedge g^{-1})\ot g)\co (\delta_H\ot \delta_H))\wedge (f^{-1}\co \mu_H)$

\item[ ]

\item[ ]$=(\mu_A\co (f\ot f))\wedge (f^{-1}\co \mu_H)$

\item[ ]

\item[ ]$=i_A\co \sigma_{A_{H}}.$

\end{itemize}

\end{dem}

\begin{cor}
\label{second-step} Let $H$ be a cocommutative weak Hopf algebra and let
$A_{H}\hookrightarrow A$, $A_{H}\hookrightarrow B$ two equivalent
$H$-cleft extensions with associated convolution invertible total
integrals $f$ and $g$ respectively. Then the corresponding  crossed systems
$(\varphi_{A_H}, \sigma_{A_H})$ and $(\phi_{A_H}, \tau_{A_H})$
are equivalent.
\end{cor}

\begin{dem} If
$A_{H}\hookrightarrow A$ and $A_{H}\hookrightarrow B$ are equivalent, there exists an isomorphism of right $H$-comodule algebras
$T:A\rightarrow B$ such that $i_{B}=T\co i_{A}$, and, as a consequence, $l=T\co f$ is a convolution invertible total
integral for $A_{H}\hookrightarrow B$ with inverse $l^{-1}=T\co f^{-1}$. Therefore, by Proposition
\ref{cleftimplicasistemacruz},  the crossed system $(\psi_{A_{H}}, \omega_{A_{H}})$ associated to
$A_{H}\hookrightarrow B$ for $l$ is equivalent to $(\phi_{A_H}, \tau_{A_H})$. Moreover, if
$p_{A}$ is the factorization through $i_{A}$ of the morphism $q_{A}=\mu_{A}\co (A\ot f^{-1})\co \rho_{A}$ and
$p_{B}$ is the factorization through $i_{B}$ of the morphism $q_{B}=\mu_{A}\co
(B\ot l^{-1})\co \rho_{B}$, we have the following: By (\ref{alternativamodulo})
$$\psi_{A_{H}}=p_{B}\co \mu_{B}\co (l\ot i_{B})$$
and then
$$\psi_{A_{H}}=p_{B}\co \mu_{B}\co (T\ot T)\co (f\ot i_{A})=p_{B}\co T\co \mu_{A}\co (f\ot i_{A})=
p_{A}\co \mu_{A}\co (f\ot i_{A})=\varphi_{A_{H}}.$$
On the other hand, by (\ref{sigma-f})
$$\omega_{A_{H}}=p_{B}\co \mu_{B}\co (l\ot l)$$
and, as a consequence,
$$\omega_{A_{H}}=p_{B}\co \mu_{B}\co (T\ot T)\co  (f\ot f)=p_{B}\co T\co \mu_{A}\co (f\ot f)=
p_{A}\co \mu_{A}\co (f\ot f)=\sigma_{A_{H}}.$$
Therefore $(\varphi_{A_H}, \sigma_{A_H})$ and $(\phi_{A_H}, \tau_{A_H})$
are equivalent.
\end{dem}

\begin{teo}
\label{siysolosiiso}
Let $H$ be a cocommutative weak Hopf algebra. Two $H$-cleft extensions $A_{H}\hookrightarrow A$,
$A_{H}\hookrightarrow B$ are equivalent if and only if so are their respective associated crossed
systems.
\end{teo}

\begin{dem} The "if" part is a consequence of the previous corollary. Moreover, if $A_{H}\hookrightarrow A$,
$A_{H}\hookrightarrow B$ are two weak $H$-cleft extensions with equivalent crossed systems
$(\varphi_{A_H}, \sigma_{A_H})$ and $(\phi_{A_H}, \tau_{A_H})$, by Proposition
\ref{first-step} we know that the associated $H$-cleft extensions
$A_H\hookrightarrow A_H\times_{\sigma_{A_H}}H$
and $A_H\hookrightarrow A_H\times_{\tau_{A_H}}H$ are equivalent. Therefore, by Theorem
\ref{cleftimplicasistemacruz-previo},
we obtain
$$A_{H}\hookrightarrow A \approx A_{H}\hookrightarrow A_{H}\times_{\sigma_{A_{H}}}H \approx
A_H\hookrightarrow A_H\times_{\tau_{A_H}}H \approx A_{H}\hookrightarrow B$$
which proves the Theorem.
\end{dem}

Now we can give the main result of this section which is a generalization of Theorem 2.7 of \cite{doi1}.

\begin{teo}
\label{correspondenciabiyectiva}
Let $H$ be a cocommutative weak Hopf algebra and $(A,\rho_{A})$  a right $H$-comodule algebra. There exists a
bijective correspondence between the equivalence  classes of $H$-cleft extensions
$A_{H}\hookrightarrow B$ and
the equivalence classes of crossed systems for $H$ over $A_{H}$.

\end{teo}

\begin{dem} If $CS(H,A_{H})$ denotes the set of equivalence classes of
crossed systems of $H$ over $A_{H}$ and $Cleft(A_{H})$ the set of equivalence classes of $H$-cleft extensions
$A_{H}\hookrightarrow B$, by Proposition \ref{first-step} and Corollary \ref{second-step} we
have two maps
$$F:CS(H,A_{H})\rightarrow Cleft(A_{H}),\;\;\; G:Cleft(A_{H})\rightarrow CS(H,A_{H})$$
defined by
$$F([(\varphi_{A_{H}}, \sigma_{A_{H}})])=[A_{H}\hookrightarrow A_{H}\times_{\sigma_{A_{H}}}H]$$
and
$$G([A_{H}\hookrightarrow B])=[(\phi_{A_H}, \tau_{A_H})].$$

The map $G$ is the inverse of $F$, because, by Proposition \ref{third-step}, we have
$$(G\co F)([(\varphi_{A_{H}},\sigma_{A_{H}})])=G([A_{H}\hookrightarrow A_{H}\times_{\sigma_{A_{H}}}H])=
[(\varphi_{A_{H}},\sigma_{A_{H}})],$$
and by Theorem \ref{cleftimplicasistemacruz-previo}
$$(F\co G)([A_{H}\hookrightarrow B])=F([(\phi_{A_{H}}, \tau_{A_{H}})])=
[A_{H}\hookrightarrow A_{H}\times_{\tau_{A_{H}}}H]=[A_{H}\hookrightarrow B].$$
\end{dem}

\section{Crossed systems and cohomology}

In \cite{NikaRamon6} we have developed a cohomology theory of algebras over weak Hopf algebras which generalizes the one given in \cite{Moss} for Hopf algebras. The main result contained in \cite{NikaRamon6} (see Theorem 3.11) asserts that if $(A,\varphi_{A})$ is a commutative left $H$-module algebra, there exists a bijection between the second cohomology group, denoted by $ H^{2}_{\varphi_{A}}(H,A)$, and the equivalence classes of  weak crossed products $A\ot_{\alpha} H$ where $\alpha:H\ot H\rightarrow A$ satisfies the 2-cocycle  and the normal conditions. In this section, for a cocommutative weak Hopf algebra and an $H$-cleft extension $A_{H} \hookrightarrow A$, we will establish a bijection between the set of equivalence classes of  crossed systems with a fixed weak $H$-module algebra structure and the second cohomology group $H_{\varphi_{\mathcal{Z}(A_{H})}}^2(H, Z(A_H))$, being  $Z(A_H)$ the center of the subalgebra of coinvariants $A_H$.  Our results generalizes to the weak Hopf algebra setting   the ones proved by Doi for Hopf algebras in  \cite{doi1}.

\begin{prop}
\label{centroesHmoduloalgebra}

Let $H$ be a cocommutative weak Hopf algebra and let $A_{H} \hookrightarrow A$ be an $H$-cleft extension. We denote by $(\varphi_{A_H}, \sigma_{A_H})$ the corresponding crossed system defined by the convolution invertible total integral $f:H\rightarrow A$. Then $(Z(A_H), \varphi_{Z(A_H)})$ is a left $H$-module algebra, where $\varphi_{Z(A_H)}$ is the factorization through the equalizer $i_{Z(A_H)}$ of the morphism $\varphi_{A_H}\co (H\ot i_{Z(A_H)})$.

\end{prop}

\begin{dem}

We define $\psi_A:H\ot A_H\rightarrow A$ as $\psi_A=\mu_A\co (A\ot (\mu_A\co c_{A,A}))\co (((f^{-1}\ot f)\co \delta_H)\ot i_A)$. In a similar way to Proposition \ref{estructuramoduloinducida}, it is not difficult to see that $\psi_A$ factors through the equalizer $i_A$, and then there exists a morphism $\psi_{A_H}:H\ot A_H\rightarrow A_H$ such that $i_A\co \psi_{A_H}=\psi_A$.
On the other hand, the following equalities hold:

\begin{equation}
\label{primeraigualdad}
\mu_A\co (f^{-1}\ot i_A)=\mu_A\co (\psi_A\ot f^{-1})\co (H\ot c_{H,A_H})\co (\delta_H\ot A_H).
\end{equation}

\begin{equation}
\label{segundaigualdad}
\mu_A\co c_{A,A}\co (f\ot i_A)=\mu_A\co (\psi_A\ot f)\co  (\delta_H\ot A_H).
\end{equation}

Indeed, using (e3), (e1) and that $u_1$ factors through the center of $A$ (which follows because $H$ is cocommutative and then (d4) and (d5) coincide),

\begin{itemize}

\item[ ]$\hspace{0.38cm} \mu_A\co (f^{-1}\ot i_A)$

\item[ ]

\item[ ]$=\mu_A\co ((f^{-1}\wedge u_1)\ot i_A)$

\item[ ]

\item[ ]$=\mu_A\co (\mu_{A}\ot A)\co (f^{-1}\ot i_A\ot u_1))\co (H\ot c_{H,A_H})\co (\delta_H\ot A_H)$

\item[ ]

\item[ ]$=\mu_A\co (\mu_{A}\ot A)\co (f^{-1}\ot i_A\ot (f\wedge f^{-1}))\co (H\ot c_{H,A_H})\co (\delta_H\ot A_H)$

\item[ ]

\item[ ]$=\mu_A\co (\psi_A\ot f^{-1})\co (H\ot c_{H,A_H})\co (\delta_H\ot A_H).$

\end{itemize}

The proof for (\ref{segundaigualdad}) follows a similar pattern.

 Now we can prove that $\varphi_{A_H}\co (H\ot i_{Z(A_H)})$ factors through the center of $A_H$. Indeed:

\begin{itemize}

\item[ ]

\item[ ]$\hspace{0.38cm} i_A\co \mu_{A_H}\co ((\varphi_{A_H}\co (H\ot i_{Z(A_H)}))\ot A_H)$

\item[ ]

\item[ ]$=\mu_A\co ((\mu_A\co(f\ot i_A))\ot (\mu_A\co (f^{-1}\ot i_A)))\co (H\ot c_{H,A_H}\ot A_H)\co (\delta_H\ot i_{Z(A_H)}\ot A_H)$

\item[ ]

\item[ ]$=\mu_A\co ((\mu_A\co(f\ot i_A))\ot (\mu_A\co (\psi_A\ot f^{-1})\co (H\ot c_{H,A_H})\co (\delta_H\ot A_H))\co (H\ot c_{H,A_H}\ot A_H)$

\item[ ]

\item[ ]$\hspace{0.38cm}\co (\delta_H\ot i_{Z(A_H)}\ot A_H)$

\item[ ]

\item[ ]$=\mu_A\co (\mu_A\ot f^{-1})\co ((\mu_A\co(f\ot (i_A\co i_{Z(A_H)})))\ot (i_A\co \psi_{A_H})\ot H)\co (H\ot Z(A_H)\ot H\ot c_{H, A_H})$

\item[ ]

\item[ ]$\hspace{0.38cm}\co (H\ot Z(A_H)\ot \delta_H\ot A_H)\co (H\ot c_{H, Z(A_H)}\ot A_H)\co (\delta_H\ot Z(A_H)\ot A_H)$

\item[ ]

\item[ ]$=\mu_A\co (\mu_{A}\ot A)\co (f\ot (\mu_{A}\co c_{A,A}\co ((i_{A}\co i_{Z(A_H)})\ot (i_{A}\co \psi_{A_{H}})))\ot f^{-1})\co (H\ot Z(A_H)\ot H\ot c_{H, A_H})$

\item[ ]

\item[ ]$\hspace{0.38cm}\co (H\ot Z(A_H)\ot \delta_H\ot A_H)\co (H\ot c_{H, Z(A_H)}\ot A_H)\co (\delta_H\ot Z(A_H)\ot A_H)$

\item[ ]

\item[ ]$=\mu_A\co ((\mu_A\co ((\mu_A\co (f\wedge f^{-1}\ot i_A))\ot f)\co (H\ot c_{H,A_H})\co (\delta_H\ot A_H))\ot (\mu_A\co ((i_A\co i_{Z(A_H)})\ot f^{-1}))$

\item[ ]

\item[ ]$\hspace{0.38cm}\co (H\ot c_{Z(A_H),A_H}\ot H)\co (H\ot Z(A_H)\ot c_{H,A_H})\co (H\ot c_{H,Z(A_H)}\ot A_H)\co (\delta_H\ot Z(A_H)\ot A_H)$

\item[ ]

\item[ ]$=\mu_A\co ((\mu_A\co c_{A,A}\co (f\ot i_A))\ot (\mu_A\co ((i_A\co i_{Z(A_H)})\ot f^{-1}))\co (H\ot c_{Z(A_H),A_H}\ot H)$

\item[ ]

\item[ ]$\hspace{0.38cm}\co (H\ot Z(A_H)\ot c_{H,A_H})\co (H\ot c_{H,Z(A_H)}\ot A_H)\co (\delta_H\ot Z(A_H)\ot A_H)$

\item[ ]

\item[ ]$=i_A\co \mu_{A_H}\co c_{{A_H},{A_H}}\co ((\varphi_{A_H}\co (H\ot i_{Z(A_H)}))\ot A_H)$

\end{itemize}

In the foregoing computations, the first equality uses the definition of $\varphi_A$; the second one relies on (\ref{primeraigualdad}); the third and fifth ones are consequence of the definition of $\psi_A$, the fourth follows by the properties of the center of $A_H$; finally in the last one we apply (\ref{segundaigualdad})

Then there exists a morphism $\varphi_{Z(A_H)}:H\ot Z(A_H)\rightarrow Z(A_H)$ such that $i_{Z(A_H)}\co \varphi_{Z(A_H)}=\varphi_{A_H}\co (H\ot i_{Z(A_H)})$.
Trivially, $\varphi_{Z(A_H)}$ satisfies the conditions of Definition \ref{weak-H-mod}. In order to show that $(Z(A_H), \varphi_{Z_(A_H)})$ is a left $H$-module algebra, we only need to prove the equality

\begin{equation}
\label{centroesHmodulo}
\varphi_{Z(A_H)}\co (H\ot \varphi_{Z(A_H)})=\varphi_{Z(A_H)}\co (\mu_H\ot Z(A_H))
\end{equation}

which follows composing with the equalizer $i_{Z(A_H)}$. Indeed:

\begin{itemize}

\item[ ]

\item[ ]$\hspace{0.38cm} i_{Z(A_H)}\co \varphi_{Z(A_H)}\co (H\ot \varphi_{Z(A_H)})$

\item[ ]

\item[ ]$=\varphi_{A_H}\co (H\ot \varphi_{A_H})\co (H\ot H\ot i_{Z(A_H)})$

\item[ ]

\item[ ]$=\mu_{A_H}\co (A_H\ot (\mu_{A_H}\co ((i_{Z(A_H)}\co \varphi_{Z(A_H)})\ot A_H)))\co (A_H\ot H\ot c_{A_H,Z(A_H}))\co (\sigma\ot \mu_H\ot \sigma^{-1}\ot Z(A_H))$

\item[ ]

\item[ ]$\hspace{0.38cm}\co (H\ot H\ot \delta_{H\ot H}\ot Z(A_H))\co (\delta_{H\ot H}\ot Z(A_H))$

\item[ ]

\item[ ]$=\mu_{A_H}\co (\sigma\wedge \sigma^{-1}\ot (i_{Z(A_H)}\co \varphi_{Z(A_H)}))\co (H\ot H\ot \mu_H\ot Z(A_H))\co (\delta_{H\ot H}\ot Z(A_H))$

\item[ ]

\item[ ]$=\mu_{A_H}\co (u_2\ot (i_{Z(A_H)}\co \varphi_{Z(A_H)}))\co (H\ot H\ot \mu_H\ot Z(A_H))\co (\delta_{H\ot H}\ot Z(A_H))$

\item[ ]

\item[ ]$=\mu_{A_H}\co (\varphi_{A_H}\ot \varphi_{A_H})\co (H\ot c_{H,A_H}\ot A_H)\co ((\delta_H\co \mu_H)\ot \eta_{A_H}\ot i_{Z(A_H)})$

\item[ ]

\item[ ]$=\varphi_{A_H}\co (\mu_H\ot i_{Z(A_H)})$

\end{itemize}

In the above calculations, the first equality uses that $\varphi_{A_H}$ factors through $i_{Z(A_H)}$; the second follows by (g1); the third one by the properties of the center of $A_H$; the fourth because $\sigma$ is a morphism in $Reg_{\varphi_{A_H}}(H,A_H)$; in the fifth equality we apply that $H$ is a weak Hopf algebra; finally, the last one follows because $(A_H, \varphi_{A_H})$ is a weak $H$-module.

Taking into account that $i_{Z(A_H)}$ is a monomorphism, $\varphi_{Z(A_H)}\co (H\ot \varphi_{Z(A_H)})=\varphi_{Z(A_H)}\co (\mu_H\ot Z(A_H))$ and we conclude the proof.

\end{dem}

The following technical Lemma will be useful for the last Theorem of this paper.

\begin{lem}
\label{igualdadutil}
Let $H$ be a cocommutative weak Hopf algebra and let $A_{H} \hookrightarrow A$ be an $H$-cleft extension. We denote by $\varphi_{A_H}$ the weak $H$-module algebra structure defined for $A_H$ by a convolution invertible total integral $f:H\rightarrow A$, and by  $\psi_{A_H}$ the morphism defined in Proposition \ref{centroesHmoduloalgebra}. Then the equality

\begin{equation}
\label{varphipsi}
\varphi_{A_H}\co (H\ot \psi_{A_H})\co (\delta_H\ot A_H)=\mu_{A_H}\co ((\varphi_{A_H}\co (H\ot \eta_{A_H}))\ot A_H)
\end{equation}

holds.

\end{lem}

\begin{dem}
We compose the left part of the equality with the monomorphism $i_A$. Using that $H$ is cocommutative and (\ref{primeraigualdad}),

\begin{itemize}

\item[ ]

\item[ ]$\hspace{0.38cm} i_A\co \varphi_{A_H}\co (H\ot \psi_{A_H})\co (\delta_H\ot A_H)$

\item[ ]

\item[ ]$=\mu_A\co (\mu_A\ot f\wedge f^{-1})\co (f\wedge f^{-1}\ot c_{H,A})\co (\delta_H\ot i_A)$

\item[ ]

\item[ ]$=\mu_A\co (A\ot \mu_A)\co (f\ot (\mu_A\co (f^{-1}\ot i_A))\ot f\wedge f^{-1})\co (\delta_H\ot c_{H,A})\co (\delta_H\ot A_H)$

\item[ ]

\item[ ]$=\mu_A\co (f\wedge f^{-1}\ot i_A)$

\item[ ]

\item[ ]$=\mu_A\co (u_1\ot i_A)$

\item[ ]

\item[ ]$=\mu_A\co ((i_A\co \varphi_{A_H}\co (H\ot \eta_{A_H}))\ot i_A)$

\item[ ]

\item[ ]$=i_A\co \mu_{A_H}\co ((\varphi_{A_H}\co (H\ot \eta_{A_H}))\ot A_H)$

\end{itemize}

\end{dem}

Now we will show the main result of this section.

\begin{teo}
\label{biyeccionconcohomologia}
Let $H$ be a cocommutative weak Hopf algebra and let $A_{H} \hookrightarrow A$ be an $H$-cleft extension. We denote by $(\varphi_{A_H}, \sigma_{A_H})$ the corresponding crossed system defined by the convolution invertible integral $f:H\rightarrow A$. Then there is a bijective correspondence between the second cohomology group $H_{\varphi_{\mathcal{Z}(A_{H})}}^2(H, Z(A_H))$ and the equivalence classes of crossed systems for $H$ over $A$ having $\varphi_{A_H}$ as weak $H$-module algebra structure.
\end{teo}

\begin{dem}

Let $[\tau]$ be in $H_{\varphi_{A_H}}^2(H, Z(A_H))$. Using the properties of the center of $A_H$, it is not difficult to prove that the morphism $\sigma_{A_H}\wedge (i_{Z(A_H)}\co \tau)$ satisfies conditions (g1) and (g2). As far as (g3), note that

\begin{itemize}

\item[ ]

\item[ ]$\hspace{0.38cm} (\sigma_{A_H}\wedge (i_{Z(A_H)}\co \tau))\co (\Pi_{H}^{L}\ot H)\co \delta_H$

\item[ ]

\item[ ]$=\mu_{A_H}\co (\sigma_{A_H}\ot (i_{Z(A_H)}\co \tau))\co (H\ot c_{H,H}\ot H)\co ((\delta_H\co \Pi_{H}^{L})\ot \delta_H)\co \delta_H$

\item[ ]

\item[ ]$=\mu_{A_H}\co (\sigma_{A_H}\ot (i_{Z(A_H)}\co \tau))\co (H\ot c_{H,H}\ot H)\co (((\Pi_{H}^{L}\ot \Pi_{H}^{L})\co \delta_H)\ot \delta_H)\co \delta_H$

\item[ ]

\item[ ]$=\mu_{A_H}\co ((\sigma_{A_H}\co (\Pi_{H}^{L}\ot H)\co \delta_H)\ot ((i_{Z(A_H)}\co \tau)\co (\Pi_{H}^{L}\ot H)\co \delta_H)) \co \delta_{H}$

\item[ ]

\item[ ]$=u_1\wedge u_1$

\item[ ]

\item[ ]$=u_1$,

\end{itemize}

where the first equality follows by the definition of the convolution product, in the second one we apply that $H$ is cocommutative and therefore $\delta_H\co \Pi_{H}^{L}=(\Pi_{H}^{L}\ot \Pi_{H}^{L})\co \delta_H$; the third uses that $H$ is cocommutative; the fourth relies on (\ref{sigma-PiL-l}); finally, the last one follows because $(A, \varphi_A)$ is a weak left $H$-module algebra.
By the equivalence between (\ref{sigma-eta-l}) and (\ref{sigma-PiL-l}) we have that $u_{1}=(\sigma_{A_H}\wedge (i_{Z(A_H)}\co \tau))\co (\eta_H\ot H)$.
In a similar way  $u_{1}=(\sigma_{A_H}\wedge (i_{Z(A_H)}\co \tau))\co (H\ot \eta_H)$ and then $(\varphi_{A_H}, \sigma_{A_H}\wedge (i_{Z(A_H)}\co \tau))$ is a crossed system for $A$ over $H$.

Conversely let $(\varphi_{A_H}, \gamma)$ be a crossed system for $H$ over $A_{H}$. Then the morphism $\sigma_{A_H}^{-1}\wedge \gamma$ factors through the equalizer $i_{Z(A_H)}$. Indeed:

\begin{itemize}

\item[ ]

\item[ ]$\hspace{0.38cm} \mu_{A_H}\co(A_H\ot \sigma_{A_H}^{-1}\wedge \gamma)$

\item[ ]

\item[ ]$=\mu_{A_H}\co(A_H\ot (u_2\wedge u_2\wedge \sigma_{A_H}^{-1}\wedge \gamma))$

\item[ ]

\item[ ]$=\mu_{A_H}\co (\mu_{A_H}\ot A_H )\co (u_2\wedge u_2\ot A_H\ot \sigma_{A_H}^{-1}\wedge \gamma)\co (H\ot c_{A_H,H}\ot H\ot H)\co (c_{A_H,H}\ot H\ot H\ot H)$

\item[ ]

\item[ ]$\hspace{0.38cm}\co (A_H\ot \delta_{H\ot H}) $

\item[ ]

\item[ ]$=\mu_{A_H}\co (\mu_{A_H}\ot H)\co (\sigma_{A_H}^{-1}\wedge \sigma_{A_H}\wedge u_2\ot A_H\ot \sigma_{A_H}^{-1}\wedge \gamma)\co (H\ot c_{A_H,H}\ot H\ot H)\co (c_{A_H,H}\ot H\ot H\ot H)$

\item[ ]

\item[ ]$\hspace{0.38cm}\co (A_H\ot \delta_{H\ot H}) $

\item[ ]

\item[ ]$=\mu_{A_H}\co (\mu_{A_H}\ot H)\co (\sigma_{A_H}^{-1}\wedge \sigma_{A_H}\ot (
\varphi_{A_H}\co (H\ot \psi_{A_H})\co ((\delta_H \co \mu_H)\ot A_H))\ot \sigma_{A_H}^{-1}\wedge \gamma)$

\item[ ]

\item[ ]$\hspace{0.38cm}\co (\delta_{H\ot H}\ot A_H\ot H\ot H) \co (H\ot c_{A_H,H}\ot H\ot H) \co (c_{A_H,H}\ot H\ot H\ot H)\co (A_H\ot \delta_{H\ot H})$

\item[ ]

\item[ ]$=\mu_{A_H}\co (\mu_{A_H}\ot A_H)\co ((A_H\ot (\mu_{A_H}\co (\sigma_{A_H} \ot (\varphi_{A_H}\co (\mu_H\ot A_H))\co (\delta_{H\ot H}\ot A_H))\ot A_H)$

\item[ ]

\item[ ]$\hspace{0.38cm}\co (\sigma_{A_H}^{-1}\ot H\ot H\ot (\varphi_{A_{H}}\co (\mu_H\ot A_H))\ot \sigma^{-1}\wedge \gamma)\co (H\ot H\ot \delta_{H\ot H}\ot A_H\ot H\ot H)$

\item[ ]

\item[ ]$\hspace{0.38cm}\co (\delta_{H\ot H}\ot A_H\ot H\ot H)\co (H\ot c_{A_H,H}\ot H\ot H)\co (c_{A_H,H}\ot H\ot H\ot H)\co (A_H\ot \delta_{H\ot H})$

\item[ ]

\item[ ]$=\mu_{A_H}\co (\mu_{A_H}\ot A_H)\co ((A_H\ot (\mu_{A_H}\co ((\varphi_{A_H}\co (H\ot \varphi_{A_H}))\ot A_H)\co (H\ot H\ot c_{A_H,A_H})$

\item[ ]

\item[ ]$\hspace{0.38cm}\co (H\ot H\ot \sigma_{A_H}\ot A_H)\co(\delta_{H\ot H}\ot A_H))\ot A_H)\co (\sigma_{A_H}^{-1}\ot H\ot H\ot (\varphi_{A_{H}}\co (\mu_H\ot A_H))\ot \sigma_{A_H}^{-1}\wedge \gamma)$

\item[ ]

\item[ ]$\hspace{0.38cm}\co (H\ot H\ot \delta_{H\ot H}\ot A_H\ot H\ot H)\co (\delta_{H\ot H}\ot A_H\ot H\ot H)\co (H\ot c_{A_H,H}\ot H\ot H)$

\item[ ]

\item[ ]$\hspace{0.38cm}\co (c_{A_H,H}\ot H\ot H\ot H)\co (A_H\ot \delta_{H\ot H})$

\item[ ]

\item[ ]$=\mu_{A_H}\co (\mu_{A_H}\ot A_H)\co (A_H\ot (\varphi_{A_H}\co (H\ot \varphi_{A_H}))\ot A_H)$

\item[ ]

\item[ ]$\hspace{0.38cm}\co (\sigma_{A_H}^{-1}\ot H\ot H\ot (\varphi_{A_{H}}\co (\mu_H\ot A_H))\ot \sigma_{A_H}\wedge \sigma_{A_H}^{-1}\wedge \gamma)\co (H\ot H\ot \delta_{H\ot H}\ot A_H\ot H\ot H)$

\item[ ]

\item[ ]$\hspace{0.38cm}\co (\delta_{H\ot H}\ot A_H\ot H\ot H)\co (H\ot c_{A_H,H}\ot H\ot H)\co (c_{A_H,H}\ot H\ot H\ot H)\co (A_H\ot \delta_{H\ot H})$

\item[ ]

\item[ ]$=\mu_{A_H}\co (A_H\ot (\mu_{A_H}\co ((\varphi_{A_H}\co (H\ot \varphi_{A_H}))\ot A_H)\co (H\ot H\ot c_{A_H,A_H})\co (H\ot H\ot \gamma\ot A_H)$

\item[ ]

\item[ ]$\hspace{0.38cm}\co(\delta_{H\ot H}\ot A_H)))\co (A_H\ot H\ot H\ot (\varphi_{A_{H}}\co (\mu_H\ot A_H)))\co (\sigma_{A_H}^{-1}\ot \delta_{H\ot H}\ot A_H)$

\item[ ]

\item[ ]$\hspace{0.38cm}\co (H\ot H\ot H\ot c_{A_H,H})\co (H\ot H\ot c_{A_H,H}\ot H)\co (H\ot c_{A_H,H}\ot H\ot H)\co (c_{A_H,H}\ot H\ot H\ot H)$

\item[ ]

\item[ ]$\hspace{0.38cm}\co (A_H\ot \delta_{H\ot H})$

\item[ ]

\item[ ]$=\mu_{A_H}\co (A_H\ot (\mu_{A_H}\co (\gamma \ot (\varphi_{A_H}\co (\mu_{H}\ot A_{H})))))\co(\delta_{H\ot H}\ot A_H))\co (A_H\ot H\ot H\ot (\varphi_{A_{H}}\co (\mu_H\ot A_H)))$

\item[ ]

\item[ ]$\hspace{0.38cm}\co (\sigma_{A_H}^{-1}\ot \delta_{H\ot H}\ot A_H)\co (H\ot H\ot H\ot c_{A_H,H})\co (H\ot H\ot c_{A_H,H}\ot H)\co (H\ot c_{A_H,H}\ot H\ot H)$

\item[ ]

\item[ ]$\hspace{0.38cm}\co (c_{A_H,H}\ot H\ot H\ot H)\co (A_H\ot \delta_{H\ot H}) $

\item[ ]

\item[ ]$=\mu_{A_H}\co (\sigma_{A_H}^{-1}\wedge \gamma\ot (\varphi_{A_H}\co (H\ot \varphi_{A_{H}})\co (\delta_H \co \mu_H)\ot A_H))\ot (\delta_{H\ot H}\ot A_H)\co (H\ot c_{A_H,H})\co (c_{A_H,H}\ot H)$

\item[ ]

\item[ ]$=\mu_{A_H}\co (\sigma_{A_H}^{-1}\wedge \gamma\wedge u_2\ot A_H)\co (H\ot c_{A_H,H})\co (c_{A_H,H}\ot H)$

\item[ ]

\item[ ]$=\mu_{A_H}\co (\sigma_{A_H}^{-1}\wedge \gamma\ot A_H)\co (H\ot c_{A_H,H})\co (c_{A_H,H}\ot H)$

\end{itemize}

In the above computations, the first and the third equalities follow because $\sigma_{A_H}^{-1}$ is in $Reg_{\varphi_{A_H}}(H,A_H)$; the second one because $u_2$ factors through the center of $A_H$; in the fourth and the eleventh ones we use (\ref{centroesHmodulo}); in the fifth and the tenth equalities we apply that $H$ is a weak Hopf algebra; the sixth and ninth rely on (g1) for $\sigma_{A_H}$ and $\gamma$, respectively; the seventh one follows by cocommutativity; the eighth uses that $H$ is cocommutative and $\sigma_{A_H}$ is a morphism in $Reg_{\varphi_{A_H}}(H,A_H)$; finally, the twelfth equality follows by (f1) and (f2) for $\gamma$.

As a consequence, using that $H$ is cocommutative, $\sigma_{A_H}^{-1}\wedge \gamma\wedge \sigma_{A_H}=\sigma_{A_H}\wedge\sigma_{A_H}^{-1}\wedge \gamma=\gamma$, and therefore $\sigma_{A_H}^{-1}\wedge \gamma=\gamma\wedge \sigma_{A_H}^{-1}$.

The proof for the condition (g2) follows a similar pattern to the one developed in \cite{doi1} and will be omitted. As far as (g3) the proof follows in a similar way to the one giving for $\sigma_{A_H}\wedge (i_{Z(A_H)}\co \tau)$ using Proposition \ref{sigma-sigma-1}.

Finally, we have to show that the correspondence is well defined.
Let $[\tau]$ and $[\tau^{\prime}]$ be in $H^2(H, Z(A_H))$ such that
the crossed systems $(\varphi_{A_H}, \sigma_{A_H}\wedge
(i_{Z(A_H)}\co \tau))$ and $(\varphi_{A_H}, \sigma\wedge
(i_{Z(A_H)}\co \tau^{\prime}))$ are equivalent. Let $h$ be the
morphism in $Reg_{\varphi_{A_H}}(H,A_H)$ satisfying conditions
(\ref{relacionsigmas}) and (\ref{relacionfis}). Then $h$ factors
through the center of $A_H$. Indeed:

\begin{itemize}

\item[ ]

\item[ ]$\hspace{0.38cm} \mu_{A_H}\co(h\ot A_H)$

\item[ ]

\item[ ]$=\mu_{A_H}\co((h\wedge u_1\wedge u_1)\ot A_H)$

\item[ ]

\item[ ]$=\mu_{A_H}\co((h\wedge u_1)\ot (\mu_{A_H}\co c_{A_H,A_H}))\co (H\ot u_1\ot A_H)\co (\delta_H\ot A_H)$

\item[ ]

\item[ ]$=\mu_{A_H}\co (\mu_{A_H}\ot u_1)\co (h\ot (\mu_{A_H}\co ((\varphi_{A_H}\co (H\ot \eta_{A_H}))\ot A_H)))\ot H)\co (\delta_H\ot c_{H,A_H})\co (\delta_H\ot A_H)$

\item[ ]

\item[ ]$=\mu_{A_H}\co (\mu_{A_H}\ot u_1)\co (h\ot (\varphi_{A_H}\co (H\ot \varphi_{A_{H}})\co (\delta_H\ot A_H))\ot H)\co (\delta_H\ot c_{H,A_H})\co (\delta_H\ot A_H)$

\item[ ]

\item[ ]$=\mu_{A_H}\co ((\mu_{A_H}\co (\mu_{A_H}\ot A_H)\co (h\ot \varphi_{A_H}\ot h^{-1})\co (\delta_H\ot c_{H, A_H})\co (\delta_H\ot A_H))\ot h)$

\item[ ]

\item[ ]$\hspace{0.38cm}\co (H\ot \varphi_{A_{H}}\ot H)\co (\delta_H\ot c_{H,A_H})\co (\delta_H\ot A_H)$

\item[ ]

\item[ ]$=\mu_{A_H}\co (\varphi_{A_H}\ot h)\co (H\ot \varphi_{A_{H}}\ot H)\co (\delta_H\ot c_{H,A_H})\co (\delta_H\ot A_H)$

\item[ ]

\item[ ]$=\mu_{A_H}\co ((\mu_{A_H}\co (u_1\ot A_H))\ot h)\co (H\ot c_{H,A_H})\co (\delta_H\ot A_H)$

\item[ ]

\item[ ]$=\mu_{A_H}\co (A_H\ot (u_1\wedge h))\co c_{H,A_H}$

\item[ ]

\item[ ]$=\mu_{A_H}\co (A_H\ot h)\co c_{H,A_H}.$

\end{itemize}

In the foregoing computations, the first and the last equalities follow because $h$ is in $Reg_{\varphi_{A_H}}(H,A_H)$; the second, third and eight ones use the definition of $u_1$ and that this morphism factors through the center of $A_H$; the fourth and seventh equalities rely on (\ref{varphipsi}); the fifth is a consequence of cocommutativity of $H$; finally, the sixth one follows by (\ref{relacionfis}).

Using that $h$ factors through the center of $A_H$ and by (\ref{relacionsigmas}) it is not difficult to see that $\tau$ and $\tau^{\prime}$ are cohomologous.

Conversely, if $\tau$ and $\tau^{\prime}$ are cohomologous, by the properties of the center of $A_H$ we get that the corresponding crossed systems $(\varphi_{A_H}, \sigma_{A_H}\wedge (i_{Z(A_H)}\co \tau))$ and $(\varphi_{A_H}, \sigma_{A_H}\wedge (i_{Z(A_H)}\co \tau^{\prime}))$ are equivalent and we conclude the proof.

\end{dem}

\section*{Acknowledgements}
The authors were supported by Ministerio de Ciencia e
Innovaci\'on
(Project: MTM2010-15634) and by FEDER.

\end{document}